\documentclass[onecolumn,final,a4paper]{elsarticle}
\graphicspath{{Figures/}}
\usepackage{comment}
\usepackage{placeins}

\usepackage[utf8]{inputenc}
\usepackage{amsmath,amssymb,amsfonts,latexsym,stmaryrd}
\usepackage{xcolor}
\SetSymbolFont{stmry}{bold}{U}{stmry}{m}{n}
\usepackage{algorithmic}
\usepackage{booktabs}

\usepackage[colorlinks=true,breaklinks=true,linkcolor=lightblue,citecolor=lightblue,hypertexnames=false]{hyperref}
\usepackage{stackengine}
\usepackage{mathrsfs}
\usepackage{rotating} 
\usepackage{amssymb,amsmath,bbm,amsfonts,latexsym}
\usepackage{amsthm}
\usepackage{graphicx,float,epsfig,color,fancyhdr}

\newif\ifauthorcolors
\authorcolorsfalse
\definecolor{colDib}{rgb}{0,0,1}         
\definecolor{colRado}{rgb}{0.85,0,0}     
\definecolor{colDuc}{rgb}{0,0.55,0}      
\definecolor{colMario}{rgb}{1,0.5,0}     
\definecolor{colOleg}{rgb}{0.5,0,0.6}    
\definecolor{colTed}{rgb}{0,0.5,0.5}     
\definecolor{colBoian}{rgb}{0.55,0.3,0}  

\usepackage{placeins}
\usepackage{amsopn}

\usepackage[utf8]{inputenc}
\usepackage{graphicx}
\usepackage{wrapfig}

\usepackage{amsfonts}
\usepackage{amsmath}
\numberwithin{equation}{section}
\usepackage{mathtools}

\usepackage[ruled,vlined,linesnumbered]{algorithm2e}
\let\oldnl\nl
\newcommand{\nonl}{\renewcommand{\nl}{\let\nl\oldnl}}%

\usepackage[utf8]{inputenc}
\usepackage{graphicx}
\usepackage{wrapfig}

\usepackage{amsfonts}
\usepackage{amsmath}
\numberwithin{equation}{section}
\usepackage{mathtools}

\usepackage{makecell}

\usepackage[ruled,vlined,linesnumbered]{algorithm2e}
\let\oldnl\nl

\renewcommand{\vec}[1]{\mathbf{#1}}

\newcommand{\HAT}[1]{\widetilde{}}

\usepackage{lipsum}

\newtheorem{lemma}{Lemma}[section]

\newtheorem{proposition}{Proposition}[section]
\theoremstyle{definition}
\newtheorem{definition}{Definition}[section]
\theoremstyle{plain}

\theoremstyle{remark}

\theoremstyle{plain}

\def\bm{\mathbf{m}}

\newcommand\bn{\mathbf{n}}

\def\l{\lambda}

\newcommand{\vertiii}[1]{{\left\vert\kern-0.25ex\left\vert\kern-0.25ex\left\vert #1 
    \right\vert\kern-0.25ex\right\vert\kern-0.25ex\right\vert}}

\newcommand*{\jump}[1]{\lbrack\hspace{-1.5pt}\lbrack #1\rbrack\hspace{-1.5pt}\rbrack}

\usepackage[utf8]{inputenc}
\usepackage[T1]{fontenc}
\usepackage{amsmath,amssymb,amsthm}
\usepackage{bm}
\usepackage{graphicx}
\usepackage{geometry}
\usepackage{hyperref}
\usepackage{booktabs}
\usepackage{siunitx}
\usepackage{pgfplots}
\usetikzlibrary{decorations.pathreplacing}
\pgfplotsset{compat=1.18}

\hypersetup{
    pdftitle={A hybrid TT/QTT discontinuous Galerkin method with fixed-matrix upwinding for three-dimensional transient multigroup neutron transport},
    pdfauthor={D. Adak, R. Vuchkov, O. Korobkin, T. Sariyski, B. S. Alexandrov},
    colorlinks=true,
    linkcolor=blue,
    citecolor=blue,
    urlcolor=blue
}

\newcommand{\rvec}{\vec{r}}
\newcommand{\Om}{\hat{\Omega}}
\newcommand{\dOm}{d\hat{\Omega}}
\newcommand{\pd}[2]{\frac{\partial #1}{\partial #2}}

\usepackage[normalem]{ulem} 

\usepackage{multirow}

\begin{document}

\begin{frontmatter}

  \title{A hybrid TT/QTT discontinuous Galerkin method with fixed-matrix
         upwinding for three-dimensional transient multigroup neutron transport}

  \author[IIT-KGP]{Dibyendu Adak\fnref{equal}}
  \author[SNL]{Radoslav Vuchkov\fnref{equal}}
  \fntext[equal]{These authors contributed equally to this work.}
  \author[SLIC]{Oleg Korobkin}
  \author[SLIC]{Ted Sariyski}
  \author[Tensor,UNM]{Boian S. Alexandrov\corref{cor1}}
  \ead{boian@tensormedium.ai}

  \cortext[cor1]{Corresponding author.}

  \address[IIT-KGP]{Indian Institute of Technology--Kharagpur, WB, India}
  \address[SNL]{Sandia National Laboratories, Albuquerque, NM 87185, USA}
  \address[SLIC]{SLIC.AI LLC, Santa Fe, NM, USA}
  \address[Tensor]{Tensor Medium Inc, Fort Worth, TX, USA}
  \address[UNM]{University of New Mexico, Albuquerque, NM, USA}

  \begin{abstract}
The time-dependent Boltzmann equation describes the angular neutron flux
over position, direction and energy as it evolves in time. Full
tensor-product discretizations produce prohibitively large solution
arrays. We develop and analyze an upwind discontinuous Galerkin method
in a hybrid tensor-train/quantized tensor-train (TT/QTT) format. Our
central contribution is a fixed-matrix upwind construction for this
format: the same pair of one-dimensional DG matrices serves all
directions and energy groups, with signed and absolute direction-cosine
factors carrying the angular dependence and preserving the Kronecker
structure. For homogeneous materials on uniform Cartesian
grids, we establish operator rank bounds independent of the mesh size,
number of energy groups and number of directions at fixed scattering
order. We derive a conditional final-time error bound in terms of
accumulated step defects and enforce group-integrated balance for the
compressed system. Crank--Nicolson experiments confirm second-order
convergence on a separable manufactured solution whose finest final
state occupies $13$\,kB in the hybrid format, compared with $35$\,TB as a
full tensor---a compression factor of approximately $2.7\times10^9$.
A modified rank-adaptive basis-update and Galerkin integrator (TT-BUG)
advances passive and active neutron bursts with approximately half a
billion unknowns. At comparable verified scalar-flux errors, selected
TT-BUG runs are approximately $15$--$20$ times faster than the tested
AMEn baseline. For the reported anisotropic active-burst configurations,
TT-BUG computes the angular flux on one processor approximately $50$
times faster than OpenSn running on $24$ processors, corresponding to a
$1200$-fold reduction in processor-time cost, defined as processors times
wall time. OpenSn supplies an independent discretization check; the
processor-cost ratio compares the reported configurations rather than
matched physical errors. The modified TT-BUG source-onset and saturation
treatment is validated numerically; its order, convergence and general
tolerance-to-error relation remain unproved.
  \end{abstract}

\end{frontmatter}

\section{Introduction}\label{sec:intro}

Neutrons in a reactor core, a radiation shield or a critical assembly
are far fewer than the nuclei around them, and neutron--neutron
collisions are neglected at the densities considered here. Their population is therefore
governed by a linear kinetic equation, the linear Boltzmann equation
for neutrons, which states the conservation of neutrons in the phase
space of position, direction of flight and
energy~\cite{bell1970nuclear,lewis1993computational}. Two families of methods solve the Boltzmann equation. Monte Carlo
methods~\cite{romano2015openmc} follow the histories of sampled
neutrons with the cross sections as probabilities, need no mesh in
phase space and treat any geometry, but their error decreases only
with the square root of the number of histories, so that a
time-dependent flux resolved in energy throughout a volume is
expensive to obtain from them. Deterministic methods discretize every
coordinate: the direction by discrete ordinates or by spherical
harmonics~\cite{lewis1993computational}, the energy by groups, space
by finite differences, characteristics or finite elements, among them
the upwind discontinuous Galerkin scheme introduced for this equation
by Reed and Hill~\cite{reed1973triangular,lesaint1974finite}, and time
by implicit stepping. Such methods typically combine upwinding in space
with Krylov solvers for the linear systems of each time
step~\cite{adams2002fast}, and production codes such as
PARTISN~\cite{alcouffe2017partisn} and OpenSn~\cite{andrs2026opensn}
parallelize them over directions, groups and subdomains. The difficulty of the deterministic approach is the number of unknowns,
which is the product of the numbers of degrees of freedom in position,
direction and energy. Let $N$ be
the number of cells per axis of a cube, $M$ the number of discrete
directions and $G$ the number of energy groups. A scheme of degree one
has two coefficients per cell along each axis, so there are $n = 2N$
coefficients per axis, and the discrete flux at one time step has
$n^{3} M G$ entries, one for each choice of the three position
coefficients, the direction and the group. Realistic problems, such as
the pulsed slowing down of fast neutrons, can require hundreds of energy
groups and more than a hundred directions. The manufactured-solution
verification in Section~\ref{sec:manufactured} uses $G=256$ and $M=128$:
at $N=32$ the discrete flux has $8.6\times10^9$ entries, and at $N=256$
the finest final state occupies $13$\,kB in hybrid TT/QTT format instead
of $35$\,TB as a full tensor, a compression factor of approximately
$2.7\times10^9$ (Table~\ref{tab:mms}). This example tests a deliberately
separable solution. The burst experiments use
$G=16$ and $M=128$ and assess the cost of evolving non-separable states.
These tests do not establish runtime scaling with $G$ or $M$.

Low-rank and reduced order models (ROMs)~\cite{benner2015survey} are a popular alternative, in which one seeks to reduce the storage and computational effort by exploiting low rank structure. The dynamical
low-rank approximation has been applied to time-dependent radiation
transport in~\cite{peng2020lowrank,einkemmer2021asymptotic},
see~\cite{einkemmer2025review} for a review. Rank-adaptive basis-update
and Galerkin (BUG)
integrators~\cite{ceruti2022unconventional,ceruti2022rankadaptive}
provide one approach, and recent work adds exact
conservation~\cite{coughlin2024conservative} and multilevel uncertainty
quantification~\cite{patwardhan2026mlmc}. These integrators evolve a
factorization that separates position from direction or velocity, and
they are usually combined with explicit substeps under a CFL
condition~\cite{patwardhan2026mlmc,coughlin2024conservative}. Here the
tensor train separates position, direction and energy group and, in
quantized form, the binary digits of each space index. We use this same
representation with a modified rank-adaptive TT-BUG integrator with
implicit projected substeps as the principal time-stepping method in
our performance tests, and implicit linear systems solved in the
compressed format as a baseline \cite{coughlin2024conservative,oseledets2011tensor,hackbusch2019tensor,oseledets2010quantized,khoromskij2011quantics}. 
Tensor trains have been applied to a variety of physical problems
modeled by partial differential equations, among them plasma
physics~\cite{kormann2015semilagrangian,ye2024quantized,chinomona2026tensortrain},
the Boltzmann equation of kinetic theory~\cite{huang2026homogeneous},
thermal radiation transport~\cite{gorodetsky2025thermal,deshpande2026multigroup},
space-time finite elements for convection-diffusion-reaction
problems~\cite{adak2025spacetime}, and neutron criticality and steady
transport~\cite{truong2024tensor,ortega2026reactor,myers2026tensorized}.
Our work is most similar to~\cite{truong2024tensor}, whose authors also
use the mixed TT and QTT formulation, limited however to the
steady-state problem. 

To the best of our knowledge, we present the first hybrid TT/QTT
discontinuous Galerkin solver for the three-dimensional, time-dependent
multigroup Boltzmann neutron transport equation, with both the discrete
operators and the angular flux represented in compressed form. 
The central methodological contribution is the fixed-matrix upwind
construction for hybrid TT/QTT transient transport. Two one-dimensional
DG matrices, $\mathbf{C}$ and $\mathbf{D}$, are assembled once and serve
all spatial axes, discrete directions and energy groups on the grids
considered here. The direction dependence enters through the analytic
factors $\mu$ and $|\mu|$ in
$\mathbf{A}(\mu)=\mu\mathbf{C}+|\mu|\mathbf{D}$
(Section~\ref{sec:dg-space}). This construction provides a reusable
Kronecker representation of the upwind DG operator, for which we derive
explicit quantized rank bounds independent of mesh size, group count
and direction count at fixed scattering order. Tensor upwinding and BUG
integration already exist; the contribution here is this compressed DG
operator construction and its analysis, using the classical upwind flux. Our
modified TT-BUG implementation adds source-onset and saturation-triggered
basis enrichment. Its numerical performance is assessed against AMEn,
while its order, convergence and general tolerance-to-error relation
remain unproved. The step-defect analysis quantifies the effect of
defects once they are bounded; it does not supply those missing BUG
guarantees. Full-array runs provide tightly converged discrete references,
and OpenSn supplies an independent discretization check. Comparisons
against the reference-generation costs are distinguished from
comparisons at similar measured errors.

The burst experiments demonstrate a substantial reduction in cost for
the tested TT-BUG implementation. At comparable verified scalar-flux
errors, selected TT-BUG runs are approximately $15$--$20$ times faster
than the tested AMEn baseline. For the reported anisotropic active-burst
configurations, TT-BUG computes the angular flux in $14.3$\,s on one
processor, compared with $716$\,s for OpenSn on $24$ processors
(Table~\ref{tab:active}). This is approximately $50$ times less wall
time and $1200$ times less processor-time cost, where processor-time
cost is the number of allocated processors times wall time. The codes
use different spatial discretizations and discrete references; this
configuration-specific ratio does not establish a speedup at matched
physical or angular-flux error.

The paper is organized as follows. Section~\ref{sec:tt} reviews the
tensor-train format and its quantized form. Section~\ref{sec:problem}
states the multigroup transport problem. Section~\ref{sec:dg}
discretizes the direction by discrete
ordinates~\cite{lewis1993computational}, space by the upwind
discontinuous Galerkin scheme and time by the $\theta$-method, and
describes the two time-stepping methods, implicit solves with the
alternating minimal energy method (AMEn)~\cite{dolgov2014amen} and the
modified TT-BUG integrator. Section~\ref{sec:theory} bounds the ranks
of the operators, bounds the error at the final time by the residuals
of the steps and introduces the group balance projection.
Section~\ref{sec:examples} verifies the solver on a manufactured
solution and, against full-tensor references and the transport code
OpenSn, on passive and active neutron bursts, whose evolving flux
couples position, direction and energy, and compares the cost of TT-BUG
with AMEn, with the full-tensor solver and with OpenSn. Section~\ref{sec:conclusions}
concludes.

\section{Tensor trains and the quantized form}
\label{sec:tt}
\label{sec:qtt}
\label{sec:rounding}

The tensor train (TT)
format~\cite{oseledets2011tensor,hackbusch2019tensor} stores a tensor
$\mathcal{X}(i_1,\dots,i_p)$ with $p$ indices of sizes $n_1,\dots,n_p$
as a product of three-way cores,
\begin{equation}
\mathcal{X}(i_1,\dots,i_p)
= \sum_{\alpha_0,\dots,\alpha_p}
\mathcal{U}_1(\alpha_0,i_1,\alpha_1)\,
\mathcal{U}_2(\alpha_1,i_2,\alpha_2)\cdots
\mathcal{U}_p(\alpha_{p-1},i_p,\alpha_p),
\label{eq:tt}
\end{equation}
where the rank index $\alpha_j$ runs from $1$ to $r_j$ and
$r_0 = r_p = 1$. The numbers $r_1,\dots,r_{p-1}$ are the TT-ranks of the
train, $r_j$ sitting on the link $j$ between the cores
$\mathcal{U}_j$ and $\mathcal{U}_{j+1}$, and the smallest admissible
$r_j$ is the rank of the matrix obtained by unfolding $\mathcal{X}$
with the indices $i_1,\dots,i_j$ as rows and the remaining indices as
columns~\cite{oseledets2011tensor}. The train takes $\sum_j r_{j-1}\,n_j\,r_j$ numbers, linear in
the number of indices and quadratic in the ranks, against
$n_1\cdots n_p$ for the full tensor. A matrix acting on such a tensor is
stored in the same way, as a TT matrix whose cores carry a row and a
column index, and a Kronecker product of $p$ factors, one acting on
each index, is a TT matrix with every TT-rank equal to one. The product
of a TT matrix with a train and the sum of two trains are computed
core by core, with ranks that multiply and add, and rounding brings
the ranks back: a pass of truncated singular value decompositions over
the links keeps the tensor to within a relative tolerance $\varepsilon$
in the Frobenius norm~\cite{oseledets2011tensor}, optionally with a cap
$r_{\max}$ on the ranks. A cap can force additional truncation, in
which case the tolerance bound need not hold; Table~\ref{tab:mms-diag}
records such events. The tolerance is the main parameter of the
method.

In this work the tensors have five indices, the position along each of
the three axes, the direction and the group, in that order: the
coefficient vector $\bm{\Psi}$ of the scheme of Section~\ref{sec:dg} is
read as a tensor of sizes $n$, $n$, $n$, $M$, $G$, with $n = 2N$
coefficients along an axis of $N$ cells. The space indices are
quantized~\cite{oseledets2010quantized,khoromskij2011quantics}: with
$N$ a power of two, $n = 2^{d}$, a position index is written in its
$d$ binary digits,
\begin{equation}
i_x - 1 = \sum_{j=1}^{d} 2^{j-1}\beta_j , \qquad \beta_j \in \{0,1\},
\label{eq:qtt-digits}
\end{equation}
and likewise $i_y$ and $i_z$, and the train carries one core per digit.
The least significant digit selects the mean or the slope within a cell
and the others are the digits of the cell number, which is the
cell-major, basis-minor ordering of Section~\ref{sec:dg}. The direction
and the group index are left whole, and this mixed form, quantized in
space and plain in direction and group, is what we call the hybrid of
tensor train and quantized tensor train (QTT). The implementation
used in Section~\ref{sec:examples} puts the group and the direction
first, in the order group, direction, $x$, $y$, $z$. A single Kronecker
product has rank one across every link between different indices in
either order, but the ranks of a sum of products, and those of the
solution, depend on the order, because each link separates a different
set of indices; \ref{app:rank-order} gives the bound that holds
in every order. The ranks across the
links between two digits of the same axis are the quantized ranks, and
the storage of an axis falls from $n\,r^2$ to $2d\,r^2$, logarithmic
in the number of cells. In this form a Kronecker product of five
factors has rank one across the four links between different indices
and, inside the digits of an axis, the quantized rank of its space
factor, which is at most three for any matrix with constant
coefficients built from the identity and the shift along the
axis~\cite{kazeev2012laplace}.

\paragraph{Dependence on groups and directions}
For the implemented order, group, direction, then the $3d$ spatial
digits, the number of stored entries of a state is exactly
\begin{equation}
S_{\mathrm{state}}
= G r_1 + M r_1r_2 + 2\sum_{k=3}^{3d+2}r_{k-1}r_k,
\qquad r_0=r_{3d+2}=1.
\label{eq:storage-gm}
\end{equation}
A TT matrix with dense cores and corresponding ranks $R_k$ stores
\begin{equation}
S_{\mathrm{operator}}
= G^2 R_1 + M^2 R_1R_2
+ 4\sum_{k=3}^{3d+2}R_{k-1}R_k,
\qquad R_0=R_{3d+2}=1.
\label{eq:storage-gm-op}
\end{equation}
Thus even fixed ranks do not make storage independent of $G$ or $M$;
only the space indices are quantized. These counts describe the stored
cores, including their zero entries, rather than the least possible
storage of the operator. A matrix-free implementation can instead keep
diagonal direction factors as vectors and the isotropic angular factor
as an outer product. Neither count predicts runtime or peak
memory, which also depend on the algorithm, on the intermediate
quantities it holds and on the ranks of its iterates. Comparisons as $G$ and $M$ increase must therefore
measure both accuracy and cost rather than infer a speedup from the
operator-rank bounds.

\section{The transport problem}
\label{sec:problem}

\subsection{The linear Boltzmann equation}

Throughout the paper, the measure $\dOm$ on the unit sphere is
normalized to total mass one, so that the scalar flux is the angular
average of the angular flux, and the scattering kernel and the fission
spectrum are taken per unit of this normalized measure.
The unknown is the angular flux $\psi(\rvec,\Om,E,t)$, and it obeys
\begin{equation}
\begin{aligned}
\underbrace{\frac{1}{v(E)}\,\pd{\psi}{t}}_{\text{change in time}}
+ \underbrace{\Om\cdot\nabla\psi}_{\text{streaming}}
+ \underbrace{\sigma(\rvec,E)\,\psi}_{\text{removal by collision}}
&= \\
&\quad \underbrace{\int_{0}^{\infty}\!\!\int_{\Om'}
   \sigma_s(\rvec,E'\!\to\!E,\Om'\!\cdot\!\Om)\,
   \psi(\rvec,\Om',E',t)\,\dOm'\,dE'}_{\text{gain by scattering}} \\
&\quad + \underbrace{\chi(E)\int_{0}^{\infty}
   \nu\sigma_f(\rvec,E')\,\phi(\rvec,E',t)\,dE'}_{\text{gain by fission}}
 + \underbrace{q(\rvec,\Om,E,t)}_{\text{external source}}
\end{aligned}
\label{eq:transport}
\end{equation}
inside the body $\mathcal{D}$. On the left-hand side, the arguments
$(\rvec,\Om,E,t)$ of $\psi$ are suppressed. The scalar flux is

$\phi(\rvec,E,t) = \int \psi(\rvec,\Om,E,t)\,\dOm$. The position $\rvec$ ranges over $\mathcal{D}\subset\mathbb{R}^3$,
the direction of flight $\Om$ over the unit sphere, the kinetic energy
$E$ over the positive reals and the time $t$ over an interval that starts
at zero, and the three variables $(\rvec,\Om,E)$ together form the phase
space. The angular flux is the number of neutrons per unit volume, per
unit normalized solid angle and per unit energy at $(\rvec,\Om,E)$,
multiplied by their speed $v(E)$: the quotient $\psi/v$ is the neutron
density, and $\psi$ itself is the rate at which neutrons cross a unit
area placed normal to $\Om$. The total macroscopic cross section
$\sigma(\rvec,E)$ is the probability per unit path length that a neutron
of energy $E$ collides with a nucleus, in any way, and its inverse is the
mean free path. The differential scattering cross section
$\sigma_s(\rvec,E'\!\to\!E,\Om'\!\cdot\!\Om)$ is the probability per unit
path that a neutron of energy $E'$ and direction $\Om'$ leaves a collision
with energy $E$ and direction $\Om$, per unit energy and per unit
normalized solid angle of arrival. In an isotropic medium it depends on
the two directions only through their cosine $\Om'\cdot\Om$. The quantity
$\nu\sigma_f(\rvec,E')$ is the number of fission neutrons born per unit
path of a neutron of energy $E'$, and $\chi(E)$ is the energy spectrum of
the neutrons born, with $\int_{0}^{\infty}\chi\,dE = 1$. Fission neutrons
are born isotropically, which is why the fission gain carries no
direction argument and depends on $\psi$ only through the scalar flux.
Finally, $q(\rvec,\Om,E,t)$ is an external source, neutrons put into the
body by something other than a collision. Each term of
\eqref{eq:transport} is a rate of neutrons per unit volume of phase
space. The left side is what happens to the neutrons that are in the
state $(\rvec,\Om,E)$: their density changes in time, they move in
straight lines along $\Om$, and any collision takes them out of the
state. The right side is what puts neutrons into the state: scattering
from other energies and directions, fission, and the external source.
The equation needs the state at the start,
\begin{equation}
\psi(\rvec,\Om,E,0) = \psi_0(\rvec,\Om,E)
\qquad \text{for } \rvec \in \mathcal{D} ,
\label{eq:transport-ic}
\end{equation}
and, on the boundary of $\mathcal{D}$, the flux of the neutrons that
come in,
\begin{equation}
\psi(\rvec,\Om,E,t) = \psi_{\mathrm{in}}(\rvec,\Om,E,t)
\qquad \text{for } \rvec \in \partial\mathcal{D}
\text{ and } \Om\cdot\bn(\rvec) < 0 ,
\label{eq:transport-bc}
\end{equation}
where $\bn(\rvec)$ is the outward unit normal at the boundary point.
Only the incoming directions can be prescribed, and what leaves is
decided by the equation. A vacuum boundary has $\psi_{\mathrm{in}} = 0$:
nothing comes back in. A beam is a prescribed $\psi_{\mathrm{in}}$ on
the part of the boundary through which it enters. In all our examples
$\mathcal{D}$ is the cube $(0,L)^3$ of side $L$.

\subsection{Energy groups and group constants}

The energy axis is split into $G$ intervals $[E_g, E_{g-1}]$,
$g = 1,\dots,G$, with $E_0 > E_1 > \dots > E_G$, numbered from the
highest energy group to the lowest energy group, so that the first group holds the fastest neutrons and the $G$ th
group the slowest. The unknown of group $g$ is the angular flux
integrated over its interval,
\begin{equation}
\psi_g(\rvec,\Om,t) = \int_{E_g}^{E_{g-1}} \psi(\rvec,\Om,E,t)\,dE ,
\qquad g = 1,\dots,G .
\label{eq:group-flux}
\end{equation}
Integrating \eqref{eq:transport} over a group requires the shape of the
flux in energy inside the group, which is replaced by a fixed weight
$W(E)$. This defines the group constants, in which the
dependence on $\rvec$ is not written. The medium is homogeneous in
all our examples, so the group constants do not depend on
position. The total cross section and the
inverse speed of group $g$ are the weighted averages
\begin{equation}
\sigma^{g} =
\frac{\int_{E_g}^{E_{g-1}} \sigma(E)\,W(E)\,dE}
     {\int_{E_g}^{E_{g-1}} W(E)\,dE} ,
\qquad
\frac{1}{v_g} =
\frac{\int_{E_g}^{E_{g-1}} \frac{1}{v(E)}\,W(E)\,dE}
     {\int_{E_g}^{E_{g-1}} W(E)\,dE} .
\label{eq:group-total}
\end{equation}
The scattering transfer cross section from group $g'$ to group $g$ is
\begin{equation}
\sigma_s^{g'\to g}(\Om'\!\cdot\!\Om) =
\frac{\int_{E_g}^{E_{g-1}}\!\int_{E_{g'}}^{E_{g'-1}}
      \sigma_s(E'\!\to\!E,\Om'\!\cdot\!\Om)\,W(E')\,dE'\,dE}
     {\int_{E_{g'}}^{E_{g'-1}} W(E')\,dE'} ,
\label{eq:group-scatter}
\end{equation}
the cross section for a neutron that enters a collision in group $g'$
to leave it in group $g$. The arrow reads from the group before the
collision to the group after it. The fission production cross section
and the fission spectrum of a group are
\begin{equation}
\nu\sigma_f^{g} =
\frac{\int_{E_g}^{E_{g-1}} \nu\sigma_f(E)\,W(E)\,dE}
     {\int_{E_g}^{E_{g-1}} W(E)\,dE} ,
\qquad
\chi_g = \int_{E_g}^{E_{g-1}} \chi(E)\,dE ,
\label{eq:group-fission}
\end{equation}
where $\nu\sigma_f^{g}$ is the number of neutrons produced by fission
per unit path travelled by a neutron of group $g$, and $\chi_g$ is the
fraction of the fission neutrons born in group $g$. 
The number of groups and weights $W(E)$ are prescribed inputs, and we
study the numerical solution of the resulting multigroup model.

\subsection{The multigroup equations}

With these constants the transport equation \eqref{eq:transport} becomes
the multigroup system
\begin{equation}
\frac{1}{v_g}\pd{\psi_g}{t} + \Om\cdot\nabla\psi_g + \sigma^{g}\psi_g
= \sum_{g'=1}^{G} \int_{\Om'}
  \sigma_s^{g'\to g}(\Om'\!\cdot\!\Om)\,\psi_{g'}(\rvec,\Om',t)\,\dOm'
+ \chi_g \sum_{g'=1}^{G} \nu\sigma_f^{g'}\,\phi_{g'}
+ q_g ,
\label{eq:multigroup}
\end{equation}
for $g = 1,\dots,G$, where $q_g$ is the integral of $q$ over the group
and $\phi_g$ is the scalar flux of group $g$, the angular average
\begin{equation}
\phi_g(\rvec,t) = \int_{\Om} \psi_g(\rvec,\Om,t)\,\dOm .
\label{eq:twogroup-phi}
\end{equation}
The quotient $\phi_g/v_g$ is the neutron density of the group. The right
side couples the groups: a collision in group $g'$ feeds group $g$, and
a fission in any group feeds every group in the proportion $\chi_g$.

When the scattering is isotropic, $\sigma_s^{g'\to g}$ does not depend
on the angle, the scattering gain of group $g$ reduces to
$\sum_{g'} \sigma_s^{g'\to g}\,\phi_{g'}$, and the numbers
$\sigma_s^{g'\to g}$ form the $G\times G$ scattering matrix
\begin{equation}
\bm{\Sigma}_s = \bigl[\sigma_s^{g'\to g}\bigr]_{g,g'=1}^{G} ,
\label{eq:scatter-matrix}
\end{equation}
with the row index $g$ the group after the collision and the column
index $g'$ the group before it, so that the scattering gain of all
groups at once is the product $\bm{\Sigma}_s\,\bm{\phi}$ with
$\bm{\phi} = (\phi_1,\dots,\phi_G)^{\top}$. In a collision with a free
nucleus at rest the neutron can only lose energy, so such collisions
send neutrons from $g'$ to $g \ge g'$ and $\bm{\Sigma}_s$ is lower
triangular. The thermal motion of the nuclei of the medium lets a slow
neutron gain energy, so the thermal groups add an upscatter block above
the diagonal.

When the scattering is anisotropic, the dependence of the transfer
cross section on the cosine $\Om'\cdot\Om$ is expanded in Legendre
polynomials and truncated at an order $\ell_{\max}$,
\begin{equation}
\sigma_s^{g'\to g}(\Om'\!\cdot\!\Om)
= \sum_{\ell=0}^{\ell_{\max}} (2\ell+1)\,
  \sigma_{s,\ell}^{g'\to g}\,P_\ell(\Om'\!\cdot\!\Om) ,
\label{eq:legendre-moments}
\end{equation}
where $P_\ell$ is the Legendre polynomial of degree $\ell$ and the
moments are
\begin{equation*}
\sigma_{s,\ell}^{g'\to g}
= \int_{\Om'} \sigma_s^{g'\to g}(\Om'\!\cdot\!\Om)\,P_\ell(\Om'\!\cdot\!\Om)\,\dOm'
= \frac12 \int_{-1}^{1} \sigma_s^{g'\to g}(\mu)\,P_\ell(\mu)\,d\mu ,
\end{equation*}
the second form because under the normalized measure the cosine
$\mu = \Om'\cdot\Om$ is distributed uniformly on $[-1,1]$. The same
measure makes the Legendre polynomials orthogonal with
$\int_{\Om'} P_\ell(\Om'\!\cdot\!\Om)^2\,\dOm' = 1/(2\ell+1)$, which is where
the factor $2\ell+1$ of \eqref{eq:legendre-moments} comes from. The
moments of one order form the $G\times G$ matrix
$\bm{\Sigma}_{s,\ell} = \bigl[\sigma_{s,\ell}^{g'\to g}\bigr]_{g,g'=1}^{G}$,
read as \eqref{eq:scatter-matrix}. The moment of order zero is
the transfer cross section averaged over the direction of arrival, so
the matrix of order zero is the isotropic matrix,
$\bm{\Sigma}_{s,0} = \bm{\Sigma}_s$. Our examples use isotropic scattering, $\ell_{\max} = 0$, except for
the active anisotropic burst, which uses $\ell_{\max}=3$.
The conditions \eqref{eq:transport-ic} and \eqref{eq:transport-bc}
carry over group by group. The initial state is
\begin{equation}
\psi_g(\rvec,\Om,0) = \psi_g^{0}(\rvec,\Om)
\qquad \text{for } \rvec \in \mathcal{D} ,
\label{eq:manygroup-ic}
\end{equation}
and on the boundary, the incoming flux of each group is prescribed by
the group integral $\psi_{\mathrm{in},g}$ of the inflow datum
$\psi_{\mathrm{in}}$,
\begin{equation}
\psi_g(\rvec,\Om,t) = \psi_{\mathrm{in},g}(\rvec,\Om,t)
\qquad \text{for } \rvec \in \partial\mathcal{D}
\text{ and } \Om\cdot\bn(\rvec) < 0 ,
\label{eq:group-bc}
\end{equation}
imposed at every time $t$. For a vacuum boundary this is
\begin{equation}
\psi_g(\rvec,\Om,t) = 0
\qquad \text{for } \rvec \in \partial\mathcal{D}
\text{ and } \Om\cdot\bn(\rvec) < 0 .
\label{eq:manygroup-bc}
\end{equation}

\section{Discretization}
\label{sec:dg}

The multigroup problem of Section~\ref{sec:problem} is discretized in
direction, in space and in time, in this order. The unknown throughout is
the group flux $\psi_g$ of \eqref{eq:multigroup}.

\subsection{Directions}
\label{sec:dg-directions}

The direction $\Om$ is restricted to a quadrature set
$\{\Om_m, w_m\}_{m=1}^{M}$ on the unit sphere, with the weights normalized
so that $\sum_{m} w_m = 1$, like the angular measure of
Section~\ref{sec:problem}. The set is a product rule: Gauss-Legendre points
and weights in the polar cosine, equally spaced azimuths, taken over all
eight octants.
The unknowns are the fluxes along these directions,
\begin{equation*}
\psi_{g,m}(\rvec,t) = \psi_g(\rvec,\Om_m,t) ,
\qquad g = 1,\dots,G , \quad m = 1,\dots,M ,
\end{equation*}
and every integral over the direction becomes the quadrature sum. This
turns \eqref{eq:multigroup} into $GM$ equations that are still continuous
in space and time,
\begin{equation}
\begin{aligned}
&\frac{1}{v_g}\,\pd{\psi_{g,m}}{t}
+ \Om_m\cdot\nabla\psi_{g,m}
+ \sigma^{g}\,\psi_{g,m} \\
&\qquad
= \sum_{g'=1}^{G} \sum_{\ell=0}^{\ell_{\max}} \sigma_{s,\ell}^{g'\to g}
  \sum_{m'=1}^{M} (\mathbf{P}_\ell)_{mm'}\,\psi_{g',m'}
+ \chi_g \sum_{g'=1}^{G} \nu\sigma_f^{g'}
  \sum_{m'=1}^{M} w_{m'}\,\psi_{g',m'}
+ q_{g,m} ,
\end{aligned}
\label{eq:graphite-sn}
\end{equation}
for $g = 1,\dots,G$ and $m = 1,\dots,M$, where
$q_{g,m}(\rvec,t) = q_g(\rvec,\Om_m,t)$ is the source along the direction
$\Om_m$ and
\begin{equation*}
(\mathbf{P}_\ell)_{mm'} = (2\ell+1)\,P_\ell(\Om_m\cdot\Om_{m'})\,w_{m'} ,
\qquad m, m' = 1,\dots,M ,
\end{equation*}
is the $M\times M$ matrix that the Legendre expansion
\eqref{eq:legendre-moments} of the transfer cross section leaves once the
integral over $\Om'$ is replaced by the quadrature sum. For $\ell = 0$ the
Legendre polynomial is one, so with $\mathbf{w} = (w_1,\dots,w_M)^{\top}$
the vector of the weights and $\mathbf{1}$ the vector of $M$ ones,
$\mathbf{P}_0 = \mathbf{1}\mathbf{w}^{\top}$ replaces the flux in every
direction by its angular average. For isotropic scattering the scattering
gain of \eqref{eq:graphite-sn} is therefore
$\sum_{g'} \sigma_s^{g'\to g} \sum_{m'} w_{m'}\,\psi_{g',m'}$, the transfer
cross section times the angular average of the flux of group $g'$, and this is the form used by the isotropic examples. Each equation of
\eqref{eq:graphite-sn} is a transport equation along one fixed direction
$\Om_m$, a first-order equation in $\rvec$ and $t$ with a constant
coefficient in front of the gradient, and the $GM$ equations are coupled
only through their right-hand sides. The inflow condition
\eqref{eq:group-bc} is imposed direction by direction: on a face of the
boundary, $\psi_{g,m} = \psi_{\mathrm{in},g}(\rvec,\Om_m,t)$ for every
$\Om_m$ that points into the domain there, $\Om_m\cdot\mathbf{n} < 0$, and
the datum is zero on a vacuum face.

\subsection{Space}
\label{sec:dg-space}

The cube $\mathcal{D} = (0,L)^3$ is divided into $N$ equal cells per axis,
of width $h = L/N$, $N^3$ cells in all. On each cell the flux $\psi_{g,m}$ is
approximated by a polynomial of degree one in each coordinate, written in
the Legendre basis: along each axis a cell carries two coefficients, the
mean of the flux over the cell and its slope, so a cell carries eight
coefficients and an axis $n = 2N$ of them, $n^3$ per group and direction.
The coefficients of an axis are ordered cell by cell, the mean before the
slope inside each cell, the cell-major, basis-minor ordering that the
quantization of Section~\ref{sec:qtt} reads digit by digit. The
approximation is discontinuous across the cell faces: neighbouring cells
need not agree on their common face. In one dimension we write $\varphi_i$,
$i = 1,\dots,n$, for the basis functions, each supported on one cell, where
it is the constant or the linear Legendre polynomial, and in three
dimensions the basis functions of a cell are the products
\begin{equation*}
\varphi_{i_x i_y i_z}(x,y,z)
= \varphi_{i_x}(x)\,\varphi_{i_y}(y)\,\varphi_{i_z}(z) ,
\end{equation*}
so that on a cell the discrete flux is
$\sum \Psi_{i_x i_y i_z}\,\varphi_{i_x i_y i_z}$, the sum running over
the eight index triples of that cell.

Each equation of \eqref{eq:graphite-sn} is multiplied by a test function
$\varphi$ of a cell $K$ of the mesh, one of its eight basis functions, and
integrated over the cell. The streaming term is integrated by parts, which
moves the derivative onto the test function and produces the flux through
the faces of the cell:
\begin{align}
\frac{1}{v_g}\int_K \pd{\psi_{g,m}}{t}\,\varphi\,dV
&- \int_K \psi_{g,m}\,\Om_m\!\cdot\!\nabla\varphi\,dV
+ \int_{\partial K} \widehat{(\Om_m\!\cdot\!\mathbf{n})\,\psi_{g,m}}\,\varphi\,dS
\notag\\
&+ \int_K \sigma^{g}\psi_{g,m}\,\varphi\,dV
= \int_K Q_{g,m}\,\varphi\,dV ,
\label{eq:weak}
\end{align}
where $\mathbf{n}$ is the outward unit normal on $\partial K$ and
$Q_{g,m}$ is the right-hand side of \eqref{eq:graphite-sn}, the scattering
and fission gains and the source $q_{g,m}$. On a face the flux is taken
upwind,
\begin{equation}
\widehat{(\Om_m\!\cdot\!\mathbf{n})\,\psi_{g,m}} =
\begin{cases}
(\Om_m\!\cdot\!\mathbf{n})\,\psi_{g,m}^{\text{int}}, & \Om_m\!\cdot\!\mathbf{n} \geq 0,\\[4pt]
(\Om_m\!\cdot\!\mathbf{n})\,\psi_{g,m}^{\text{ext}}, & \Om_m\!\cdot\!\mathbf{n} < 0,
\end{cases}
\label{eq:upwind}
\end{equation}
where int and ext denote the traces on the face from the interior and from
the exterior of $K$. On a face of the cube there is no exterior cell, and
the exterior trace is the prescribed inflow datum of
\eqref{eq:group-bc}, zero for vacuum. This is how boundary data enter
the discrete system. The upwind flux takes the face value from the cell the
neutrons come from, which is what makes the discrete streaming operator
stable.

The cells are boxes and the direction
$\Om_m = (\Omega_{m,x},\Omega_{m,y},\Omega_{m,z})$ is constant, so the
streaming term splits into its three coordinate parts,
$\Om_m\cdot\nabla
= \Omega_{m,x}\,\partial_x + \Omega_{m,y}\,\partial_y
+ \Omega_{m,z}\,\partial_z$,
and each part acts on the coefficients of one axis only. With the product
basis, the volume integral of \eqref{eq:weak} for the $x$ part is the
product of a one-dimensional integral in $x$, which carries the derivative,
and two one-dimensional mass integrals in $y$ and $z$. The face integrals
factor in the same way: on a face normal to $x$ the integrand is the upwind
value in $x$ times the product of the $y$ and $z$ basis functions, and the
integral over the face is a product of one-dimensional integrals. 
On the face between a cell and its right neighbour, with int the trace
from the cell and ext the trace from the neighbour, the upwind flux
\eqref{eq:upwind} of $\mu\psi$ is $\mu$ times the average of the two
traces plus $|\mu|/2$ times their jump,
\begin{equation*}
\widehat{\mu\psi} = \mu\,\{\psi\} + \frac{|\mu|}{2}\,\jump{\psi}
= \frac{\mu}{2}\bigl(\psi^{\text{int}} + \psi^{\text{ext}}\bigr)
+ \frac{|\mu|}{2}\bigl(\psi^{\text{int}} - \psi^{\text{ext}}\bigr) ,
\end{equation*}
with $\mu$ the direction cosine $\Omega_{m,x}$: the average changes sign
with the direction, the jump does not. The same formula serves the left
face of the cell once int is read as the trace from the side of smaller
$x$ and ext as the trace from the side of larger $x$, so that for
$\mu > 0$ it always returns the trace from the left, where the neutrons
come from. Hence a derivative $\mu\,\partial_x$,
with its face terms, becomes the $n\times n$ matrix
\begin{equation}
\mathbf{A}(\mu) = \mu\,\mathbf{C} + |\mu|\,\mathbf{D} ,
\label{eq:graphite-A}
\end{equation}
where $\mathbf{C}$ is the part of the upwind operator that changes sign
with the direction (the derivative with the average of the two face values)
and $\mathbf{D}$ is the part that does not (the jump of the flux across the
faces, which is the dissipation of the upwind flux). The cell mass matrix
is diagonal in the Legendre basis and is inverted into $\mathbf{C}$ and
$\mathbf{D}$. This is the convention of the whole paper: no mass matrix
appears below, and the space factors of every operator other than the
streaming are identities. Keeping the mass matrix explicit would put the
block-diagonal one-dimensional mass matrix in place of those identities,
and since it too has quantized rank one, nothing in
Section~\ref{sec:theory} would change. The same two matrices serve all
three axes, all directions and all groups.

With the trace vectors $\mathbf{l} = \varphi(-1)$ and
$\mathbf{r} = \varphi(+1)$, which hold the values of the two basis
functions of a cell at its left and right ends in the local coordinate,
$\mathbf{C}$ and $\mathbf{D}$ are block tridiagonal in the cell index with
constant $2\times 2$ blocks. The diagonal block holds the volume term and
the face terms of the cell itself, built from $\mathbf{r}\mathbf{r}^{\top}$
and $\mathbf{l}\mathbf{l}^{\top}$, and the off-diagonal blocks couple a
cell to its two neighbours through $\mathbf{r}\mathbf{l}^{\top}$ and
$\mathbf{l}\mathbf{r}^{\top}$. At the two ends of the axis there is no
neighbour, the exterior trace is the inflow datum, and for vacuum the
inflow block of the boundary cell is absent. Naming the three blocks of
$\mathbf{C}$ as $\mathbf{C}_0$, $\mathbf{C}_+$ and $\mathbf{C}_-$, the
diagonal block and the couplings to the right and to the left
neighbour, and those of $\mathbf{D}$ likewise.  Let's define the shift matrix,
$\mathbf{Z}_N$, for $N$ cells:
$\mathbf{Z}_N\mathbf{e}_e = \mathbf{e}_{e+1}$ for $e < N$ and
$\mathbf{Z}_N\mathbf{e}_N = \mathbf{0}$. In this notation, we have,
\begin{equation}
\mathbf{C} = \mathbf{C}_0\otimes\mathbf{I}_N
+ \mathbf{C}_+\otimes\mathbf{Z}_N
+ \mathbf{C}_-\otimes\mathbf{Z}_N^{\top},
\qquad
\mathbf{D} = \mathbf{D}_0\otimes\mathbf{I}_N
+ \mathbf{D}_+\otimes\mathbf{Z}_N
+ \mathbf{D}_-\otimes\mathbf{Z}_N^{\top},
\label{eq:graphite-CD-blocks}
\end{equation}
Throughout the
paper a Kronecker product lists its factors in the order of the indices
they act on, and a coefficient tensor is vectorized with its first index
varying fastest, the column-major order of the implementation; in the
convention of MATLAB's \texttt{kron}, in which the last factor acts on
the fastest index, \eqref{eq:graphite-CD-blocks} reads
$\mathbf{C} = \mathbf{I}_N\otimes\mathbf{C}_0 + \mathbf{Z}_N\otimes\mathbf{C}_+
+ \mathbf{Z}_N^{\top}\otimes\mathbf{C}_-$, and the products of
Section~\ref{sec:dg-semidiscrete} are reversed in the same way. 
The blocks do not depend
on the cell, so $\mathbf{C}$ and $\mathbf{D}$ have the same structure
on every mesh, which is what the rank bound of
Proposition~\ref{prop:rank-operator} rests on.

Sources and inflow data enter the discrete system through their
coefficients, which are computed cell by cell. Along one axis, the cell
$K_e = (a_e, a_e + h)$, $e = 1,\dots,N$, carries the two basis functions
$1$ and $\xi$, where $\xi = 2(x - a_e)/h - 1$ is the local coordinate in
$(-1,1)$. They are orthogonal on the cell. The best approximation in the
mean-square sense of a function $u(x)$ on $K_e$ by $c_0 + c_1\,\xi$ has
therefore the coefficients
\begin{equation}
c_0 = \frac{1}{2}\int_{-1}^{1} u\,d\xi ,
\qquad
c_1 = \frac{3}{2}\int_{-1}^{1} u\,\xi\,d\xi ,
\label{eq:graphite-proj}
\end{equation}
the mean of $u$ over the cell and its slope. Collecting the two
coefficients of every cell gives the vector of $n$ coefficients of $u$
along the axis. In three dimensions the basis functions are the products
$\varphi_{i_x i_y i_z}$, and when the function is itself a product of a
function of $x$, one of $y$ and one of $z$, the three-dimensional
integrals that define its coefficients split into products of
one-dimensional integrals of the form \eqref{eq:graphite-proj}, so that
the coefficient with index $(i_x,i_y,i_z)$ is the product of the three
one-dimensional coefficients. The tensor upwind flux of~\cite{gorodetsky2025thermal} uses rank-one
direction masks and is also algebraic. Our construction instead uses
the same fixed pair of DG matrices for every direction and group,
with analytic direction-cosine factors. Section~\ref{sec:rank} derives
explicit QTT rank bounds for this construction on uniform Cartesian
grids in a homogeneous medium.

\subsection{The semi-discrete system}
\label{sec:dg-semidiscrete}

All coefficients are collected in one vector $\bm{\Psi}(t)$ of length
$n^3 M G$, indexed by $(i_x, i_y, i_z, m, g)$: the position of the
coefficient along each of the three axes, the direction and the group. In
this ordering every term of \eqref{eq:graphite-sn}, once integrated against
the basis as in \eqref{eq:weak}, becomes a sum of a few Kronecker products
of five small matrices, one per index, and each product acts on
$\bm{\Psi}$ by applying its matrices to their own indices. This product
structure rests on the uniform Cartesian grid, with the same cells, the
same basis and the same operators throughout the cube; geometric changes
such as a boundary patch or a local refinement can increase the
complexity of the representation. With
$\mathbf{I}_n$, $\mathbf{I}_M$, $\mathbf{I}_G$ the identity matrices of the
three sizes,
$\bm{\Omega}_x = \operatorname{diag}(\Omega_{1,x},\dots,\Omega_{M,x})$
the diagonal matrix of the $x$ direction cosines, $|\bm{\Omega}_x|$ the
diagonal matrix of their absolute values, and likewise for $y$ and $z$, the
streaming term is the sum of six products,
\begin{equation}
\begin{aligned}
\mathbf{L} =\ &
\mathbf{C}\otimes\mathbf{I}_n\otimes\mathbf{I}_n\otimes\bm{\Omega}_x\otimes\mathbf{I}_G
+ \mathbf{D}\otimes\mathbf{I}_n\otimes\mathbf{I}_n\otimes|\bm{\Omega}_x|\otimes\mathbf{I}_G
\\
&+ \mathbf{I}_n\otimes\mathbf{C}\otimes\mathbf{I}_n\otimes\bm{\Omega}_y\otimes\mathbf{I}_G
+ \mathbf{I}_n\otimes\mathbf{D}\otimes\mathbf{I}_n\otimes|\bm{\Omega}_y|\otimes\mathbf{I}_G
\\
&+ \mathbf{I}_n\otimes\mathbf{I}_n\otimes\mathbf{C}\otimes\bm{\Omega}_z\otimes\mathbf{I}_G
+ \mathbf{I}_n\otimes\mathbf{I}_n\otimes\mathbf{D}\otimes|\bm{\Omega}_z|\otimes\mathbf{I}_G ,
\end{aligned}
\label{eq:graphite-L}
\end{equation}
which is \eqref{eq:graphite-A} applied on each axis. The medium is
homogeneous in all our examples, so the group constants do not
depend on position and the space factors of the time, removal, scattering
and fission operators are identities. The time term, the removal and the
isotropic scattering are single products,
\begin{equation}
\begin{aligned}
\mathbf{V}^{-1} &=
\mathbf{I}_n\otimes\mathbf{I}_n\otimes\mathbf{I}_n\otimes\mathbf{I}_M
\otimes\operatorname{diag}(1/v_1,\dots,1/v_G) , \\
\bm{\Sigma}_t &=
\mathbf{I}_n\otimes\mathbf{I}_n\otimes\mathbf{I}_n\otimes\mathbf{I}_M
\otimes\operatorname{diag}(\sigma^{1},\dots,\sigma^{G}) , \\
\mathbf{S} &=
\mathbf{I}_n\otimes\mathbf{I}_n\otimes\mathbf{I}_n
\otimes\mathbf{1}\mathbf{w}^{\top}\otimes\bm{\Sigma}_s ,
\end{aligned}
\label{eq:graphite-VSS}
\end{equation}
where $\mathbf{V}$ is the diagonal matrix that carries the speed $v_g$ on
every coefficient of group $g$, so that $\mathbf{V}^{-1}$ carries $1/v_g$,
$\bm{\Sigma}_t$ carries the total cross section $\sigma^{g}$ in the same
way, $\mathbf{1}\mathbf{w}^{\top}$ is the angular average of
Section~\ref{sec:dg-directions}, and $\bm{\Sigma}_s$ is the $G\times G$
scattering matrix \eqref{eq:scatter-matrix}. {For anisotropic scattering of
Legendre order $\ell_{\max}$ the scattering operator discretizes
\eqref{eq:legendre-moments} and is the sum
\begin{equation}
\mathbf{S} = \sum_{\ell=0}^{\ell_{\max}}
\mathbf{I}_n\otimes\mathbf{I}_n\otimes\mathbf{I}_n\otimes\mathbf{P}_\ell
\otimes\bm{\Sigma}_{s,\ell} ,
\label{eq:rank-K-scatter}
\end{equation}
with $\mathbf{P}_\ell$ the angular matrices of \eqref{eq:graphite-sn} and
$\bm{\Sigma}_{s,\ell}$ the $G\times G$ matrices of the Legendre moments of
the transfer cross sections}, and for $\ell_{\max} = 0$ it reduces to the
single product of \eqref{eq:graphite-VSS}. The fission gain is one
product,
\begin{equation}
\mathbf{F} = \mathbf{I}_n\otimes\mathbf{I}_n\otimes\mathbf{I}_n\otimes
\mathbf{1}\mathbf{w}^{\top}\otimes\bm{\chi}\,(\nu\bm{\sigma}_f)^{\top} ,
\label{eq:rank-K-fission}
\end{equation}
where $\bm{\chi} = (\chi_1,\dots,\chi_G)^{\top}$ is the fission spectrum
and $\nu\bm{\sigma}_f = (\nu\sigma_f^{1},\dots,\nu\sigma_f^{G})^{\top}$
the vector of the production cross sections over the groups: the angular
average forms the scalar flux of each group, which sets the fission rate,
the vector of ones emits the new neutrons alike in every direction, and
the group factor, a matrix of rank one, distributes them over the groups
with the spectrum. The discretization in direction and space therefore
turns
\eqref{eq:multigroup} into a system of ordinary differential equations in
time,
\begin{equation}
\mathbf{V}^{-1}\,\frac{d\bm{\Psi}}{dt}
+ \bigl(\mathbf{L} + \bm{\Sigma}_t - \mathbf{S} - \mathbf{F}\bigr)\,\bm{\Psi}
= \mathbf{Q}(t) ,
\qquad \bm{\Psi}(0) = \bm{\Psi}^{0} ,
\label{eq:semidiscrete}
\end{equation}
where $\bm{\Psi}^{0}$ holds the coefficients of the initial state and
$\mathbf{Q}(t)$ holds the coefficients of the source $q_{g,m}$, obtained
by \eqref{eq:graphite-proj}, together with the contribution of any nonzero
inflow data, which enter the face integral of \eqref{eq:weak} through the
exterior trace of \eqref{eq:upwind} and are moved to the right-hand side
because they are known. 

\subsection{Time discretization}
\label{sec:dg-time}
\label{sec:dg-ttmatrices}

The streaming and removal terms of \eqref{eq:semidiscrete} are stiff and
the scattering couples all directions and groups, so an explicit
integration would be restricted to very small steps. We therefore
integrate in time by the $\theta$-method. Let
$t_0 = 0 < t_1 < t_2 < \dots$ be the time levels,
$\Delta t_n = t_{n+1} - t_n$ the step, $\bm{\Psi}^{n}$ the approximation
of $\bm{\Psi}(t_n)$, and
$\mathcal{A} = \mathbf{L} + \bm{\Sigma}_t - \mathbf{S} - \mathbf{F}$ the
operator of \eqref{eq:semidiscrete}. For a parameter $\theta \in [0,1]$
the step from $t_n$ to $t_{n+1}$ reads
\begin{equation}
\mathbf{V}^{-1}\,\frac{\bm{\Psi}^{n+1} - \bm{\Psi}^{n}}{\Delta t_n}
+ \mathcal{A}\bigl(\theta\,\bm{\Psi}^{n+1} + (1-\theta)\,\bm{\Psi}^{n}\bigr)
= \theta\,\mathbf{Q}(t_{n+1}) + (1-\theta)\,\mathbf{Q}(t_n) ,
\label{eq:theta-step}
\end{equation}
for $n = 0, 1, 2, \dots$, starting from $\bm{\Psi}^{0}$. The choice
$\theta = 0$ is the explicit forward Euler method, which the stiffness
excludes, $\theta = 1$ is backward Euler, of first order in
$\Delta t_n$, and $\theta = 1/2$ is Crank--Nicolson, the trapezoidal
rule, of second order. The spatial scheme is of second order, so the
examples of Section~\ref{sec:examples} use $\theta = 1/2$.
Multiplied by $\Delta t_n$, \eqref{eq:theta-step} is one linear system
per step,
\begin{equation}
\mathbf{K}\,\bm{\Psi}^{n+1} = \mathbf{B}\,\bm{\Psi}^{n} + \mathbf{q}^{n},
\qquad
\mathbf{K} = \mathbf{V}^{-1} + \theta\,\Delta t_n\,\mathcal{A},
\qquad
\mathbf{B} = \mathbf{V}^{-1} - (1-\theta)\,\Delta t_n\,\mathcal{A},
\label{eq:theta-matrices}
\end{equation}
with the source term
$\mathbf{q}^{n} = \Delta t_n\bigl(\theta\,\mathbf{Q}(t_{n+1})
+ (1-\theta)\,\mathbf{Q}(t_n)\bigr)$. We call $\mathbf{K}$ the physics
matrix of the time step. \ref{rem:term-count} gives its
Kronecker term count.

\subsection{Time stepping}
\label{sec:solution-procedures}

Both time-stepping methods use the hybrid TT/QTT operators of
Section~\ref{sec:dg-ttmatrices}, rounded at $10^{-12}$. The
manufactured-solution tolerance diagnostic explicitly uses
$10^{-10}$ instead. The tolerance
$\varepsilon$ is the rounding budget of one step. Lemma~\ref{lem:recursion}
of Section~\ref{sec:error-propagation} carries the residual of every
step, the rounding included, to the error at the final time, and
Lemma~\ref{lem:onestep} bounds the factor by which one step amplifies
it. TT-BUG is the principal method in the performance comparisons of
Section~\ref{sec:examples} and AMEn is the comparison solver.

\subsubsection{The implicit step by AMEn}
\label{sec:dg-procedure}

The first method solves the linear system \eqref{eq:theta-matrices} of
every step with the AMEn algorithm of~\cite{dolgov2014amen}, started
from the state of the previous step, with rank enrichment two, at most
30 sweeps and a solver tolerance of $10^{-8}$, and rounds the result to
the tolerance $\varepsilon$. These are the default settings of the
solver. Every row labelled AMEn in Section~\ref{sec:examples} uses them,
without the group balance projection of
Section~\ref{sec:error-propagation}, and Table~\ref{tab:mms-diag}
reports the stopping behaviour and the true residuals for the
manufactured solution. The largest rank inside AMEn, which sets the cost
of a step, is larger than the rank of the rounded state, and the tables
report both.

\subsubsection{The rank-adaptive TT-BUG integrator}
\label{sec:dg-bug}

TT-BUG is the principal advancement method in the burst performance
tests. We adapt the rank-adaptive basis-update and Galerkin integrator
for tree tensor networks~\cite{ceruti2023tree,ceruti2022rankadaptive};
see also~\cite{lubich2015time} for tensor-train time integration.
The implementation works with the velocity-scaled matrices
$\widetilde{\mathbf K}_s$, $\widetilde{\mathbf B}_s$ and source
$\widetilde{\mathbf q}_s^n$ of Table~\ref{tab:notation}.

A basic sweep visits the cores from last to first. Let $\mathbf F_k$
be the orthonormal frame embedding the local coefficients into the
full coefficient space. Its local Galerkin step is
\begin{equation*}
(\mathbf F_k^\top\widetilde{\mathbf K}_s\mathbf F_k)\mathbf x_1
= (\mathbf F_k^\top\widetilde{\mathbf B}_s\mathbf F_k)\mathbf x_0
  +\mathbf F_k^\top\widetilde{\mathbf q}_s^n .
\end{equation*}
GMRES solves this reduced system to relative residual $10^{-8}$.
The old and updated local states augment the basis, and new directions
below $0.1\varepsilon$ times the local updated-state norm are dropped.
The root is advanced by a Galerkin step; the result is rounded and
the group balance correction is applied. Local augmentation at most
doubles its input basis dimension, but prior source and predictor
enrichment can increase that input dimension.

Source onset and basis saturation trigger additional passes, as
specified in Algorithm~\ref{alg:bug}. The auxiliary enrichment predictor
uses $\widetilde{\mathbf B}_s(\widetilde{\bm\Psi}^n+
\widetilde{\mathbf q}_s^n)$; this is deliberately different from the
physical right-hand side
$\widetilde{\mathbf B}_s\widetilde{\bm\Psi}^n+
\widetilde{\mathbf q}_s^n$. Repeated passes stop on a relative iterate-change
test or after $64$ passes. This stopping rule does not certify the
full-system step residual. The onset and saturation treatments are
numerical modifications of BUG. Their order, convergence and general
tolerance-to-error relation are not proved in this work. The defect
analysis of Section~\ref{sec:error-propagation} applies to a BUG trajectory
only with the complete step defect, including the effects of the
evolving bases and projected integration.

\section{Low-rank tensor approximations}
\label{sec:theory}

\subsection{Rank structure of the operator}
\label{sec:rank}

\paragraph{The physics matrix}
Every operator in the physics matrix $\mathbf{K}$ of
\eqref{eq:theta-matrices}, Section~\ref{sec:dg}, is a sum of Kronecker
products of five factors, one per index, so that
\begin{equation}
\mathbf{K} = \sum_{s=1}^{q}
\mathbf{K}_s^{(x)} \otimes \mathbf{K}_s^{(y)} \otimes \mathbf{K}_s^{(z)}
\otimes \mathbf{K}_s^{(m)} \otimes \mathbf{K}_s^{(g)} .
\label{eq:rank-K-sum}
\end{equation}
Here $\mathbf{K}_s^{(x)}$, $\mathbf{K}_s^{(y)}$ and $\mathbf{K}_s^{(z)}$
are $n\times n$ matrices acting on the $n = 2N$ coefficients along one
axis, $\mathbf{K}_s^{(m)}$ is an $M\times M$ matrix acting on the
direction index and $\mathbf{K}_s^{(g)}$ a $G\times G$ matrix acting on
the group index. The Kronecker product of the five acts on a coefficient
vector indexed by $(i_x,i_y,i_z,m,g)$ by applying each matrix to its own
index. 
Here $q$ is the number of products in this representation, hence an
upper bound for the minimal Kronecker rank.
As in Section~\ref{sec:dg}, the cell mass matrices are inverted into
the one-dimensional operators of \eqref{eq:graphite-A}, so the space
factors of $\mathbf{V}^{-1}$, $\bm{\Sigma}_t$, $\mathbf{S}$ and
$\mathbf{F}$ are identities. With $\mathbf{I}_n$, $\mathbf{I}_M$ and
$\mathbf{I}_G$ the identities, the time term and the removal of
\eqref{eq:theta-matrices} are
\begin{equation}
\mathbf{V}^{-1} + \theta\,\Delta t_n\,\bm{\Sigma}_t
= \mathbf{I}_n\otimes\mathbf{I}_n\otimes\mathbf{I}_n\otimes\mathbf{I}_M
   \otimes\operatorname{diag}\bigl(1/v_g + \theta\,\Delta t_n\,\sigma^{g}\bigr) .
\label{eq:rank-K-mass}
\end{equation}
The streaming operator $\mathbf{L}$ is the sum of the six
products \eqref{eq:graphite-L}, in which $\mathbf{C}$ and $\mathbf{D}$
are the sign-changing and the dissipative part of the one-dimensional
upwind derivative, the split \eqref{eq:graphite-A}. The scattering
operator $\mathbf{S}$ is the sum of the $\ell_{\max}+1$ products
\eqref{eq:rank-K-scatter}, one per Legendre order, and for $\ell = 0$
its direction factor is $\mathbf{P}_0 = \mathbf{1}\mathbf{w}^{\top}$,
the angular average. The fission operator $\mathbf{F}$ is the single
product \eqref{eq:rank-K-fission}. All three are stated in
Section~\ref{sec:dg}, and Table~\ref{tab:rank-K-factors} lists the five
factors of every term of \eqref{eq:rank-K-sum}. The term count $q$
takes one of three values. With fission and anisotropic
scattering of Legendre order $\ell_{\max}$ it is $q = 9 + \ell_{\max}$:
one term for the time derivative and the removal, six streaming terms,
$\ell_{\max} + 1$ scattering terms and one fission term. Without
fission it is $q = 8 + \ell_{\max}$. For isotropic scattering without
fission it is $q = 8$, the scattering rows of the table reduce to the
single row $s = 8$ with $\mathbf{P}_0$, the product
$\mathbf{1}\mathbf{w}^{\top}\otimes\bm{\Sigma}_s$ of
\eqref{eq:graphite-VSS}, and this is the case of the example of
Section~\ref{sec:manufactured}. Figure~\ref{fig:rank-links} shows the links
of the train and the bound of Proposition~\ref{prop:rank-operator} that
applies to each.

\begin{table}[htbp]
\centering
\caption{The five Kronecker factors of each term of the physics matrix
\eqref{eq:rank-K-sum}. The scattering rows run over the Legendre orders
$\ell = 0,\dots,\ell_{\max}$.} 
\label{tab:rank-K-factors}
\begin{tabular}{cccccc}
\toprule
$s$ & $\mathbf{K}_s^{(x)}$ & $\mathbf{K}_s^{(y)}$ & $\mathbf{K}_s^{(z)}$
    & $\mathbf{K}_s^{(m)}$ & $\mathbf{K}_s^{(g)}$ \\
\midrule
$1$ & $\mathbf{I}_n$ & $\mathbf{I}_n$ & $\mathbf{I}_n$ & $\mathbf{I}_M$
    & $\operatorname{diag}\bigl(1/v_g + \theta\,\Delta t_n\,\sigma^{g}\bigr)$ \\
\addlinespace
$2$ & $\theta\,\Delta t_n\,\mathbf{C}$ & $\mathbf{I}_n$ & $\mathbf{I}_n$ & $\bm{\Omega}_x$ & $\mathbf{I}_G$ \\
$3$ & $\theta\,\Delta t_n\,\mathbf{D}$ & $\mathbf{I}_n$ & $\mathbf{I}_n$ & $|\bm{\Omega}_x|$ & $\mathbf{I}_G$ \\
$4$ & $\mathbf{I}_n$ & $\theta\,\Delta t_n\,\mathbf{C}$ & $\mathbf{I}_n$ & $\bm{\Omega}_y$ & $\mathbf{I}_G$ \\
$5$ & $\mathbf{I}_n$ & $\theta\,\Delta t_n\,\mathbf{D}$ & $\mathbf{I}_n$ & $|\bm{\Omega}_y|$ & $\mathbf{I}_G$ \\
$6$ & $\mathbf{I}_n$ & $\mathbf{I}_n$ & $\theta\,\Delta t_n\,\mathbf{C}$ & $\bm{\Omega}_z$ & $\mathbf{I}_G$ \\
$7$ & $\mathbf{I}_n$ & $\mathbf{I}_n$ & $\theta\,\Delta t_n\,\mathbf{D}$ & $|\bm{\Omega}_z|$ & $\mathbf{I}_G$ \\
\addlinespace
$8+\ell$ & $\mathbf{I}_n$ & $\mathbf{I}_n$ & $\mathbf{I}_n$
    & $-\theta\,\Delta t_n\,\mathbf{P}_\ell$ & $\bm{\Sigma}_{s,\ell}$ \\
\addlinespace
$9+\ell_{\max}$ & $\mathbf{I}_n$ & $\mathbf{I}_n$ & $\mathbf{I}_n$
    & $-\theta\,\Delta t_n\,\mathbf{1}\mathbf{w}^{\top}$ & $\bm{\chi}\,(\nu\bm{\sigma}_f)^{\top}$ \\
\bottomrule
\end{tabular}
\end{table}

We use the following standard tensor-rank definitions~\cite{kolda2009tensor,oseledets2011tensor}.
\begin{definition}[unfolding matrices]
\label{def:unfolding}
Let $\mathcal{X} \in \mathbb{R}^{n_1\times n_2\times\cdots\times n_p}$ be a
tensor with $p$ modes of sizes $n_1,\dots,n_p$, with entries
$\mathcal{X}(i_1,\dots,i_p)$ for $1 \le i_j \le n_j$. Then for
$k = 1,\dots,p-1$, the $k$-th unfolding matrix of $\mathcal{X}$ is the
matrix $\mathbf{X}_k \in \mathbb{R}^{(n_1\cdots n_k)\times(n_{k+1}\cdots n_p)}$ with
\begin{equation}
\mathbf{X}_k = \mathbf{X}_k(i_1,\dots,i_k;\,i_{k+1},\dots,i_p)
= \mathcal{X}(i_1,\dots,i_p),
\label{eq:unfolding}
\end{equation}
the matrix whose rows are enumerated by the first $k$ indices and whose
columns by the last $p-k$. For ease of notation we also write
$\mathbf{X}_k(\mathcal{Y})$ for the $k$-th unfolding matrix of another
tensor $\mathcal{Y}$ with the same modes, so that
$\mathbf{X}_k = \mathbf{X}_k(\mathcal{X})$.
\end{definition}

\begin{definition}[TT-ranks]
\label{def:rank}
For $k = 1,\dots,p-1$, the $k$-th TT-rank of $\mathcal{X}$ is
$r_k = \operatorname{rank}\mathbf{X}_k$, i.e., the rank of its $k$-th
unfolding.
\end{definition}

\paragraph{The modes}
From here on every tensor has the $3d + 2$ modes of the quantized form
of Section~\ref{sec:qtt}, the $d$ digits of $i_x$, the $d$ digits of
$i_y$, the $d$ digits of $i_z$, the direction $m$ and the group $g$, in
this order. Four of its unfoldings have complete indices in their rows
and in their columns,
\begin{equation}
\begin{aligned}
\mathbf{X}_{d}    &= \mathbf{X}_{d}(i_x;\, i_y, i_z, m, g), &
\mathbf{X}_{2d}   &= \mathbf{X}_{2d}(i_x, i_y;\, i_z, m, g), \\
\mathbf{X}_{3d}   &= \mathbf{X}_{3d}(i_x, i_y, i_z;\, m, g), &
\mathbf{X}_{3d+1} &= \mathbf{X}_{3d+1}(i_x, i_y, i_z, m;\, g),
\end{aligned}
\label{eq:unfold-indices}
\end{equation}
and we say that these four values of $k$ separate two different
indices. Every other $k$ lies strictly inside the digits of one
position index, $k = (j-1)d + \ell$ with $j \in \{1,2,3\}$ and
$1 \le \ell \le d-1$, and then the first $\ell$ digits of the $j$-th
position index are in the rows and its remaining $d-\ell$ digits in the
columns; for the first axis,
\begin{equation}
\mathbf{X}_{\ell} = \mathbf{X}_{\ell}(\beta_1,\dots,\beta_\ell;\,
\beta_{\ell+1},\dots,\beta_d, i_y, i_z, m, g),
\qquad 1 \le \ell \le d-1,
\label{eq:unfold-digits}
\end{equation}
and we say that such a $k$ separates two digits of that index. An
$n\times n$ matrix acting on one axis, $n = 2^{d}$, is written with the
$d$ digit modes of that axis alone, and its $d-1$ TT-ranks are its
quantized ranks (Section~\ref{sec:qtt}). A matrix acting on the tensor
is read, as in Section~\ref{sec:tt}, as a tensor whose modes carry the
pair of a row value and a column value of the corresponding
index~\cite[Section~3.1]{oseledets2011tensor}, so that
Definitions~\ref{def:unfolding} and \ref{def:rank} apply to it with
the pairs in place of the indices. For a Kronecker product
$\mathbf{A}\otimes\mathbf{B}$ in which $\mathbf{A}$ acts on the modes
in the rows of $\mathbf{X}_k$ and $\mathbf{B}$ on the modes in its
columns, the entry at a row tuple and a column tuple is an entry of
$\mathbf{A}$ times an entry of $\mathbf{B}$, so the unfolding is the
outer product
\begin{equation}
\mathbf{X}_k(\mathbf{A}\otimes\mathbf{B})
= \operatorname{vec}(\mathbf{A})\,\operatorname{vec}(\mathbf{B})^{\top},
\label{eq:unfold-kron}
\end{equation}
where $\operatorname{vec}(\mathbf{A})$ lists the entries of $\mathbf{A}$
in the order of the row pairs of $\mathbf{X}_k$ and
$\operatorname{vec}(\mathbf{B})$ those of $\mathbf{B}$ in the order of
its column pairs; in particular $\mathbf{X}_k(\mathbf{A}\otimes\mathbf{B})$
has rank one. The following facts are
standard~\cite{oseledets2011tensor}; we record them in the notation of
the modes above because the propositions below cite them part by part.

\begin{lemma}
\label{lem:rank-facts}
Let $\mathcal{X}$ and $\mathcal{Y}$ be tensors with the modes above and
$r_k$ their TT-ranks.
\begin{enumerate}
\renewcommand{\theenumi}{(\alph{enumi})}\renewcommand{\labelenumi}{\theenumi}
\item\label{lem:rank-facts-a} For every $k$,
  $r_k(\mathcal{X} + \mathcal{Y}) \le r_k(\mathcal{X}) + r_k(\mathcal{Y})$.
\item\label{lem:rank-facts-b} If $\mathcal{X} = \mathbf{A}^{(x)}\otimes\mathbf{A}^{(y)}\otimes
  \mathbf{A}^{(z)}\otimes\mathbf{A}^{(m)}\otimes\mathbf{A}^{(g)}$ is a
  Kronecker product of five factors, one acting on each index, then
  \begin{equation}
  r_k(\mathcal{X}) = 1
  \qquad\text{for } k = d,\ 2d,\ 3d,\ 3d+1,
  \label{eq:rank-kron-indices}
  \end{equation}
  the four values of $k$ that separate two different indices.
\item\label{lem:rank-facts-c} For the same $\mathcal{X}$ and every other
  $k$, writing $r_\ell(\mathbf{A})$ for the $\ell$-th quantized rank of
  a factor $\mathbf{A}$ and $\mathbf{A}^{(1)}, \mathbf{A}^{(2)},
  \mathbf{A}^{(3)}$ for $\mathbf{A}^{(x)}, \mathbf{A}^{(y)},
  \mathbf{A}^{(z)}$,
  \begin{equation}
  r_k(\mathcal{X}) = r_\ell\bigl(\mathbf{A}^{(j)}\bigr)
  \qquad\text{for } k = (j-1)d + \ell,\quad 1 \le \ell \le d-1,
  \label{eq:rank-kron-digits}
  \end{equation}
  the $\ell$-th quantized rank of the factor on the axis whose digits
  $k$ separates.
\item\label{lem:rank-facts-d} If $\mathbf{A}^{(x)} = \mathbf{I}_n$, then every quantized rank of
  $\mathbf{A}^{(x)}$ is one, and likewise on the other axes.
\end{enumerate}
\end{lemma}
We omit the proof, since these are general properties of tensor trains;
see \cite{oseledets2011tensor}.

The following proposition is an elementary consequence of the Kronecker
form \eqref{eq:rank-K-sum} and of Lemma~\ref{lem:rank-facts}; we state it
with its proof because the bound it gives is explicit and does not depend
on the number of cells, directions or groups, nor on the time step. The
bound prevents the operator TT-ranks from increasing under refinement;
the cost of a solve still depends on the ranks of the state and inside
the solver, the mode sizes, the number of sweeps and the number of
cores.

\begin{proposition}
\label{prop:rank-operator}
Let $\mathbf{K}$ be the physics matrix \eqref{eq:theta-matrices} for a
homogeneous medium, written as the sum \eqref{eq:rank-K-sum} of $q$
Kronecker products, in the quantized format of Section~\ref{sec:qtt},
with the modes ordered as above, and let $r_k$, $k = 1,\dots,3d+1$,
denote its TT-ranks. Then
\begin{enumerate}
\renewcommand{\theenumi}{(\roman{enumi})}\renewcommand{\labelenumi}{\theenumi}
\item\label{prop:rank-operator-i} $r_k \le q$ for $k = d$, $2d$, $3d$ and $3d+1$, the four values of
  $k$ at which $\mathbf{X}_k$ separates two different indices, $q$ the
  term count;
\item\label{prop:rank-operator-ii} for every other $k$, at which
  $\mathbf{X}_k$ separates two digits of the $j$-th axis,
  \begin{equation}
  r_k \le r_\ell(\mathbf{C}) + r_\ell(\mathbf{D}) + 2(j-1) + 1
  \qquad\text{for } k = (j-1)d + \ell,\quad 1 \le \ell \le d-1,
  \label{eq:rank-K-digits}
  \end{equation}
  where $r_\ell(\mathbf{C})$ and $r_\ell(\mathbf{D})$ are the $\ell$-th
  quantized ranks of the one-dimensional matrices $\mathbf{C}$ and
  $\mathbf{D}$ of \eqref{eq:graphite-A};
\item\label{prop:rank-operator-iii} $r_\ell(\mathbf{C}) \le 3$ and
  $r_\ell(\mathbf{D}) \le 3$ for every $\ell$ and every number of cells.
\end{enumerate}
In particular every TT-rank of $\mathbf{K}$ is at most $\max(q, 11)$,
independently of the number of cells per axis, the number of directions,
the number of groups and the time step.
\end{proposition}

\begin{figure}[htbp]
\centering
\begin{tikzpicture}[x=0.95cm,
  core/.style={draw, minimum size=5.5mm, inner sep=0pt, font=\scriptsize},
  inlink/.style={thin},
  modelink/.style={line width=1.4pt},
  bound/.style={font=\scriptsize},
  brace/.style={decorate, decoration={brace, amplitude=4pt}}]
\node[core] (x1) at (0,0) {$\beta_1$};
\node[core] (x2) at (1,0) {$\beta_2$};
\node       (x3) at (2,0) {$\cdots$};
\node[core] (x4) at (3,0) {$\beta_d$};
\node[core] (y1) at (4,0) {$\gamma_1$};
\node[core] (y2) at (5,0) {$\gamma_2$};
\node       (y3) at (6,0) {$\cdots$};
\node[core] (y4) at (7,0) {$\gamma_d$};
\node[core] (z1) at (8,0) {$\delta_1$};
\node[core] (z2) at (9,0) {$\delta_2$};
\node       (z3) at (10,0) {$\cdots$};
\node[core] (z4) at (11,0) {$\delta_d$};
\node[core] (m)  at (12,0) {$m$};
\node[core] (g)  at (13,0) {$g$};
\foreach \a/\b in {x1/x2, x2/x3, x3/x4, y1/y2, y2/y3, y3/y4, z1/z2, z2/z3, z3/z4}
  \draw[inlink] (\a) -- (\b);
\foreach \a/\b in {x4/y1, y4/z1, z4/m, m/g}
  \draw[modelink] (\a) -- (\b);
\draw[brace] ([yshift=3mm]x1.north west) -- ([yshift=3mm]x4.north east)
  node[midway, above=5pt, font=\small] {$x$, $d$ digits};
\draw[brace] ([yshift=3mm]y1.north west) -- ([yshift=3mm]y4.north east)
  node[midway, above=5pt, font=\small] {$y$};
\draw[brace] ([yshift=3mm]z1.north west) -- ([yshift=3mm]z4.north east)
  node[midway, above=5pt, font=\small] {$z$};
\node[font=\small, rotate=60, anchor=south west, inner sep=1pt] at ([yshift=2mm]m.north) {direction};
\node[font=\small, rotate=60, anchor=south west, inner sep=1pt] at ([yshift=2mm]g.north) {group};
\node[bound] at (1.5,-0.75) {$r_{\mathbf{C}} + r_{\mathbf{D}} + 1 \ (\le 7)$};
\node[bound] at (5.5,-0.75) {$r_{\mathbf{C}} + r_{\mathbf{D}} + 3 \ (\le 9)$};
\node[bound] at (9.5,-0.75) {$r_{\mathbf{C}} + r_{\mathbf{D}} + 5 \ (\le 11)$};
\foreach \xx in {3.5, 7.5, 11.5, 12.5}
  \node[bound] at (\xx,-0.45) {$q$};
\node[font=\scriptsize, anchor=north] at (6.5,-1.25)
  {Thin links: inside an axis, bound \ref{prop:rank-operator-ii}. \quad Thick links: between two different indices, bound \ref{prop:rank-operator-i}.};
\end{tikzpicture}
\caption{The quantized tensor train of the physics matrix $\mathbf{K}$
and the rank bounds of Proposition~\ref{prop:rank-operator}. Boxes are
the cores: the $d$ binary digits of each axis, then the direction and
the group. The TT-rank $r_k$ sits on the link after the $k$-th core. Thin links join
two digits of the same axis and carry the bound \ref{prop:rank-operator-ii}, which grows by two for
every axis to the left. Thick links join two different indices and carry
the bound \ref{prop:rank-operator-i}, the term count $q$. The numbers in parentheses use $r_{\mathbf{C}}, r_{\mathbf{D}}
\le 3$ from \ref{prop:rank-operator-iii}.}
\label{fig:rank-links}
\end{figure}
\FloatBarrier

\begin{proof}
We first prove \ref{prop:rank-operator-i}. Let $k$ be one of the four
values $d$, $2d$, $3d$ and $3d+1$. By Definition~\ref{def:rank}, by the linearity of the unfolding, by
Lemma~\ref{lem:rank-facts}\ref{lem:rank-facts-a} and by
Lemma~\ref{lem:rank-facts}\ref{lem:rank-facts-b}, equation
\eqref{eq:rank-kron-indices}, in this order,
\begin{align*}
r_k
&= \operatorname{rank}\mathbf{X}_k(\mathbf{K})
 = \operatorname{rank}\,\mathbf{X}_k\Bigl(\sum_{s=1}^{q}
   (\mathbf{K}_s^{(x)}\otimes\mathbf{K}_s^{(y)}\otimes\mathbf{K}_s^{(z)}
   \otimes\mathbf{K}_s^{(m)}\otimes\mathbf{K}_s^{(g)})\Bigr) \\
&= \operatorname{rank}\sum_{s=1}^{q}
   \mathbf{X}_k\bigl(\mathbf{K}_s^{(x)}\otimes\mathbf{K}_s^{(y)}\otimes\mathbf{K}_s^{(z)}
   \otimes\mathbf{K}_s^{(m)}\otimes\mathbf{K}_s^{(g)}\bigr) \\
&\le \sum_{s=1}^{q}
   \operatorname{rank}\mathbf{X}_k\bigl(\mathbf{K}_s^{(x)}\otimes\mathbf{K}_s^{(y)}\otimes\mathbf{K}_s^{(z)}
   \otimes\mathbf{K}_s^{(m)}\otimes\mathbf{K}_s^{(g)}\bigr)
 = \sum_{s=1}^{q} 1 = q .
\end{align*}

We now turn to \ref{prop:rank-operator-ii}. Let
$k = (j-1)d + \ell$ with $1 \le \ell \le d-1$, so that by
\eqref{eq:unfold-digits} the unfolding $\mathbf{X}_k$ has the digits of
the first $j-1$ axes and the first $\ell$ digits of the $j$-th axis in
its rows, and the remaining $d-\ell$ digits of the $j$-th axis, the
digits of the following axes, the direction and the group in its
columns. Write
$\mathbf{T}_s = \mathbf{K}_s^{(x)}\otimes\mathbf{K}_s^{(y)}\otimes
\mathbf{K}_s^{(z)}\otimes\mathbf{K}_s^{(m)}\otimes\mathbf{K}_s^{(g)}$
for the $s$-th term of \eqref{eq:rank-K-sum} and split the rows of
Table~\ref{tab:rank-K-factors} into three sets: $A_1 = \{2j, 2j+1\}$,
the two streaming terms acting on the $j$-th axis; $A_2 = \{2,\dots,2j-1\}$,
the $2(j-1)$ streaming terms acting on one of the preceding axes; and
$A_3$, every remaining row, the time and removal term, the streaming
terms acting on a following axis, the scattering terms and the fission
term. By Definition~\ref{def:rank}, by the linearity of the unfolding
and by Lemma~\ref{lem:rank-facts}\ref{lem:rank-facts-a} applied to the
terms of $A_1$ and $A_2$ one by one and to the sum over $A_3$ as a
whole,
\begin{align}
r_k
&= \operatorname{rank}\mathbf{X}_k(\mathbf{K})
 = \operatorname{rank}\Bigl(\sum_{s\in A_1}\mathbf{X}_k(\mathbf{T}_s)
   + \sum_{s\in A_2}\mathbf{X}_k(\mathbf{T}_s)
   + \sum_{s\in A_3}\mathbf{X}_k(\mathbf{T}_s)\Bigr) \notag \\
&\le \sum_{s\in A_1}\operatorname{rank}\mathbf{X}_k(\mathbf{T}_s)
   + \sum_{s\in A_2}\operatorname{rank}\mathbf{X}_k(\mathbf{T}_s)
   + \operatorname{rank}\sum_{s\in A_3}\mathbf{X}_k(\mathbf{T}_s) .
\label{eq:rank-K-three-groups}
\end{align}
The first sum has the two terms $\mathbf{T}_{2j}$ and
$\mathbf{T}_{2j+1}$, which carry $\Delta t_n\,\mathbf{C}$ and
$\Delta t_n\,\mathbf{D}$ on the $j$-th axis. By
Lemma~\ref{lem:rank-facts}\ref{lem:rank-facts-c}, equation
\eqref{eq:rank-kron-digits}, their $k$-th TT-ranks are
$r_\ell(\mathbf{C})$ and $r_\ell(\mathbf{D})$, so the first sum is
$r_\ell(\mathbf{C}) + r_\ell(\mathbf{D})$. Each of the $2(j-1)$ terms $\mathbf{T}_s$, $s \in A_2$, of the second
sum carries $\mathbf{I}_n$ on the $j$-th axis, so by
Lemma~\ref{lem:rank-facts}\ref{lem:rank-facts-c} and
\ref{lem:rank-facts-d} its $k$-th TT-rank is one, and the second sum is
$2(j-1)$. These terms are counted one by one: the factor $\mathbf{C}$
or $\mathbf{D}$ that each carries on a preceding axis lies in the rows
of $\mathbf{X}_k$ and differs from term to term. Each term $\mathbf{T}_s$, $s \in A_3$, of the third sum carries
$\mathbf{I}_n$ on the $j$-th axis and on every preceding axis, so nothing but identities lies in the rows of
$\mathbf{X}_k(\mathbf{T}_s)$. Writing
$\mathbf{I}_n = \mathbf{I}_2\otimes\cdots\otimes\mathbf{I}_2$, $d$
copies, one per digit, the unfolding of such a term is by
\eqref{eq:unfold-kron} the rank-one matrix
\begin{equation}
\mathbf{X}_k(\mathbf{T}_s) = \mathbf{a}\,\mathbf{v}_s^{\top},
\qquad
\mathbf{a} = \operatorname{vec}\bigl(\mathbf{I}_2\otimes\cdots\otimes\mathbf{I}_2\bigr),
\label{eq:rank-K-identity-group}
\end{equation}
with $(j-1)d + \ell$ copies of $\mathbf{I}_2$ in $\mathbf{a}$, one per
row digit, and $\mathbf{v}_s$ the vectorization of the remaining
factors of $\mathbf{T}_s$, the $d-\ell$ copies of $\mathbf{I}_2$ on the
other digits of the $j$-th axis and the factors on the following axes,
the direction and the group. The vector $\mathbf{a}$ is indexed by the
row digits and their primed copies, equals one where every digit agrees
with its copy and zero elsewhere, and is the same for every term of
$A_3$; only $\mathbf{v}_s$ changes with $s$. Hence
$\sum_{s\in A_3}\mathbf{X}_k(\mathbf{T}_s)
= \mathbf{a}\,\bigl(\sum_{s\in A_3}\mathbf{v}_s\bigr)^{\top}$ has rank
at most one, and \eqref{eq:rank-K-three-groups} gives
$r_k \le r_\ell(\mathbf{C}) + r_\ell(\mathbf{D}) + 2(j-1) + 1$, which is
\eqref{eq:rank-K-digits}.

Finally we prove \ref{prop:rank-operator-iii}. By
\eqref{eq:graphite-CD-blocks}, $\mathbf{C}$ is the sum of three
Kronecker products of a constant $2\times 2$ block with the identity,
the shift and the transposed shift of the $N$ cells, and the modes of
the quantized form, equation \eqref{eq:qtt-digits}, are the basis digit
first and the cell digits after it, the order of the two factors in
\eqref{eq:graphite-CD-blocks}. Let $1 \le \ell \le d-1$. The rows of
the $\ell$-th unfolding $\mathbf{X}_\ell(\mathbf{C})$ hold the basis
digit and the $\ell - 1$ low digits of the cell number, which together
address one of $L = 2^{\ell-1}$ cells inside a block, and its columns
hold the $d - \ell$ high digits, which address one of $N/L$ blocks. The
shift either stays inside a block or carries from the last cell of a
block to the first cell of the next,
\begin{equation*}
\mathbf{Z}_N = \mathbf{Z}_L\otimes\mathbf{I}_{N/L}
+ \mathbf{J}_L\otimes\mathbf{Z}_{N/L},
\qquad
\mathbf{J}_L = \mathbf{e}_1\mathbf{e}_L^{\top},
\end{equation*}
with the first factor acting on the low digits and the second on the
high digits, and $\mathbf{Z}_N^{\top}$ splits in the same way with the
transposed factors. Inserting both into $\mathbf{C}$, collecting the
terms by their factor on the high digits and applying
\eqref{eq:unfold-kron} to each, the $\ell$-th unfolding of
$\mathbf{C}$ is
\begin{align}
\mathbf{X}_\ell(\mathbf{C})
&= \operatorname{vec}\bigl(\mathbf{C}_0\otimes\mathbf{I}_L
   + \mathbf{C}_+\otimes\mathbf{Z}_L
   + \mathbf{C}_-\otimes\mathbf{Z}_L^{\top}\bigr)\,
  \operatorname{vec}(\mathbf{I}_{N/L})^{\top} \notag \\
&\quad + \operatorname{vec}\bigl(\mathbf{C}_+\otimes\mathbf{J}_L\bigr)\,
  \operatorname{vec}(\mathbf{Z}_{N/L})^{\top}
+ \operatorname{vec}\bigl(\mathbf{C}_-\otimes\mathbf{J}_L^{\top}\bigr)\,
  \operatorname{vec}(\mathbf{Z}_{N/L}^{\top})^{\top},
\label{eq:rank-CD-unfolding}
\end{align}
a sum of three rank-one matrices, so $r_\ell(\mathbf{C}) \le 3$. For
$\ell = 1$ the low factors are $1\times 1$, $\mathbf{Z}_1 = 0$ and
$\mathbf{J}_1 = 1$, and \eqref{eq:rank-CD-unfolding} reduces to the
three terms $\operatorname{vec}(\mathbf{C}_0)\operatorname{vec}(\mathbf{I}_N)^{\top}$,
$\operatorname{vec}(\mathbf{C}_+)\operatorname{vec}(\mathbf{Z}_N)^{\top}$ and
$\operatorname{vec}(\mathbf{C}_-)\operatorname{vec}(\mathbf{Z}_N^{\top})^{\top}$.
The same holds for $\mathbf{D}$ with its blocks $\mathbf{D}_0$, $\mathbf{D}_+$ and $\mathbf{D}_-$. The three column factors
$\mathbf{I}_{N/L}$, $\mathbf{Z}_{N/L}$ and $\mathbf{Z}_{N/L}^{\top}$
are the three states, no carry, a carry upward and a carry downward,
behind the quantized rank three of the discrete Laplace
operator~\cite{kazeev2012laplace}. Nothing in the argument depends on
$d$, so the bound holds for every number of cells.

\end{proof}

The relation between the bound and the measured operator ranks is
explained in \ref{app:rank-counts}. A bound valid for any index
order is given in \ref{app:rank-order}.

\subsection{Rank of the discrete solution}
\label{sec:rank-solution}

Proposition~\ref{prop:rank-operator} concerns the matrix of the step.
The work of AMEn grows also with the ranks of the solution, and nothing
in the operator alone bounds those: they depend on the data and on the
time. What can be said in general is conditional: a state that lies
close to a tensor of known ranks can be rounded to those ranks, provided
the tolerance is at least of the size of that distance. The statement
below makes this precise for one rounding; it holds for any
number $p$ of modes, so for the tensor train with one core per index as
well as for the hybrid form, and $\|\cdot\|$ is the Frobenius norm, the
Euclidean norm of the entries.

\begin{proposition}
\label{prop:rank-solution}
Let $\widehat{\bm{\Psi}}$ and $\bm{\Psi}_\ast$ be tensors with the same $p$
modes, let $r_k^\ast$, $k = 1,\dots,p-1$, be the TT-ranks of
$\bm{\Psi}_\ast$, and let $\eta = \|\widehat{\bm{\Psi}} - \bm{\Psi}_\ast\|$.
Let $\widetilde{\bm{\Psi}}$ be the result of rounding $\widehat{\bm{\Psi}}$ to
the relative tolerance $\varepsilon$ by the TT-SVD, using the absolute discarded-tail tolerance
$\varepsilon\|\widehat{\bm{\Psi}}\|/\sqrt{p-1}$ at each unfolding~\cite[Algorithms~1 and 2]{oseledets2011tensor}.
Assume that no additional rank cap overrides this tolerance. If
\begin{equation}
\varepsilon\,\|\widehat{\bm{\Psi}}\| \;\ge\; \sqrt{p-1}\;\eta ,
\label{eq:rank-solution-hyp}
\end{equation}
then
\begin{enumerate}
\renewcommand{\theenumi}{(\roman{enumi})}\renewcommand{\labelenumi}{\theenumi}
\item\label{prop:rank-solution-i} the TT-ranks of
  $\widetilde{\bm{\Psi}}$ satisfy $r_k(\widetilde{\bm{\Psi}}) \le r_k^\ast$
  for every $k$, and
\item\label{prop:rank-solution-ii}
  $\|\widetilde{\bm{\Psi}} - \bm{\Psi}_\ast\| \le \eta + \varepsilon\,\|\widehat{\bm{\Psi}}\|$.
\end{enumerate}
\end{proposition}

\begin{proof}
The TT-SVD visits the unfoldings in the order $k = 1,\dots,p-1$ and
fixes one core at each step. At step $k$ the first $k-1$ cores are
fixed, and their product is a matrix $\mathbf{U}_k$ with
$n_1\cdots n_{k-1}$ rows and $r_{k-1}(\widetilde{\bm{\Psi}})$ columns,
with $\mathbf{U}_1 = 1$, which satisfies
$\mathbf{U}_k^{\top}\mathbf{U}_k = \mathbf{I}$ because the TT-SVD leaves
every fixed core left-orthogonal. For step $k$ we define the matrix
\begin{equation}
\mathbf{H}_k = \bigl(\mathbf{U}_k^{\top}\otimes\mathbf{I}_{n_k}\bigr)\,
\mathbf{X}_k(\widehat{\bm{\Psi}})
\label{eq:ttsvd-step}
\end{equation}
of size $(r_{k-1}(\widetilde{\bm{\Psi}})\,n_k) \times (n_{k+1}\cdots n_p)$.
Step $k$ computes the singular values
$\sigma_1 \ge \sigma_2 \ge \cdots \ge \sigma_m \ge 0$ of $\mathbf{H}_k$,
with $m$ the smaller of its two dimensions, and keeps the smallest
number $r$ of them for which the discarded tail satisfies
\begin{equation}
\Bigl(\sum_{j>r}\sigma_j^2\Bigr)^{1/2}
\;\le\; \frac{\varepsilon\,\|\widehat{\bm{\Psi}}\|}{\sqrt{p-1}} ,
\label{eq:ttsvd-tail}
\end{equation}
where $r$ is the TT-rank
$r_k(\widetilde{\bm{\Psi}})$ of the rounded tensor. We next show that $r \le r_k^\ast$. To do so, compare the singular
values of $\mathbf H_k$ with those of the original unfolding:
\begin{align}
(\mathbf{U}_k\mathbf{U}_k^{\top}\otimes\mathbf{I}_{n_k})\,\mathbf{X}_k(\widehat{\bm{\Psi}})
&= \bigl(\mathbf{U}_k\otimes\mathbf{I}_{n_k}\bigr)\,\mathbf{H}_k ,
\label{eq:ttsvd-identity}\\
\sigma_j(\mathbf{H}_k)
&= \sigma_j\bigl((\mathbf{U}_k\mathbf{U}_k^{\top}\otimes\mathbf{I}_{n_k})\,\mathbf{X}_k(\widehat{\bm{\Psi}})\bigr),
\qquad j = 1,\dots,m ,
\label{eq:ttsvd-singular}
\end{align}
where the second line follows from the first because the matrix
$\mathbf{U}_k\otimes\mathbf{I}_{n_k}$ has orthonormal columns, and
multiplication by such a matrix preserves singular values.
The matrix $\mathbf{U}_k\mathbf{U}_k^{\top}\otimes\mathbf{I}_{n_k}$ is an orthogonal projector, so
$\|(\mathbf{U}_k\mathbf{U}_k^{\top}\otimes\mathbf{I}_{n_k})\mathbf{A}\mathbf{x}\| \le \|\mathbf{A}\mathbf{x}\|$ for
every $\mathbf{x}$, and the minimax characterization of the singular
values~\cite[Theorem~8.6.1]{golub2013matrix} gives
$\sigma_j((\mathbf{U}_k\mathbf{U}_k^{\top}\otimes\mathbf{I}_{n_k})\mathbf{A}) \le \sigma_j(\mathbf{A})$ for every
$j$ and every matrix $\mathbf{A}$ with $n_1\cdots n_k$ rows. With
\eqref{eq:ttsvd-singular} and $\mathbf{A} = \mathbf{X}_k(\widehat{\bm{\Psi}})$
this gives, for every $r$,
\begin{equation}
\Bigl(\sum_{j>r}\sigma_j(\mathbf{H}_k)^2\Bigr)^{1/2}
\;\le\;
\Bigl(\sum_{j>r}\sigma_j\bigl(\mathbf{X}_k(\widehat{\bm{\Psi}})\bigr)^2\Bigr)^{1/2} .
\label{eq:ttsvd-projector}
\end{equation}

Take $r = r_k^\ast$ in \eqref{eq:ttsvd-projector}. The right side is a
discarded tail of singular values, and the theorem of Eckart and
Young~\cite{eckart1936approximation} relates such a tail to low-rank
approximation: for every matrix $\mathbf{A}$ and every matrix
$\mathbf{Y}$ of the same size with rank at most $r$,
\begin{equation}
\Bigl(\sum_{j>r}\sigma_j(\mathbf{A})^2\Bigr)^{1/2}
\;\le\; \|\mathbf{A} - \mathbf{Y}\| ,
\label{eq:eckart-young}
\end{equation}
with equality when $\mathbf{Y}$ is the singular value decomposition of
$\mathbf{A}$ truncated after $r$ terms. We apply \eqref{eq:eckart-young}
with $\mathbf{A} = \mathbf{X}_k(\widehat{\bm{\Psi}})$, $r = r_k^\ast$ and
$\mathbf{Y} = \mathbf{X}_k(\bm{\Psi}_\ast)$. This $\mathbf{Y}$ has
rank at most $r_k^\ast$, as \eqref{eq:eckart-young} requires, since the
rank of $\mathbf{X}_k(\bm{\Psi}_\ast)$ is the TT-rank $r_k^\ast$ by
definition, and
$\|\mathbf{X}_k(\widehat{\bm{\Psi}}) - \mathbf{X}_k(\bm{\Psi}_\ast)\|
= \|\widehat{\bm{\Psi}} - \bm{\Psi}_\ast\| = \eta$, since an unfolding only
rearranges the entries of a tensor. Together with the hypothesis
\eqref{eq:rank-solution-hyp},
\begin{equation}
\Bigl(\sum_{j>r_k^\ast}\sigma_j(\mathbf{H}_k)^2\Bigr)^{1/2}
\;\le\; \Bigl(\sum_{j>r_k^\ast}\sigma_j\bigl(\mathbf{X}_k(\widehat{\bm{\Psi}})\bigr)^2\Bigr)^{1/2}
\;\le\; \eta
\;\le\; \frac{\varepsilon\,\|\widehat{\bm{\Psi}}\|}{\sqrt{p-1}} .
\label{eq:ttsvd-eym}
\end{equation}
So $r = r_k^\ast$ satisfies \eqref{eq:ttsvd-tail}, and the smallest
number satisfying it is at most $r_k^\ast$. This proves
\ref{prop:rank-solution-i}. For \ref{prop:rank-solution-ii}, let $\tau_k$ be the tail that step $k$
discards,
\begin{equation}
\tau_k = \Bigl(\sum_{j>r_k(\widetilde{\bm{\Psi}})}
\sigma_j(\mathbf{H}_k)^2\Bigr)^{1/2},
\qquad k = 1,\dots,p-1 ,
\label{eq:ttsvd-tau}
\end{equation}
which by \eqref{eq:ttsvd-tail} is at most $\varepsilon\|\widehat{\bm{\Psi}}\|/\sqrt{p-1}$.
We control the error of the TT-SVD by these tails, using
Theorem~2.2 of Oseledets~\cite{oseledets2011tensor}
\begin{equation}
\|\widetilde{\bm{\Psi}} - \widehat{\bm{\Psi}}\|^2
\;\le\; \sum_{k=1}^{p-1}\tau_k^2
\;\le\; (p-1)\,\frac{\varepsilon^2\|\widehat{\bm{\Psi}}\|^2}{p-1}
\;=\; \varepsilon^2\|\widehat{\bm{\Psi}}\|^2 ,
\label{eq:ttsvd-error}
\end{equation}
where the second inequality is \eqref{eq:ttsvd-tail} at every step. The
triangle inequality then gives
\begin{equation*}
\|\widetilde{\bm{\Psi}} - \bm{\Psi}_\ast\|
\;\le\; \|\widetilde{\bm{\Psi}} - \widehat{\bm{\Psi}}\|
      + \|\widehat{\bm{\Psi}} - \bm{\Psi}_\ast\|
\;\le\; \varepsilon\,\|\widehat{\bm{\Psi}}\| + \eta ,
\end{equation*}
which concludes our proof of \ref{prop:rank-solution-ii}.
\end{proof}

\ref{app:rank-solution-interpretation} interprets the hypothesis,
including why it does not hold for the runs of Table~\ref{tab:mms}; the
proposition bounds neither AMEn ranks nor an entire computed trajectory.

\subsection{Propagation of the errors of a step}
\label{sec:error-propagation}

In this section we estimate the error of the computed states.
Table~\ref{tab:notation} fixes the notation. A plain symbol belongs to
the $\theta$-method in exact arithmetic, a tilde marks a quantity
computed or compressed in tensor-train arithmetic, a hat marks the state
the linear solver delivers before it is rounded, a bar marks an average,
and the subscript $s$ marks the system multiplied by $\mathbf{V}$, the
form the implementation assembles. The hats on the numerical fluxes of
the spatial discretization are the only other hats in the paper.

\begin{table}[!htbp]
\centering
\caption{Notation of the error analysis. A plain symbol belongs to the
$\theta$-method in exact arithmetic, a tilde to a quantity computed or
compressed in tensor-train arithmetic, a hat to the state the linear
solver delivers before rounding, a bar to an average, and the subscript
$s$ to the system multiplied by $\mathbf{V}$.}
\label{tab:notation}
\vspace{6pt}
\footnotesize
\begin{tabular}{l p{0.64\textwidth}}
\toprule
\textbf{Symbol} & \textbf{Meaning} \\
\midrule
$\bm{\Psi}^n$ & solution of the $\theta$-method step \eqref{eq:theta-matrices} in exact arithmetic \\
$\widetilde{\bm{\Psi}}^n$ & computed tensor-train state after step $n$, after the rounding and, when it is used, the balance projection \\
$\widehat{\bm{\Psi}}^{n+1}$ & state the linear solver delivers in step $n$, before rounding and projection \\
$\mathbf{K}$, $\mathbf{B}$, $\mathbf{q}^n$ & step matrices and source of \eqref{eq:theta-matrices} \\
$\widetilde{\mathbf{K}}$, $\widetilde{\mathbf{B}}$, $\widetilde{\mathbf{q}}^n$, $\widetilde{\mathbf{b}}^n$ & their compressed forms, and the compressed right-hand side handed to the solver \\
$\widetilde{\mathbf{K}}_s$, $\widetilde{\mathbf{B}}_s$, $\widetilde{\mathbf{q}}^n_s$ & compressed forms of $\mathbf{V}\mathbf{K}$, $\mathbf{V}\mathbf{B}$ and $\mathbf{V}\mathbf{q}^n$, the system the implementation assembles \\
$\mathbf{e}^n = \widetilde{\bm{\Psi}}^n - \bm{\Psi}^n$ & error of the computed state \\
$\bm{\rho}^n$, $\widehat{\bm{\rho}}^n$, $\bm{\rho}^n_s$ & residual of the computed step \eqref{eq:step-defect}, of the delivered state, and of the compressed scaled system \\
$\mathbf{r}^n$, $\mathbf{p}^n$ & change of the state by the rounding and by the balance projection, $\widetilde{\bm{\Psi}}^{n+1} = \widehat{\bm{\Psi}}^{n+1} + \mathbf{r}^n + \mathbf{p}^n$ \\
$\mathbf{W}$, $\|\cdot\|_E$ & weight matrix and energy norm of the scheme \\
$\gamma$, $\vartheta$, $\kappa$ & growth bound of Lemma~\ref{lem:energy} or \ref{app:anisotropic-energy}, amplification factor of one step, and $\operatorname{cond}(\mathbf{W}\mathbf{V}^{-1})^{1/2}$ \\
$\varepsilon$, $N_t$ & rounding tolerance, and number of steps to $T = N_t\Delta t$ \\
$J$, $\mathcal{N}$ & integral over the cube, the directions and the groups, and the neutron count $\mathcal{N} = J(\mathbf{V}^{-1}\,\cdot\,)$ \\
\bottomrule
\end{tabular}
\end{table}

A step of the $\theta$-method \eqref{eq:theta-matrices} in exact
arithmetic reads
$$\mathbf{K}\bm{\Psi}^{n+1} = \mathbf{B}\bm{\Psi}^{n} + \mathbf{q}^{n}.$$
Then we define the step residual and the error of the computed states
$\widetilde{\bm{\Psi}}^n$,
\begin{equation}
\bm{\rho}^n = \mathbf{K}\widetilde{\bm{\Psi}}^{n+1}
- \mathbf{B}\widetilde{\bm{\Psi}}^{n} - \mathbf{q}^{n},
\qquad
\mathbf{e}^n = \widetilde{\bm{\Psi}}^n - \bm{\Psi}^n .
\label{eq:step-defect}
\end{equation}
Subtracting the exact step from \eqref{eq:step-defect}, the error obeys
$\mathbf{e}^{n+1} = \mathbf{K}^{-1}\mathbf{B}\,\mathbf{e}^{n}
+ \mathbf{K}^{-1}\bm{\rho}^{n}$ and, with $\mathbf{e}^0 = 0$,
\begin{equation}
\mathbf{e}^{N_t} = \sum_{n=0}^{N_t-1}
\bigl(\mathbf{K}^{-1}\mathbf{B}\bigr)^{N_t-1-n}\,\mathbf{K}^{-1}\bm{\rho}^{n} .
\label{eq:error-sum}
\end{equation}
To identify the contributions to the step residual, insert the
compressed matrices and source together with the right-hand side
$\widetilde{\mathbf b}^n$ supplied to the solver:
\begin{equation*}
\begin{aligned}
\bm{\rho}^n ={}& \bigl(\mathbf{K} - \widetilde{\mathbf{K}}\bigr)\widetilde{\bm{\Psi}}^{n+1}
- \bigl(\mathbf{B} - \widetilde{\mathbf{B}}\bigr)\widetilde{\bm{\Psi}}^{n}
- \bigl(\mathbf{q}^n - \widetilde{\mathbf{q}}^n\bigr) \\
&+ \bigl(\widetilde{\mathbf{b}}^n - \widetilde{\mathbf{B}}\widetilde{\bm{\Psi}}^{n} - \widetilde{\mathbf{q}}^n\bigr)
+ \bigl(\widetilde{\mathbf{K}}\widetilde{\bm{\Psi}}^{n+1} - \widetilde{\mathbf{b}}^n\bigr) .
\end{aligned}
\end{equation*}
The first three terms are the effects of compressing the two step
matrices and the source. The fourth is the error in forming the
compressed right-hand side. The final term is the residual of the
compressed equation evaluated at the state carried forward, including
any subsequent rounding and balance correction. These terms are
algebraically distinct; the reported local solver residual alone does
not bound their sum. To propagate their combined effect, we first bound
the inverse step matrix and the exact step propagator in an energy norm. Let $\mathbf{M}_h$ be the
one-dimensional mass matrix, which is block diagonal and diagonal in
the Legendre basis. Recall
$\mathcal{A} = \mathbf{L} + \bm{\Sigma}_t - \mathbf{S} - \mathbf{F}$ for
the operator of \eqref{eq:semidiscrete}, and let
\begin{equation*}
\mathbf{W} = \mathbf{M}_h\otimes\mathbf{M}_h\otimes\mathbf{M}_h
\otimes\operatorname{diag}(\mathbf{w})\otimes\mathbf{I}_G ,
\qquad
\|\mathbf{x}\|_E^2 = \mathbf{x}^{\top}\mathbf{W}\mathbf{V}^{-1}\mathbf{x} .
\end{equation*}
$\|\cdot\|_E$ is the mean-square norm of the scheme, with the normalized
angular weights, summed over the groups with the weight $1/v_g$. Note $\mathbf{W}$ commutes with $\mathbf{V}$, and
$\mathbf{V}_G = \operatorname{diag}(v_1,\dots,v_G)$ denotes the group
factor of $\mathbf{V}$.

\begin{lemma}
\label{lem:energy}
Let $\mathbf{W}$ and $\|\cdot\|_E$ be as above, and define
\begin{equation*}
\gamma = \max\Bigl\{0,\; -\lambda_{\min}\Bigl(\mathbf{V}_G^{1/2}\,
\operatorname{sym}\bigl(\operatorname{diag}(\sigma^{1},\dots,\sigma^{G})
- \bm{\Sigma}_s - \bm{\chi}(\nu\bm{\sigma}_f)^{\top}\bigr)\,
\mathbf{V}_G^{1/2}\Bigr)\Bigr\} ,
\qquad
\operatorname{sym}(\mathbf{Y}) = \tfrac{1}{2}\bigl(\mathbf{Y} + \mathbf{Y}^{\top}\bigr) .
\end{equation*}
Then
\begin{enumerate}
\renewcommand{\theenumi}{(\roman{enumi})}\renewcommand{\labelenumi}{\theenumi}
\item\label{lem:energy-i} $\mathbf{M}_h\mathbf{C}$ is skew-symmetric and
$\mathbf{M}_h\mathbf{D}$ has a positive semidefinite symmetric part, i.e.,
 $\mathbf{x}^{\top}\mathbf{W}\mathbf{L}\,\mathbf{x} \ge 0$ for every
$\mathbf{x}$,
\item\label{lem:energy-ii} if the scattering is isotropic, then
$\mathbf{x}^{\top}\mathbf{W}\mathcal{A}\,\mathbf{x} \ge -\gamma\,\|\mathbf{x}\|_E^2$
for every $\mathbf{x}$.
\end{enumerate}
\end{lemma}

\begin{proof}
\ref{lem:energy-i} Consider the cells $(x_{i-1}, x_i)$, $i = 1,\dots,N$, with
$x_0 = 0$ and $x_N = L$, let $u$ be the function of the scheme with the
coefficients $\mathbf{x}$ and $\mu$ the component of the direction
along the line. The upwind form of the derivative $\mu\,\partial_x u$
in \eqref{eq:weak}, integrated by parts on every cell with the flux
\eqref{eq:upwind} on the faces and tested with $u$, is
\begin{equation*}
\mathbf{x}^{\top}\mathbf{M}_h\mathbf{A}(\mu)\,\mathbf{x}
= \sum_{i=1}^{N}\Bigl(
\underbrace{-\int_{x_{i-1}}^{x_i} \mu\,u\,u'\,dx}_{=\ \frac{\mu}{2}\bigl(u(x_{i-1}^+)^2 - u(x_i^-)^2\bigr)}
+ \underbrace{\widehat{\mu u}(x_i)\,u(x_i^-) - \widehat{\mu u}(x_{i-1})\,u(x_{i-1}^+)}_{\text{flux through the two faces of the cell}}
\Bigr) ,
\end{equation*}
where $u(x_i^-)$ and $u(x_i^+)$ are the traces of $u$ at $x_i$ from the
left and from the right, and
$\widehat{\mu u} = \mu\{u\} + \frac{|\mu|}{2}\jump{u}$ is the upwind
flux written with the average
$\{u\} = \frac12\bigl(u(x_i^-) + u(x_i^+)\bigr)$ and the jump
$\jump{u} = u(x_i^-) - u(x_i^+)$ of the two traces. Collect the terms
face by face. An interior face $x_i$ receives the right end of cell $i$
and the left end of cell $i+1$,
\begin{align*}
\underbrace{\frac{\mu}{2}\bigl(u(x_i^+)^2 - u(x_i^-)^2\bigr)}_{\text{from the two cell integrals}}
+ \underbrace{\widehat{\mu u}(x_i)\bigl(u(x_i^-) - u(x_i^+)\bigr)}_{\text{from the two face terms}}
&= \frac{\mu}{2}\bigl(u(x_i^+)^2 - u(x_i^-)^2\bigr)
+ \underbrace{\mu\{u\}\jump{u}}_{=\ \frac{\mu}{2}\bigl(u(x_i^-)^2 - u(x_i^+)^2\bigr)}
+ \frac{|\mu|}{2}\jump{u}^2 \\
&= \frac{|\mu|}{2}\jump{u}^2 ,
\end{align*}
since $2\{u\}\jump{u} = u(x_i^-)^2 - u(x_i^+)^2$ cancels the cell
terms. At the two ends of the line the exterior trace is the vacuum
datum zero, and the same computation with $u(x_0^-) = 0$ and
$u(x_N^+) = 0$ leaves $\frac{|\mu|}{2}u(0)^2$ and
$\frac{|\mu|}{2}u(L)^2$, with $u(0)$ and $u(L)$ the traces of $u$ at the
ends. Hence
\begin{equation*}
\mathbf{x}^{\top}\mathbf{M}_h\mathbf{A}(\mu)\,\mathbf{x}
= \frac{|\mu|}{2}\,
\underbrace{\Bigl(\sum_{\text{interior faces}} \jump{u}^2 + u(0)^2 + u(L)^2\Bigr)}_{=:\ J(u)\ \ge\ 0} ,
\end{equation*}
the classical identity of~\cite{lesaint1974finite}. Now the right-hand side $\frac{|\mu|}{2}J(u)$ is the same for $\mu$ and for
$-\mu$, while by \eqref{eq:graphite-A}
$\mathbf{A}(\pm\mu) = \pm\mu\,\mathbf{C} + |\mu|\,\mathbf{D}$.
Subtracting and adding the identity for $\mu$ and for $-\mu$,
\begin{gather*}
\underbrace{\mathbf{x}^{\top}\mathbf{M}_h\mathbf{A}(\mu)\,\mathbf{x}
- \mathbf{x}^{\top}\mathbf{M}_h\mathbf{A}(-\mu)\,\mathbf{x}}_{=\ 2\mu\,\mathbf{x}^{\top}\mathbf{M}_h\mathbf{C}\,\mathbf{x}}
= \frac{|\mu|}{2}J(u) - \frac{|\mu|}{2}J(u) = 0 , \\
\underbrace{\mathbf{x}^{\top}\mathbf{M}_h\mathbf{A}(\mu)\,\mathbf{x}
+ \mathbf{x}^{\top}\mathbf{M}_h\mathbf{A}(-\mu)\,\mathbf{x}}_{=\ 2|\mu|\,\mathbf{x}^{\top}\mathbf{M}_h\mathbf{D}\,\mathbf{x}}
= |\mu|\,J(u) .
\end{gather*}
Hence $\mathbf{x}^{\top}\mathbf{M}_h\mathbf{C}\,\mathbf{x} = 0$ for
every $\mathbf{x}$, so the symmetric part of $\mathbf{M}_h\mathbf{C}$
is zero and $\mathbf{M}_h\mathbf{C}$ is skew-symmetric, and
$\mathbf{x}^{\top}\mathbf{M}_h\mathbf{D}\,\mathbf{x} = \frac12 J(u) \ge 0$
for every $\mathbf{x}$, so the symmetric part of
$\mathbf{M}_h\mathbf{D}$ is positive semidefinite. By
\eqref{eq:graphite-L}, $\mathbf{L}$ is the sum of one $\mathbf{C}$
product and one $\mathbf{D}$ product per axis, rows $2$ to $7$ of
Table~\ref{tab:rank-K-factors} without the factor $\Delta t_n$. Name
them $\mathbf{L} = \sum_{a \in \{x,y,z\}} (\mathbf{L}_a^{C} + \mathbf{L}_a^{D})$
with, for the $x$ axis,
\begin{equation*}
\mathbf{L}_x^{C} = \mathbf{C}\otimes\mathbf{I}_n\otimes\mathbf{I}_n\otimes\bm{\Omega}_x\otimes\mathbf{I}_G ,
\qquad
\mathbf{L}_x^{D} = \mathbf{D}\otimes\mathbf{I}_n\otimes\mathbf{I}_n\otimes|\bm{\Omega}_x|\otimes\mathbf{I}_G ,
\end{equation*}
where the $y$ and $z$ axes differ only in the position of the factors.
With
$\mathbf{W} = \mathbf{M}_h\otimes\mathbf{M}_h\otimes\mathbf{M}_h\otimes\operatorname{diag}(\mathbf{w})\otimes\mathbf{I}_G$
and $\operatorname{diag}(\mathbf{w})\bm{\Omega}_x = \operatorname{diag}(\mathbf{w}\circ\Omega_{\cdot,x})$,
\begin{gather*}
\mathbf{W}\mathbf{L}_x^{C}
= \underbrace{(\mathbf{M}_h\mathbf{C})}_{\text{skew-symmetric}}
\otimes\underbrace{\mathbf{M}_h\otimes\mathbf{M}_h\otimes\operatorname{diag}(\mathbf{w}\circ\Omega_{\cdot,x})\otimes\mathbf{I}_G}_{\text{symmetric}} , \\
\mathbf{W}\mathbf{L}_x^{D}
= \underbrace{(\mathbf{M}_h\mathbf{D})}_{\text{symmetric part}\ \succeq\ 0}
\otimes\underbrace{\mathbf{M}_h\otimes\mathbf{M}_h\otimes\operatorname{diag}(\mathbf{w}\circ|\Omega_{\cdot,x}|)\otimes\mathbf{I}_G}_{\text{symmetric positive semidefinite}} .
\end{gather*}
The transpose of a Kronecker product is the Kronecker product of the
transposes. The first matrix is therefore skew-symmetric and its
quadratic form vanishes. The symmetric part of the second is
$\operatorname{sym}(\mathbf{M}_h\mathbf{D})\otimes\mathbf{M}_h\otimes\mathbf{M}_h\otimes\operatorname{diag}(\mathbf{w}\circ|\Omega_{\cdot,x}|)\otimes\mathbf{I}_G$,
a Kronecker product of positive semidefinite matrices, so its quadratic
form is nonnegative. Summing over the three axes,
\begin{equation*}
\mathbf{x}^{\top}\mathbf{W}\mathbf{L}\,\mathbf{x}
= \sum_{a \in \{x,y,z\}}\Bigl(
\underbrace{\mathbf{x}^{\top}\mathbf{W}\mathbf{L}_a^{C}\,\mathbf{x}}_{=\ 0}
+ \underbrace{\mathbf{x}^{\top}\mathbf{W}\mathbf{L}_a^{D}\,\mathbf{x}}_{\ge\ 0}
\Bigr) \ge 0 .
\end{equation*}

\ref{lem:energy-ii} By \ref{lem:energy-i} the streaming part satisfies
$\mathbf{x}^{\top}\mathbf{W}\mathbf{L}\,\mathbf{x} \ge 0$, so it remains
to bound the collision part
$\mathbf{x}^{\top}\mathbf{W}(\bm{\Sigma}_t - \mathbf{S} - \mathbf{F})\mathbf{x}$
from below. For isotropic scattering the three collision operators and
$\mathbf{V}^{-1}$ act on the direction and the group only
(Table~\ref{tab:rank-K-factors}),
\begin{gather*}
\bm{\Sigma}_t = \mathbf{I}_n\otimes\mathbf{I}_n\otimes\mathbf{I}_n\otimes\mathbf{I}_M\otimes\operatorname{diag}(\sigma^{1},\dots,\sigma^{G}),
\qquad
\mathbf{V}^{-1} = \mathbf{I}_n\otimes\mathbf{I}_n\otimes\mathbf{I}_n\otimes\mathbf{I}_M\otimes\mathbf{V}_G^{-1}, \\
\mathbf{S} = \mathbf{I}_n\otimes\mathbf{I}_n\otimes\mathbf{I}_n\otimes\mathbf{P}_0\otimes\bm{\Sigma}_s,
\qquad
\mathbf{F} = \mathbf{I}_n\otimes\mathbf{I}_n\otimes\mathbf{I}_n\otimes\mathbf{P}_0\otimes\bm{\chi}(\nu\bm{\sigma}_f)^{\top},
\end{gather*}
with $\mathbf{P}_0 = \mathbf{1}\mathbf{w}^{\top}$ the angular average,
$(\mathbf{P}_0\mathbf{y})_m = \sum_{m'} w_{m'}\,y_{m'}$ for every $m$.
Write $\mathbf{M}_3 = \mathbf{M}_h\otimes\mathbf{M}_h\otimes\mathbf{M}_h$
for the space factor of $\mathbf{W}$, so that
$\mathbf{W} = \mathbf{M}_3\otimes\operatorname{diag}(\mathbf{w})\otimes\mathbf{I}_G$.

Since
$\sum_m w_m = \mathbf{w}^{\top}\mathbf{1} = 1$, $\mathbf{P}_0$ is a
projector,
$\mathbf{P}_0^2 = \mathbf{1}(\mathbf{w}^{\top}\mathbf{1})\mathbf{w}^{\top} = \mathbf{P}_0$,
and $\operatorname{diag}(\mathbf{w})\mathbf{P}_0 = \mathbf{w}\mathbf{w}^{\top}$
is symmetric. Split
\begin{equation*}
\mathbf{x}
= \underbrace{(\mathbf{I}_n\otimes\mathbf{I}_n\otimes\mathbf{I}_n\otimes\mathbf{P}_0\otimes\mathbf{I}_G)\,\mathbf{x}}_{\bar{\mathbf{x}},\ \text{the angular mean, the same in every direction}}
+ \underbrace{(\mathbf{I}_n\otimes\mathbf{I}_n\otimes\mathbf{I}_n\otimes(\mathbf{I}_M - \mathbf{P}_0)\otimes\mathbf{I}_G)\,\mathbf{x}}_{\mathbf{x}_\perp,\ \text{the fluctuation, of weighted mean zero}} .
\end{equation*}
A mean and a fluctuation are orthogonal in the inner product of
$\mathbf{W}$,
\begin{equation*}
\bar{\mathbf{x}}^{\top}\mathbf{W}\mathbf{x}_\perp
= \mathbf{x}^{\top}\Bigl(\mathbf{M}_3\otimes
\underbrace{\mathbf{P}_0^{\top}\operatorname{diag}(\mathbf{w})(\mathbf{I}_M - \mathbf{P}_0)}_{=\ \mathbf{w}\mathbf{w}^{\top} - \mathbf{w}(\mathbf{w}^{\top}\mathbf{1})\mathbf{w}^{\top}\ =\ 0}
\otimes\mathbf{I}_G\Bigr)\mathbf{x} = 0 ,
\end{equation*}
and the computation uses nothing about $\mathbf{x}$, so every mean is
orthogonal to every fluctuation. The same holds in the inner product of
$\mathbf{W}\mathbf{V}^{-1}$, whose extra factor $\mathbf{V}_G^{-1}$
sits on the group index. In particular
$\|\mathbf{x}\|_E^2 = \|\bar{\mathbf{x}}\|_E^2 + \|\mathbf{x}_\perp\|_E^2$.

$\bm{\Sigma}_t$ and $\mathbf{V}^{-1}$ act on the group index only, so
$\bm{\Sigma}_t\bar{\mathbf{x}}$ is again a mean and
$\bm{\Sigma}_t\mathbf{x}_\perp$ again a fluctuation. $\mathbf{S}$ and
$\mathbf{F}$ carry $\mathbf{P}_0$ in the direction, so they annihilate
every fluctuation,
$\mathbf{S}\mathbf{x}_\perp = (\mathbf{I}_n\otimes\mathbf{I}_n\otimes\mathbf{I}_n\otimes\mathbf{P}_0(\mathbf{I}_M - \mathbf{P}_0)\otimes\bm{\Sigma}_s)\,\mathbf{x} = 0$
and likewise $\mathbf{F}\mathbf{x}_\perp = 0$, and they map a mean to
a mean. Expanding the quadratic form in
$\mathbf{x} = \bar{\mathbf{x}} + \mathbf{x}_\perp$, the two cross
terms pair a mean with a fluctuation and vanish, so that
\begin{equation*}
\mathbf{x}^{\top}\mathbf{W}(\bm{\Sigma}_t - \mathbf{S} - \mathbf{F})\mathbf{x}
= \underbrace{\bar{\mathbf{x}}^{\top}\mathbf{W}(\bm{\Sigma}_t - \mathbf{S} - \mathbf{F})\bar{\mathbf{x}}}_{\text{the mean part}}
+ \underbrace{\mathbf{x}_\perp^{\top}\mathbf{W}\bm{\Sigma}_t\mathbf{x}_\perp}_{\text{the fluctuation part},\ \ge 0} .
\end{equation*}
The fluctuation part is nonnegative, since
$\mathbf{W}\bm{\Sigma}_t = \mathbf{M}_3\otimes\operatorname{diag}(\mathbf{w})\otimes\operatorname{diag}(\sigma^{1},\dots,\sigma^{G})$
is a Kronecker product of symmetric positive semidefinite factors.

Write
$\bar{\mathbf{x}} = (\mathbf{I}_n\otimes\mathbf{I}_n\otimes\mathbf{I}_n\otimes\mathbf{1}\otimes\mathbf{I}_G)\,\bm{\phi}$
with
$\bm{\phi} = (\mathbf{I}_n\otimes\mathbf{I}_n\otimes\mathbf{I}_n\otimes\mathbf{w}^{\top}\otimes\mathbf{I}_G)\,\mathbf{x}$
the coefficients of the angular average, indexed by position and group.
In the direction factor
$\mathbf{1}^{\top}\operatorname{diag}(\mathbf{w})\mathbf{1} = \sum_m w_m = 1$
and
$\mathbf{1}^{\top}\operatorname{diag}(\mathbf{w})\mathbf{P}_0\mathbf{1} = (\mathbf{w}^{\top}\mathbf{1})^2 = 1$,
so the direction drops out and
\begin{gather*}
\bar{\mathbf{x}}^{\top}\mathbf{W}(\bm{\Sigma}_t - \mathbf{S} - \mathbf{F})\bar{\mathbf{x}}
= \bm{\phi}^{\top}\Bigl(\mathbf{M}_3\otimes
\underbrace{\bigl(\operatorname{diag}(\sigma^{1},\dots,\sigma^{G}) - \bm{\Sigma}_s - \bm{\chi}(\nu\bm{\sigma}_f)^{\top}\bigr)}_{=:\ \bm{\Sigma}_c,\ \text{the group matrix of the collisions}}
\Bigr)\bm{\phi}
= \bm{\phi}^{\top}\bigl(\mathbf{M}_3\otimes\operatorname{sym}(\bm{\Sigma}_c)\bigr)\bm{\phi} , \\
\|\bar{\mathbf{x}}\|_E^2
= \bar{\mathbf{x}}^{\top}\mathbf{W}\mathbf{V}^{-1}\bar{\mathbf{x}}
= \bm{\phi}^{\top}\bigl(\mathbf{M}_3\otimes\mathbf{V}_G^{-1}\bigr)\bm{\phi} ,
\end{gather*}
where only the symmetric part of $\bm{\Sigma}_c$ enters the first line
because the form is quadratic.

By the definition of $\gamma$,
\begin{equation*}
\gamma \ge -\lambda_{\min}\bigl(\mathbf{V}_G^{1/2}\operatorname{sym}(\bm{\Sigma}_c)\mathbf{V}_G^{1/2}\bigr) ,
\qquad\text{that is,}\qquad
\mathbf{V}_G^{1/2}\operatorname{sym}(\bm{\Sigma}_c)\mathbf{V}_G^{1/2} + \gamma\,\mathbf{I}_G
\ \text{is positive semidefinite,}
\end{equation*}
and multiplying by $\mathbf{V}_G^{-1/2}$ from both sides, so is
$\operatorname{sym}(\bm{\Sigma}_c) + \gamma\mathbf{V}_G^{-1}$. Its
Kronecker product with the symmetric positive definite $\mathbf{M}_3$
is positive semidefinite as well, hence
\begin{equation*}
\bm{\phi}^{\top}\bigl(\mathbf{M}_3\otimes\operatorname{sym}(\bm{\Sigma}_c)\bigr)\bm{\phi}
\ge -\gamma\,\bm{\phi}^{\top}\bigl(\mathbf{M}_3\otimes\mathbf{V}_G^{-1}\bigr)\bm{\phi}
= -\gamma\,\|\bar{\mathbf{x}}\|_E^2 .
\end{equation*}

Assembling,
\begin{align*}
\mathbf{x}^{\top}\mathbf{W}\mathcal{A}\,\mathbf{x}
&= \underbrace{\mathbf{x}^{\top}\mathbf{W}\mathbf{L}\,\mathbf{x}}_{\ge 0\ \text{by \ref{lem:energy-i}}}
+ \underbrace{\bar{\mathbf{x}}^{\top}\mathbf{W}(\bm{\Sigma}_t - \mathbf{S} - \mathbf{F})\bar{\mathbf{x}}}_{\ge\, -\gamma\|\bar{\mathbf{x}}\|_E^2}
+ \underbrace{\mathbf{x}_\perp^{\top}\mathbf{W}\bm{\Sigma}_t\mathbf{x}_\perp}_{\ge 0} \\
&\ge -\gamma\,\|\bar{\mathbf{x}}\|_E^2 \\
&\ge\; -\gamma\,\|\mathbf{x}\|_E^2 ,
\end{align*}
the last step because $\|\bar{\mathbf{x}}\|_E^2 \le \|\mathbf{x}\|_E^2$
and $\gamma \ge 0$.
\end{proof}

We now derive a global error estimate for uniform time steps of length
$\Delta t$. For isotropic scattering, let $\gamma$ be the value in
Lemma~\ref{lem:energy}; for anisotropic scattering, use
$\gamma_\Omega$ from \ref{app:anisotropic-energy} instead.
The argument below requires only the corresponding inequality
$\mathbf x^\top\mathbf W\mathcal A\mathbf x\ge-\gamma\|\mathbf x\|_E^2$.
Assume $1/2\le\theta\le1$ and $\theta\gamma\Delta t<1$, and set
\begin{equation*}
\vartheta=\frac{1+(1-\theta)\gamma\Delta t}{1-\theta\gamma\Delta t}.
\end{equation*}
Write the state delivered before the final rounding as
$\widehat{\bm\Psi}^{n+1}$ and the carried state as
$\widetilde{\bm\Psi}^{n+1}=\widehat{\bm\Psi}^{n+1}+\mathbf r^n+\mathbf p^n$.
The three contributions to the error of a step are:

\begin{itemize}
\item $\widehat{\bm{\rho}}^n = \mathbf{K}\widehat{\bm{\Psi}}^{n+1} - \mathbf{B}\widetilde{\bm{\Psi}}^{n} - \mathbf{q}^n$,
the residual of the delivered state, which collects the compression of
the operators and of the source, the rounding of the right-hand side
and the residual of the linear solve,
\item $\mathbf{r}^n$, the rounding of the delivered state,
\item $\mathbf{p}^n$, the change made by the group balance projection,
zero without it.
\end{itemize}
The residual of \eqref{eq:step-defect} is
$\bm{\rho}^n = \widehat{\bm{\rho}}^n + \mathbf{K}(\mathbf{r}^n + \mathbf{p}^n)$,
so the split separates what the solver sees from what is done to its
answer afterwards. This first lemma bounds one application of $\mathbf{K}^{-1}$ and one
exact step of the method in the energy norm.

\begin{lemma}
\label{lem:onestep}
Let $1/2 \le \theta \le 1$ and $\theta\gamma\Delta t < 1$, with $\gamma$
satisfying the energy inequality above. Then
\begin{enumerate}
\renewcommand{\theenumi}{(\roman{enumi})}\renewcommand{\labelenumi}{\theenumi}
\item\label{lem:onestep-i} $\mathbf{K}$ is invertible and
$\|\mathbf{K}^{-1}\mathbf{y}\|_E \le \|\mathbf{V}\mathbf{y}\|_E/(1 - \theta\,\gamma\Delta t)$
for every $\mathbf{y}$,
\item\label{lem:onestep-ii} $\|\mathbf{K}^{-1}\mathbf{B}\,\mathbf{x}\|_E \le \vartheta\,\|\mathbf{x}\|_E$
for every $\mathbf{x}$.
\end{enumerate}
\end{lemma}

\begin{proof}
$\mathbf{W}$ is symmetric positive definite and commutes with the
positive diagonal matrix $\mathbf{V}$, so
$\langle\mathbf{x},\mathbf{y}\rangle_E = \mathbf{x}^{\top}\mathbf{W}\mathbf{V}^{-1}\mathbf{y}$
is an inner product with norm $\|\cdot\|_E$, and
$\mathbf{x}^{\top}\mathbf{W}\mathbf{y} = \langle\mathbf{x},\mathbf{V}\mathbf{y}\rangle_E$.
The assumed energy inequality gives
$\mathbf{x}^{\top}\mathbf{W}\mathcal{A}\mathbf{x} \ge -\gamma\|\mathbf{x}\|_E^2$
for every $\mathbf{x}$. We write $s = \gamma\Delta t$.

\ref{lem:onestep-i} Both invertibility and the inverse bound follow by
testing the equation

\begin{align*}
\mathbf{K}\mathbf{z} &= \mathbf{V}^{-1}\mathbf{z} + \theta\Delta t\,\mathcal{A}\mathbf{z} = \mathbf{y} .
\end{align*}
We multiply both sides by $\mathbf{z}^{\top}\mathbf{W}$ to obtain
\begin{align*}
\mathbf{z}^{\top}\mathbf{W}\mathbf{V}^{-1}\mathbf{z} + \theta\Delta t\,\mathbf{z}^{\top}\mathbf{W}\mathcal{A}\mathbf{z}
&= \mathbf{z}^{\top}\mathbf{W}\mathbf{y} , \\
\|\mathbf{z}\|_E^2 + \theta\Delta t\,\mathbf{z}^{\top}\mathbf{W}\mathcal{A}\mathbf{z}
&= \langle\mathbf{z},\mathbf{V}\mathbf{y}\rangle_E ,
\end{align*}
where the second line is the first one written with the energy norm and
the inner product, by $\|\mathbf{z}\|_E^2 = \mathbf{z}^{\top}\mathbf{W}\mathbf{V}^{-1}\mathbf{z}$
and $\mathbf{z}^{\top}\mathbf{W}\mathbf{y} = \langle\mathbf{z},\mathbf{V}\mathbf{y}\rangle_E$.
The energy estimate
$\mathbf{z}^{\top}\mathbf{W}\mathcal{A}\mathbf{z} \ge -\gamma\|\mathbf{z}\|_E^2$
bounds the left side from below, and the Cauchy--Schwarz inequality
bounds the right side from above, so that
\begin{equation*}
(1 - \theta s)\,\|\mathbf{z}\|_E^2
\;\le\; \|\mathbf{z}\|_E^2 + \theta\Delta t\,\mathbf{z}^{\top}\mathbf{W}\mathcal{A}\mathbf{z}
\;=\; \langle\mathbf{z},\mathbf{V}\mathbf{y}\rangle_E
\;\le\; \|\mathbf{z}\|_E\,\|\mathbf{V}\mathbf{y}\|_E .
\end{equation*}
For $\mathbf{y} = 0$ the right side vanishes while the left side is
positive unless $\mathbf{z} = 0$, since $1 - \theta s > 0$. Hence
$\mathbf{z} = 0$ and $\mathbf{K}$ is invertible. For general $\mathbf{y}$
dividing by $\|\mathbf{z}\|_E$ gives
$(1 - \theta s)\|\mathbf{z}\|_E \le \|\mathbf{V}\mathbf{y}\|_E$, the bound.

\ref{lem:onestep-ii} Let $\mathbf z=\mathbf K^{-1}\mathbf B\mathbf x$,
$a=\|\mathbf z\|_E$ and $b=\|\mathbf x\|_E$. Define the
$\theta$-average $\bar{\mathbf z}=\theta\mathbf z+(1-\theta)\mathbf x$.
The step equation becomes

\begin{align*}
\mathbf{V}^{-1}\mathbf{z} + \theta\Delta t\,\mathcal{A}\mathbf{z}
&= \mathbf{V}^{-1}\mathbf{x} - (1-\theta)\Delta t\,\mathcal{A}\mathbf{x} , \\
\mathbf{V}^{-1}(\mathbf{z} - \mathbf{x}) + \Delta t\,\mathcal{A}\bigl(\theta\mathbf{z} + (1-\theta)\mathbf{x}\bigr) &= 0 , \\
\mathbf{V}^{-1}(\mathbf{z} - \mathbf{x}) + \Delta t\,\mathcal{A}\bar{\mathbf{z}} &= 0 .
\end{align*}
We multiply both sides by $\bar{\mathbf{z}}^{\top}\mathbf{W}$ and write
the first term with the inner product,
\begin{align*}
\bar{\mathbf{z}}^{\top}\mathbf{W}\mathbf{V}^{-1}(\mathbf{z} - \mathbf{x})
+ \Delta t\,\bar{\mathbf{z}}^{\top}\mathbf{W}\mathcal{A}\bar{\mathbf{z}} &= 0 , \\
\langle\bar{\mathbf{z}},\mathbf{z} - \mathbf{x}\rangle_E
&= -\Delta t\,\bar{\mathbf{z}}^{\top}\mathbf{W}\mathcal{A}\bar{\mathbf{z}} .
\end{align*}
The right side is bounded from above by the energy estimate,
$-\bar{\mathbf{z}}^{\top}\mathbf{W}\mathcal{A}\bar{\mathbf{z}} \le \gamma\|\bar{\mathbf{z}}\|_E^2$,
and then by the triangle inequality for $\bar{\mathbf{z}}$,
\begin{align*}
\langle\bar{\mathbf{z}},\mathbf{z} - \mathbf{x}\rangle_E
\;\le\; s\,\|\bar{\mathbf{z}}\|_E^2
\;\le\; s\,\bigl(\theta a + (1-\theta) b\bigr)^2 .
\end{align*}
The left side is bounded from below by expanding the inner product,
with $\langle\mathbf{z},\mathbf{z}\rangle_E = a^2$ and
$\langle\mathbf{x},\mathbf{x}\rangle_E = b^2$,
\begin{align*}
\langle\bar{\mathbf{z}},\mathbf{z} - \mathbf{x}\rangle_E
&= \theta\langle\mathbf{z},\mathbf{z}\rangle_E - \theta\langle\mathbf{z},\mathbf{x}\rangle_E
 + (1-\theta)\langle\mathbf{x},\mathbf{z}\rangle_E - (1-\theta)\langle\mathbf{x},\mathbf{x}\rangle_E \\
&= \theta a^2 - (1-\theta) b^2 + (1-2\theta)\langle\mathbf{z},\mathbf{x}\rangle_E \\
&\ge \theta a^2 - (1-\theta) b^2 + (1-2\theta)\,ab
 = (a - b)\bigl(\theta a + (1-\theta) b\bigr) ,
\end{align*}
where the inequality uses $1 - 2\theta \le 0$ together with
$\langle\mathbf{z},\mathbf{x}\rangle_E \le ab$, the Cauchy--Schwarz
inequality, so the mixed term is smallest when the two states are
aligned. Combining the two bounds,
\begin{align*}
(a - b)\bigl(\theta a + (1-\theta) b\bigr) \;\le\; s\,\bigl(\theta a + (1-\theta) b\bigr)^2 .
\end{align*}
If $\theta a + (1-\theta) b = 0$ then $a = 0$, since $\theta > 0$, and
the bound holds. Otherwise we divide by $\theta a + (1-\theta) b > 0$
and collect the terms in $a$ on the left,
\begin{align*}
a - b &\le s\,\bigl(\theta a + (1-\theta) b\bigr) , \\
(1 - \theta s)\,a &\le \bigl(1 + (1-\theta) s\bigr)\,b , \\
a &\le \vartheta\, b ,
\end{align*}
where the last line divides by $1 - \theta s > 0$ and is the definition
of $\vartheta$.
\end{proof}

Lemma~\ref{lem:onestep} bounds the inverse step matrix by
$\|\mathbf K^{-1}\mathbf y\|_E\le
\|\mathbf V\mathbf y\|_E/(1-\theta\gamma\Delta t)$ and the exact
step propagator by $\vartheta$. For Crank--Nicolson, $\theta=1/2$,
the proof is the usual energy argument. Background on energy estimates
and stability of $\theta$-methods can be found
in~\cite[Section~5.4]{quarteroni2017numerical}
and~\cite[Chapter~I]{hundsdorfer2003numerical}.

The second lemma writes the error of a computed step in terms of the
three error sources of the step.

\begin{lemma}
\label{lem:recursion}
With the notation of Table~\ref{tab:notation}, the errors of two
consecutive computed states satisfy
\begin{equation}
\mathbf{e}^{n+1} = \mathbf{K}^{-1}\mathbf{B}\,\mathbf{e}^{n}
+ \mathbf{K}^{-1}\widehat{\bm{\rho}}^{n} + \mathbf{r}^n + \mathbf{p}^n ,
\label{eq:error-recursion}
\end{equation}
and, with $\mathbf{e}^0 = 0$,
\begin{equation}
\mathbf{e}^{N_t} = \sum_{n=0}^{N_t-1}
\bigl(\mathbf{K}^{-1}\mathbf{B}\bigr)^{N_t-1-n}
\bigl(\mathbf{K}^{-1}\widehat{\bm{\rho}}^{n} + \mathbf{r}^n + \mathbf{p}^n\bigr) .
\label{eq:error-sum-split}
\end{equation}
\end{lemma}

\begin{proof}
Compare the delivered state with the exact step. By the definition of $\widehat{\bm{\rho}}^n$ the
delivered state satisfies
$\mathbf{K}\widehat{\bm{\Psi}}^{n+1} = \mathbf{B}\widetilde{\bm{\Psi}}^{n} + \mathbf{q}^n + \widehat{\bm{\rho}}^n$,
and the exact step is
$\mathbf{K}\bm{\Psi}^{n+1} = \mathbf{B}\bm{\Psi}^{n} + \mathbf{q}^n$.
Subtracting the second from the first, the source cancels and
\begin{align*}
\mathbf{K}\bigl(\widehat{\bm{\Psi}}^{n+1} - \bm{\Psi}^{n+1}\bigr)
&= \mathbf{B}\bigl(\widetilde{\bm{\Psi}}^{n} - \bm{\Psi}^{n}\bigr) + \widehat{\bm{\rho}}^n
= \mathbf{B}\,\mathbf{e}^{n} + \widehat{\bm{\rho}}^n , \\
\widehat{\bm{\Psi}}^{n+1} - \bm{\Psi}^{n+1}
&= \mathbf{K}^{-1}\mathbf{B}\,\mathbf{e}^{n} + \mathbf{K}^{-1}\widehat{\bm{\rho}}^n ,
\end{align*}
where the second line uses that $\mathbf{K}$ is invertible by
Lemma~\ref{lem:onestep}. Adding $\mathbf{r}^n + \mathbf{p}^n$ to both
sides turns the left side into
$\widetilde{\bm{\Psi}}^{n+1} - \bm{\Psi}^{n+1} = \mathbf{e}^{n+1}$,
by $\widetilde{\bm{\Psi}}^{n+1} = \widehat{\bm{\Psi}}^{n+1} + \mathbf{r}^n + \mathbf{p}^n$,
and the right side into that of \eqref{eq:error-recursion}. The sum
\eqref{eq:error-sum-split} follows by induction on $N_t$. For
$N_t = 1$ it is \eqref{eq:error-recursion} with $\mathbf{e}^0 = 0$.
If it holds for $N_t - 1$, then \eqref{eq:error-recursion} with
$n = N_t - 1$ gives
\begin{align*}
\mathbf{e}^{N_t}
&= \mathbf{K}^{-1}\mathbf{B}\,\mathbf{e}^{N_t-1}
 + \mathbf{K}^{-1}\widehat{\bm{\rho}}^{N_t-1} + \mathbf{r}^{N_t-1} + \mathbf{p}^{N_t-1} \\
&= \sum_{n=0}^{N_t-2}\bigl(\mathbf{K}^{-1}\mathbf{B}\bigr)^{N_t-1-n}
 \bigl(\mathbf{K}^{-1}\widehat{\bm{\rho}}^{n} + \mathbf{r}^n + \mathbf{p}^n\bigr)
 + \mathbf{K}^{-1}\widehat{\bm{\rho}}^{N_t-1} + \mathbf{r}^{N_t-1} + \mathbf{p}^{N_t-1} ,
\end{align*}
which is \eqref{eq:error-sum-split}, the last three terms being its
summand for $n = N_t - 1$.
\end{proof}

Substituting
$\bm\rho^n=\widehat{\bm\rho}^n+\mathbf K(\mathbf r^n+\mathbf p^n)$
into the error recursion separates the three contributions. The
pre-rounding defect is acted on by $\mathbf K^{-1}$, whereas rounding
and balance correction enter directly as changes of the state. All
three contributions then propagate through powers of
$\mathbf K^{-1}\mathbf B$. Applying Lemma~\ref{lem:onestep} to those
powers gives the following bound.

\begin{proposition}
\label{prop:global}
Let $1/2 \le \theta \le 1$, $\theta\gamma\Delta t < 1$ and
$\mathbf{e}^0 = 0$. Then
\begin{equation}
\|\mathbf{e}^{N_t}\|_E \le \sum_{n=0}^{N_t-1} \vartheta^{\,N_t-1-n}
\Bigl(\|\mathbf{r}^n\|_E + \|\mathbf{p}^n\|_E
+ \frac{\|\mathbf{V}\widehat{\bm{\rho}}^n\|_E}{1 - \theta\,\gamma\Delta t}\Bigr) .
\label{eq:global-bound}
\end{equation}
\end{proposition}

\begin{proof}
To bound each summand of \eqref{eq:error-sum-split}, set $\mathbf{d}^n = \mathbf{K}^{-1}\widehat{\bm{\rho}}^n + \mathbf{r}^n + \mathbf{p}^n$
for the summand before the powers of $\mathbf{K}^{-1}\mathbf{B}$ are
applied. Taking the norm of \eqref{eq:error-sum-split} and using the
triangle inequality,
\begin{align*}
\|\mathbf{e}^{N_t}\|_E
&\le \sum_{n=0}^{N_t-1} \bigl\|\bigl(\mathbf{K}^{-1}\mathbf{B}\bigr)^{N_t-1-n}\mathbf{d}^n\bigr\|_E
\le \sum_{n=0}^{N_t-1} \vartheta^{\,N_t-1-n}\,\|\mathbf{d}^n\|_E ,
\end{align*}
where the second inequality applies Lemma~\ref{lem:onestep}\ref{lem:onestep-ii}
once for each of the $N_t-1-n$ factors $\mathbf{K}^{-1}\mathbf{B}$.
For the summand itself, the triangle inequality and
Lemma~\ref{lem:onestep}\ref{lem:onestep-i} with $\mathbf{y} = \widehat{\bm{\rho}}^n$ give
\begin{align*}
\|\mathbf{d}^n\|_E
&\le \|\mathbf{K}^{-1}\widehat{\bm{\rho}}^n\|_E + \|\mathbf{r}^n\|_E + \|\mathbf{p}^n\|_E
\le \frac{\|\mathbf{V}\widehat{\bm{\rho}}^n\|_E}{1 - \theta\,\gamma\Delta t} + \|\mathbf{r}^n\|_E + \|\mathbf{p}^n\|_E .
\end{align*}
Inserting the second bound into the first gives \eqref{eq:global-bound}.
\end{proof}

Each contribution from step $n$ is amplified by at most
$\vartheta^{N_t-1-n}$ before the final time. If $\gamma=0$, then
$\vartheta=1$ and the estimate reduces to the sum of the step
contributions. Non-multiplication alone does not imply $\gamma=0$ in
this quadratic norm; that condition must be checked from the weighted
symmetric collision matrix. For $\gamma>0$, the bound includes possible
finite-time growth. The estimate concerns error relative to the exact
discrete $\theta$-method trajectory, not the continuous transport
solution. For TT-BUG, $\widehat{\bm\rho}^n$ must also include the basis
evolution and projected-integration defects; the local GMRES tolerance
and the rounding budget do not bound this full defect by themselves.

A known issue with rounding is that it does not conserve neutrons.
Let $J$ be the functional that integrates a state over the cube, the
directions and the groups, and $\mathcal{N}$ the neutron count,
\begin{equation*}
J(\mathbf{x}) = \sum_{g=1}^{G}\sum_{m=1}^{M} w_m \int_{\mathcal{D}} x_{g,m}(\rvec)\,dV ,
\qquad
\mathcal{N}(\mathbf{x}) = J(\mathbf{V}^{-1}\mathbf{x})
= \sum_{g=1}^{G}\frac{1}{v_g}\sum_{m=1}^{M} w_m \int_{\mathcal{D}} x_{g,m}(\rvec)\,dV ,
\end{equation*}
where $x_{g,m}$ is the function of the scheme with the coefficients of
$\mathbf{x}$ in group $g$ and direction $m$. $J$ carries no weight in
the group, and $\mathcal{N}$ divides each group by its speed, so that it counts
neutrons rather than flux. The exact discrete step conserves this count
in the sense of the balance equation. Since $\mathbf{K}$ and $\mathbf{B}$
carry the time term $\mathbf{V}^{-1}$, applying $J$ to the step gives
$J(\mathbf{K}\bm{\Psi}^{n+1} - \mathbf{B}\bm{\Psi}^{n} - \mathbf{q}^n) = 0$,
which says that the count changes through the source, leakage and net collision
production or loss. For non-multiplying scattering without fission,
the collision contribution reduces to absorption; the synthetic
active benchmark has a net production contribution. For
the computed states the same functional need not vanish.
$J(\bm\rho^n)$ is the algebraic balance defect of the uncompressed
discrete step. The implementation evaluates the analogous quantity
for the compressed scaled system, defined next. Neither quantity is
an accumulated physical-error bound.

\emph{Balance of every group.} The final carried state has compressed,
velocity-scaled residual
\begin{equation*}
\bm\rho_s^n=\widetilde{\mathbf K}_s\widetilde{\bm\Psi}^{n+1}
 -\widetilde{\mathbf B}_s\widetilde{\bm\Psi}^{n}
 -\widetilde{\mathbf q}_s^n .
\end{equation*}
The matrices and source compress $\mathbf V\mathbf K$, $\mathbf V\mathbf B$
and $\mathbf V\mathbf q^n$. Thus $\bm\rho_s^n$ differs from
$\mathbf V\bm\rho^n$ by compression errors. The uncompressed residual
$\bm\rho^n$ remains the quantity used in the error analysis.
Let $\mathcal N_g$ restrict the neutron-count functional to group $g$,
and let $\mathbf E_{g'}$ retain only the coefficients of group $g'$.
Write $\mathbf y=\widehat{\bm\Psi}^{n+1}+\mathbf r^n$ for the rounded
state before balance correction and define
\begin{equation*}
\bm\rho_{s,\mathrm{pre}}^n=\widetilde{\mathbf K}_s\mathbf y
 -\widetilde{\mathbf B}_s\widetilde{\bm\Psi}^n
 -\widetilde{\mathbf q}_s^n .
\end{equation*}
The correction coefficients solve
\begin{equation*}
\sum_{g'=1}^G\mathcal N_g(\widetilde{\mathbf K}_s\mathbf E_{g'}\mathbf y)
\alpha_{g'}=-\mathcal N_g(\bm\rho_{s,\mathrm{pre}}^n),
\qquad g=1,\ldots,G,
\end{equation*}
and the carried state is
\begin{equation*}
\widetilde{\bm\Psi}^{n+1}
=\bigl(\mathbf I_n\otimes\mathbf I_n\otimes\mathbf I_n\otimes\mathbf I_M
\otimes\operatorname{diag}(1+\alpha_1,\ldots,1+\alpha_G)\bigr)\mathbf y.
\end{equation*}
The change $\mathbf p^n=\widetilde{\bm\Psi}^{n+1}-\mathbf y$ is the
projection term in Proposition~\ref{prop:global}. If the correction
equations are satisfied, linearity gives $\mathcal N_g(\bm\rho_s^n)=0$
for each group, up to arithmetic error. This restores group-integrated
balance for the compressed step. Source, scattering and fission gains
balance leakage and total collision removal; the total cross section
already includes out-scattering. The uncompressed residual is not set
to zero.

The exact cancellation assumes that the correction system is solvable.
If a group slice is identically zero, its scaling cannot change the
state; the implementation omits the corresponding column and uses a
least-squares solve on the remaining columns. Exact cancellation then
requires compatibility of that reduced system. Ill-conditioning can
amplify the correction, so algebraic balance should be checked after
the projection. This correction neither enforces cell-local balance
nor guarantees nonnegative flux.

\section{Numerical examples}
\label{sec:examples}

We explore two examples, the first is a manufactured solution with a
known exact flux, the second is a burst of neutrons released in a cube.
We use the first as a verification test, where we verify the well-known
convergence rates of the AMEn-based Crank--Nicolson DG scheme and
compare measured operator ranks with their analytical bounds. For the second example we explore the two time integrators
numerically and compare them against a standard full-tensor solver and
against OpenSn~\cite{andrs2026opensn}. The machine and the software
behind every time and memory reported below are described in
\ref{app:machine}.

\FloatBarrier
\subsection{Manufactured solution}
\label{sec:manufactured}

On the cube $\mathcal{D} = (0,1)^3$ we solve the multigroup system
\eqref{eq:multigroup} with $G = 256$ groups, isotropic scattering and no
fission, with a source chosen so that
\begin{equation}
\psi_g(\rvec,\Om,t) = \beta_g\, t^2 e^{-t}\,
(1+\Omega_x^2)(1+\Omega_y^2)(1+\Omega_z^2)\,
\sin(\pi x)\sin(\pi y)\sin(\pi z) ,
\qquad \beta_g = e^{-3(g-1)/(G-1)} ,
\label{eq:mms-exact}
\end{equation}
is the exact solution. The group data are synthetic, with $v_g = 1$,
$\sigma^{g} = 1 + (g-1)/(2(G-1))$ and downscattering over four groups,
$\sigma_s^{g'\to g} = 0.6\,\pi\, e^{-(g-g')}$ for $0 \le g - g' \le 4$.
The exact flux vanishes on the boundary of the cube and at $t = 0$, so
the inflow data are those of vacuum and the initial state is zero. We use the product quadrature of
Section~\ref{sec:dg-directions} with eight polar points and sixteen
azimuths, $M = 128$ directions, and the Crank--Nicolson step to $T = 1$.
The error at $T = 1$ is
\begin{equation*}
\mathcal{E} = \Bigl( \sum_{g=1}^{G} \sum_{m=1}^{M} 4\pi w_m \int_{\mathcal{D}}
\bigl( \psi_{g,m}(\rvec,T) - \Pi\psi_g(\rvec,\Om_m,T) \bigr)^2 \, dV \Bigr)^{1/2},
\end{equation*}
where $\psi_{g,m}$ is the computed flux and $\Pi$ the cell-by-cell
projection of Section~\ref{sec:dg-space} onto the finite-element space,
whose coefficients are products of the one-dimensional coefficients
\eqref{eq:graphite-proj}, since the exact flux is a product of functions
of $x$, $y$ and $z$. In terms of the coefficient vectors, with $\widetilde{\bm{\Psi}}^{N_t}$ the
computed state after the $N_t$ steps to $T$ and $\bm{\psi}_\Pi$ the coefficients of
$\Pi\psi(\cdot,T)$,
\begin{equation*}
\mathcal{E} = \bigl( 4\pi\, (\widetilde{\bm{\Psi}}^{N_t} - \bm{\psi}_\Pi)^{\top}\mathbf{W}\,
(\widetilde{\bm{\Psi}}^{N_t} - \bm{\psi}_\Pi) \bigr)^{1/2},
\qquad
\mathbf{W} = \mathbf{M}_h\otimes\mathbf{M}_h\otimes\mathbf{M}_h
\otimes\operatorname{diag}(\mathbf{w})\otimes\mathbf{I}_G ,
\end{equation*}
the matrix of Lemma~\ref{lem:energy}. With $v_g = 1$ this is $\sqrt{4\pi}$
times the norm $\|\cdot\|_E$ of Section~\ref{sec:error-propagation}. We
measure the error this way here because the exact solution is known.
The bursts of Section~\ref{sec:bursts} have no exact solution and
measure their error against a reference solution instead. In
Table~\ref{tab:mms} we provide a mesh refinement study in space and
time. One can see that our upwind DG scheme exhibits the expected
convergence rates in both the TT and the hybrid TT/QTT format, while
the full tensor is listed for its size only. This shows that our
framework is flexible, but in this specific case the hybrid form
performs significantly better in terms of final memory and time. At
the finest mesh the final state takes 13\,kB against 35\,TB for the
full tensor.

\begin{table}[!htbp]
\centering
\caption{Manufactured solution, refinement in space and time with
$\Delta t = h/10$, $\varepsilon = 10^{-6}$, ranks capped at $16$.
Columns: $N$ and $\Delta t$, the error $\mathcal{E}$ at $T = 1$ and its
rate between consecutive rows, the entries and the memory of the full
tensor, and for the tensor train with one core per index (TT) and for
the hybrid form the largest rank of the state, the memory of the final
state and the wall time. The errors are those of the hybrid form, the
tensor train differs by at most $0.6\%$.}
\label{tab:mms}
\vspace{6pt}
\footnotesize
\setlength{\tabcolsep}{3pt}
\begin{tabular}{rccc@{\hskip 6pt}rr@{\hskip 6pt}rrr@{\hskip 6pt}rrr}
\toprule
& & & & \multicolumn{2}{c}{\textbf{Full tensor}} & \multicolumn{3}{c}{\textbf{TT}} & \multicolumn{3}{c}{\textbf{Hybrid}} \\
\cmidrule(lr){5-6}\cmidrule(lr){7-9}\cmidrule(lr){10-12}
$\bm{N}$ & $\bm{\Delta t}$ & $\bm{\mathcal{E}}$ & \textbf{Rate} & \textbf{Entries} & \textbf{Memory} & \textbf{Rank} & \textbf{Memory} & \textbf{Time [s]} & \textbf{Rank} & \textbf{Memory} & \textbf{Time [s]} \\
\midrule
16 & $1/160$ & $1.950\times 10^{-2}$ & -- & $1.1\times 10^{9}$ & 8.6\,GB & 16 & 146\,kB & 168 & 16 & 97\,kB & 134 \\
32 & $1/320$ & $5.012\times 10^{-3}$ & 1.96 & $8.6\times 10^{9}$ & 69\,GB & 14 & 111\,kB & 124 & 16 & 69\,kB & 89 \\
64 & $1/640$ & $1.271\times 10^{-3}$ & 1.98 & $6.9\times 10^{10}$ & 550\,GB & 4 & 27\,kB & 50 & 12 & 25\,kB & 53 \\
128 & $1/1280$ & $3.199\times 10^{-4}$ & 1.99 & $5.5\times 10^{11}$ & 4.4\,TB & 4 & 47\,kB & 227 & 7 & 13\,kB & 56 \\
256 & $1/2560$ & $8.000\times 10^{-5}$ & 2.00 & $4.4\times 10^{12}$ & 35\,TB & 4 & 88\,kB & 1030 & 7 & 13\,kB & 101 \\
\bottomrule
\end{tabular}
\end{table}

Table~\ref{tab:mms-dt} fixes the spatial discretization at $N = 64$
and refines the step. The error settles at the spatial error of this
mesh from $\Delta t = h/5$ on. Along the path of Table~\ref{tab:mms}
the time error is therefore invisible, and the rate of two there is the
spatial order. The order in time is read from the difference between
the final states of consecutive steps. The ranks are
larger than in Table~\ref{tab:mms}, since the tolerance is smaller and
the time error of a large step is itself a field that the state has to
carry.

\begin{table}[!htbp]
\centering
\caption{Manufactured solution, refinement in time on the fixed mesh
$N = 64$, hybrid form, $\varepsilon = 10^{-8}$, no cap on the ranks.
Columns: $\Delta t$, the error $\mathcal{E}$ at $T = 1$, the difference
between the final states for the steps $\Delta t$ and $\Delta t/2$
relative to the norm of the state, the order read from that difference,
and the largest rank of the state.}
\label{tab:mms-dt}
\vspace{6pt}
\footnotesize
\begin{tabular}{rcccr}
\toprule
$\bm{\Delta t}$ & $\bm{\mathcal{E}}$ & $\bm{\|\bm{\Psi}_{\Delta t}-\bm{\Psi}_{\Delta t/2}\|/\|\bm{\Psi}\|}$ & \textbf{Order} & \textbf{Rank} \\
\midrule
$1/40$ & $3.807\times 10^{-3}$ & $4.02\times 10^{-4}$ & -- & 387 \\
$1/80$ & $1.565\times 10^{-3}$ & $1.01\times 10^{-4}$ & 2.00 & 236 \\
$1/160$ & $1.294\times 10^{-3}$ & $2.51\times 10^{-5}$ & 2.00 & 126 \\
$1/320$ & $1.273\times 10^{-3}$ & $6.30\times 10^{-6}$ & 2.00 & 77 \\
$1/640$ & $1.271\times 10^{-3}$ & $1.75\times 10^{-6}$ & 1.85 & 58 \\
$1/1280$ & $1.270\times 10^{-3}$ & -- & -- & 45 \\
\bottomrule
\end{tabular}
\end{table}

Table~\ref{tab:mms-eps} fixes the mesh at $N = 64$ and the step at
$\Delta t = 1/640$ and varies the tolerance from $\varepsilon = 10^{-2}$
to $10^{-8}$, with the operators and the source compressed at
$10^{-10}$, so that $\varepsilon$ acts on the state and the right-hand
sides only. From $\varepsilon = 10^{-5}$ on the error is the spatial
error of the mesh, and a smaller tolerance only raises the rank and the
memory of the state, as for any low-rank approximation. Below the
threshold of about $10^{-5}$ a smaller tolerance increases the cost but
cannot improve the error, which is set by the mesh in space and time.

\begin{table}[!htbp]
\centering
\caption{Manufactured solution on the fixed mesh $N = 64$ with
$\Delta t = 1/640$, hybrid form, operators and source compressed at
$10^{-10}$. Columns: the tolerance $\varepsilon$ of the state and the
right-hand sides, the error $\mathcal{E}$ at $T = 1$, the largest rank
of the state, the memory of the final state, and the compression, the
ratio of the full tensor ($6.9\times 10^{10}$ entries, $550$\,GB) to
that memory.}
\label{tab:mms-eps}
\vspace{6pt}
\footnotesize
\begin{tabular}{c rrrr}
\toprule
$\bm{\varepsilon}$ & $\bm{\mathcal{E}}$ & \textbf{Rank} & \textbf{Memory} & \textbf{Compression} \\
\midrule
$10^{-2}$ & $5.639\times 10^{0}$ & 9 & 21\,kB & $2.6\times 10^{7}$ \\
$10^{-3}$ & $1.496\times 10^{-2}$ & 6 & 7.6\,kB & $7.2\times 10^{7}$ \\
$10^{-4}$ & $1.710\times 10^{-3}$ & 7 & 11\,kB & $5.1\times 10^{7}$ \\
$10^{-5}$ & $1.269\times 10^{-3}$ & 8 & 15\,kB & $3.7\times 10^{7}$ \\
$10^{-6}$ & $1.271\times 10^{-3}$ & 12 & 25\,kB & $2.2\times 10^{7}$ \\
$10^{-7}$ & $1.271\times 10^{-3}$ & 16 & 76\,kB & $7.3\times 10^{6}$ \\
$10^{-8}$ & $1.271\times 10^{-3}$ & 16 & 115\,kB & $4.8\times 10^{6}$ \\
\bottomrule
\end{tabular}
\end{table}

Table~\ref{tab:mms-ranks} reports the largest measured operator ranks
in the implemented order, group, direction, $x$, $y$, $z$. The plain TT
operator has largest rank $7$ and the rounded hybrid operator also has
largest rank $7$ on every tested mesh. These values lie below the term
count $q=8$ and the order-independent hybrid bound $q+4=12$ of
\ref{app:rank-order}. The finer bounds $7$, $9$ and $11$ in
Figure~\ref{fig:rank-links} concern the different analysis order,
$x$, $y$, $z$, direction, group. Combining time and removal into one
term is already included in $q=8$; further dependencies can reduce an
observed TT-rank below that count. This experiment supports bounded
operator ranks for the present construction. It does not establish an
unavoidable rank-growth disadvantage for another discretization or
material model.

\begin{table}[htbp]
\centering
\caption{Largest measured TT-rank of the manufactured-solution step
matrix over all bonds, in the implemented group--direction--space
order. Operators are rounded at $10^{-12}$. The bounds use $q=8$ and
\ref{app:rank-order}; Figure~\ref{fig:rank-links} uses a different
mode order.}
\label{tab:mms-ranks}
\vspace{6pt}
\begin{tabular}{lrrrrrc}
\toprule
Representation & $N=16$ & $32$ & $64$ & $128$ & $256$ & Bound \\
\midrule
TT, one core per index & 7 & 7 & 7 & 7 & 7 & $q=8$ \\
Hybrid TT/QTT & 7 & 7 & 7 & 7 & 7 & $q+4=12$ \\
\bottomrule
\end{tabular}
\end{table}

\subsection{Active and passive bursts}
\label{sec:bursts}

On the cube $\mathcal{D} = (0,1)^3$ we solve the multigroup system
\eqref{eq:multigroup} with $G = 16$ groups and no fission, driven by an
isotropic source in the first group that is a
Gaussian in space and a single pulse in time,
\begin{equation}
q_g(\rvec,\Om,t) = \delta_{g1}\, b(t)\,
\exp\Bigl(-\frac{|\rvec - \rvec_0|^2}{2 s^2}\Bigr),
\qquad
b(t) =
\begin{cases}
\sin^2(\pi t/\tau), & 0 \le t < \tau, \\
0, & t \ge \tau,
\end{cases}
\label{eq:burst}
\end{equation}
with $\tau = 0.1$ and $\rvec_0 = (0.3,\,0.4,\,0.55)$, vacuum inflow
and zero initial state. The pulse $b(t)$ rises from zero and returns to
zero at $t = \tau$, after which the source is off. The group data are
$v_g = 1$, $\sigma^{g} = 1 + (g-1)/(2(G-1))$ and downscattering over
four groups,
\begin{equation*}
\sigma_s^{g'\to g} = c_{g'}\, e^{-(g-g')} \quad\text{for } 0 \le g - g' \le 4,
\end{equation*}
and zero otherwise. Physical group speeds grow with the square root of
the energy. The equal speeds are a deliberate simplification of these
synthetic data, and the third case below gives the groups a decade of
speeds. We examine two scenarios with isotropic scattering,
a passive and an active burst, which differ in the scale $c_{g'}$ and in
the source width $s$,
\begin{equation*}
c_{g'} = \frac{0.9\,\sigma^{g'}}{\sum_{d=0}^{\min(4,\,G-g')} e^{-d}},
\quad s = 0.1 \quad\text{(passive)},
\qquad\qquad
c_{g'} = 0.6\,\pi, \quad s = 0.06 \quad\text{(active)}.
\end{equation*}
In the passive burst the scattering out of a group is $0.9\,\sigma^{g'}$,
so a tenth of the collisions absorb and the neutron population decays
once the source stops. In the active burst the data are those of
Section~\ref{sec:manufactured} at $G = 16$, a collision in the top group
produces about three outgoing neutrons in the synthetic model.
These coefficients are yield-weighted production data, not
probabilities of ordinary one-neutron scattering. The population
keeps growing after
the source stops. A third case, examined with the active burst in
Section~\ref{sec:nonseparable}, keeps the active data but makes the transfer anisotropic and the groups
unequal in speed. The transfer out of group $g'$ follows a
Henyey--Greenstein law of Legendre order $3$, whose moments are the
powers of an asymmetry $\alpha_{g'}$,
\begin{equation*}
\sigma_{s,\ell}^{g'\to g} = \sigma_s^{g'\to g}\,\alpha_{g'}^{\,\ell},
\qquad \ell = 0,\dots,3,
\end{equation*}
with $\alpha_{g'}$ falling linearly from $0.4$ in the first group to
$0.1$ in the last. These synthetic moments keep the truncated
angular kernel nonnegative; at $\alpha=0.4$ its minimum is $0.152$.
Larger asymmetries can make the order-three truncation negative
near backscattering. The speeds $v_g = 10^{-(g-1)/(G-1)}$ fall by one
decade from the first group to the last. In all cases, we use $\theta = \tfrac12$, that is,
Crank--Nicolson in time, and $\Delta t = h/10$.
In the rest of this section both bursts are solved on the same mesh,
$N = 32$ cells per axis, with $M = 128$ directions and $G = 16$ groups,
so that the state has $5.4\times 10^{8}$ unknowns, and the same settings
serve throughout. The reference is the full tensor with GMRES solved to
$10^{-10}$, whose row in every table gives its own time and memory, on
the machine of \ref{app:machine}.
OpenSn is measured against its own run at tolerance $10^{-10}$.
It uses $24$ MPI ranks on the same machine. The TT runs use one
computational thread. The isotropic active full-array reference uses
one thread; the passive and anisotropic active references use the
multithreaded MATLAB configuration. Tables~\ref{tab:passive} and
\ref{tab:active} identify these execution modes explicitly. Times are
wall seconds without thread or rank multipliers, and memory is peak
process memory in GB, summed over MPI ranks for OpenSn.

The reference rows report the cost of generating a tightly converged
discrete accuracy target. No full-array tolerance sweep at this mesh
is included, so those costs do not establish full-array speedups at
matched error. OpenSn and the DG methods use different spatial
discretizations and different discrete references; equal entries in
the tolerance column do not imply equal physical accuracy. We report
their wall times and memory as configuration-specific measurements.
For the anisotropic active burst, the OpenSn run at tolerance $10^{-2}$
takes $716$\,s on $24$ processors, while TT-BUG at
$\varepsilon=10^{-4}$ takes $14.3$\,s on one processor
(Table~\ref{tab:active}). Thus the wall-time ratio is
$716/14.3\approx50$, and the processor-time cost ratio is
\begin{equation*}
\frac{24\times716}{1\times14.3}\approx1200.
\end{equation*}
Here processor-time cost means the number of allocated processors
multiplied by wall time, not measured CPU consumption. This is a
comparison of the reported configurations, not a speedup at matched
physical or angular-flux error: the displayed accuracy measures are
scalar-flux and neutron-count errors against each code's own discrete
reference.

\FloatBarrier
\subsubsection{Passive burst}
\label{sec:reference}

In Table~\ref{tab:passive}, we examine the passive burst, $160$ steps to
$T = 0.5$, solved with OpenSn and with the two TT methods, AMEn with the
default settings of Section~\ref{sec:dg-procedure} and BUG. For OpenSn,
whose spatial discretization differs from ours, the reference is its own
run. We report two relative errors at $T$, of the
scalar flux $\phi_g = \sum_m w_m \psi_{g,m}$ in the mean-square norm
over the cube summed over the groups, and of the neutron count $\mathcal{N}$ of
Section~\ref{sec:error-propagation},
\begin{equation*}
\mathcal{E}_\phi = \frac{\bigl( \sum_{g} \int_{\mathcal{D}} (\phi_g - \phi_g^{\mathrm{ref}})^2 \, dV \bigr)^{1/2}}
{\bigl( \sum_{g} \int_{\mathcal{D}} (\phi_g^{\mathrm{ref}})^2 \, dV \bigr)^{1/2}},
\qquad
\mathcal{E}_{\mathcal{N}} = \frac{|\mathcal{N}(\widetilde{\bm{\Psi}}^{N_t}) - \mathcal{N}(\bm{\Psi}^{\mathrm{ref}})|}{\mathcal{N}(\bm{\Psi}^{\mathrm{ref}})}.
\end{equation*}
In terms of the coefficient vectors, the scalar flux is
$\bm{\phi} = (\mathbf{I}_n\otimes\mathbf{I}_n\otimes\mathbf{I}_n\otimes\mathbf{w}^{\top}\otimes\mathbf{I}_G)\,\widetilde{\bm{\Psi}}^{N_t}$,
its squared norm is
$\bm{\phi}^{\top}(\mathbf{M}_h\otimes\mathbf{M}_h\otimes\mathbf{M}_h\otimes\mathbf{I}_G)\,\bm{\phi}$,
and the count is the linear functional
\begin{equation*}
\mathcal{N}(\bm{\Psi}) = (\mathbf{m}\otimes\mathbf{m}\otimes\mathbf{m}\otimes\mathbf{w}\otimes\mathbf{V}_G^{-1}\mathbf{1})^{\top}\,\bm{\Psi},
\end{equation*}
where $\mathbf{m}$ holds the integrals of the $n$ basis functions of one
axis over their cells, $h$ for the mean and $0$ for the slope, and
$\mathbf{V}_G^{-1}\mathbf{1}$ the inverse speeds of the groups. The first
error says where the neutrons are, the second whether the run has kept
the right number of them.

Table~\ref{tab:passive} shows a substantial reduction in cost relative
to the tested AMEn baseline. At $\varepsilon=10^{-6}$, BUG gives scalar-flux
error $2.3\times10^{-5}$ in $858$\,s, compared with $2.4\times10^{-5}$
in $16725$\,s for AMEn, approximately $19.5$ times faster at similar
error. A basic BUG step uses one sweep of small projected solves;
additional passes occur when the onset or saturation rules trigger.
The AMEn settings allow up to $30$ sweeps. The difference in sweep
structure explains the observed advantage for these implementations;
it is not a representation-only complexity guarantee.

BUG also applies the group balance correction of
Section~\ref{sec:error-propagation}. With all other settings unchanged,
omitting the correction raises its final count error to
$1.3\times10^{-4}$ at $\varepsilon=10^{-4}$ and $1.6\times10^{-5}$ at
$10^{-6}$, compared with $1.6\times10^{-5}$ and $1.5\times10^{-8}$
with correction. The corresponding reductions are approximately eight
and one thousand. The additional run at $\varepsilon=5\times10^{-7}$
reduces scalar-flux error to $1.7\times10^{-5}$, although the count
error is not monotone in the tolerance. A rounding tolerance is not a
bound on either final error.

The full-array reference uses $113$\,GB, compared with
$9.2$--$17.7$\,GB for AMEn and $7.2$--$9.9$\,GB for BUG.
These are process-memory measurements for the stated implementations,
including their working storage. OpenSn uses $18.6$--$19.9$\,GB and
$1060$--$1355$\,s in its non-reference runs on $24$ ranks. Its
iteration error falls from $1.4\times10^{-5}$ to $1.0\times10^{-7}$
against its own converged run. After division by the $4\pi$ of its
angular convention, its final count differs from the full-array DG
reference by $1.6\times10^{-3}$. This difference includes the effects
of the distinct discretizations; the table does not establish a
common physical-error ranking of the codes.

\begin{table}[!htbp]
\centering
\caption{Passive burst. \textbf{Solver} denotes the method, \textbf{Tol}
its tolerance, the relative residual of GMRES for the full tensor, the
iteration tolerance of OpenSn and the rounding tolerance $\varepsilon$
of the TT methods, $\bm{\mathcal{E}_\phi}$ and $\bm{\mathcal{E}_{\mathcal{N}}}$
the relative errors of the scalar flux and of the neutron count at $T$,
\textbf{Rank} the largest rank of the state with the peak rank inside
the solver in parentheses, \textbf{Time} the wall time and
\textbf{Memory} the peak memory. \textbf{Exec} identifies one
MATLAB computational thread (1T), multithreaded MATLAB (MT), or
24 MPI ranks (24R). Reference-generation costs are not
matched-accuracy baselines; OpenSn errors use its own reference.}
\label{tab:passive}
\vspace{6pt}
\footnotesize
\setlength{\tabcolsep}{3.2pt}
\begin{tabular}{ll cc c c rr}
\toprule
\textbf{Solver} & \textbf{Tol} & $\bm{\mathcal{E}_\phi}$ & $\bm{\mathcal{E}_{\mathcal{N}}}$ & \textbf{Rank} & \textbf{Exec} & \textbf{Wall [s]} & \textbf{Memory [GB]} \\
\midrule
full tensor & $10^{-10}$ (ref.) & -- & -- & -- & MT & 10809 & 113 \\
\midrule
OpenSn & $10^{-10}$, reference & -- & -- & -- & 24R & 1874 & 19.9 \\
OpenSn & $10^{-2}$ & $1.4\times 10^{-5}$ & $9.3\times 10^{-7}$ & -- & 24R & 1060 & 18.6 \\
OpenSn & $10^{-3}$ & $1.3\times 10^{-5}$ & $1.4\times 10^{-6}$ & -- & 24R & 1096 & 18.6 \\
OpenSn & $10^{-6}$ & $1.0\times 10^{-7}$ & $3.2\times 10^{-8}$ & -- & 24R & 1355 & 19.9 \\
\midrule
AMEn & $\varepsilon = 10^{-4}$ & $2.2\times 10^{-3}$ & $2.9\times 10^{-4}$ & 367 (425) & 1T & 3152 & 9.2 \\
AMEn & $\varepsilon = 10^{-6}$ & $2.4\times 10^{-5}$ & $1.6\times 10^{-5}$ & 772 (828) & 1T & 16725 & 17.7 \\
\midrule
BUG & $\varepsilon = 10^{-4}$ & $2.5\times 10^{-3}$ & $1.6\times 10^{-5}$ & 387 (512) & 1T & 135 & 7.7 \\
BUG & $\varepsilon = 10^{-5}$ & $1.8\times 10^{-4}$ & $3.8\times 10^{-7}$ & 592 (813) & 1T & 337 & 7.2 \\
BUG & $\varepsilon = 10^{-6}$ & $2.3\times 10^{-5}$ & $1.5\times 10^{-8}$ & 870 (1207) & 1T & 858 & 9.3 \\
BUG & $\varepsilon = 5\times 10^{-7}$ & $1.7\times 10^{-5}$ & $2.0\times 10^{-7}$ & 966 (1342) & 1T & 1088 & 9.9 \\
\bottomrule
\end{tabular}
\end{table}

\begin{figure}[!htbp]
\centering
\includegraphics[width=\textwidth]{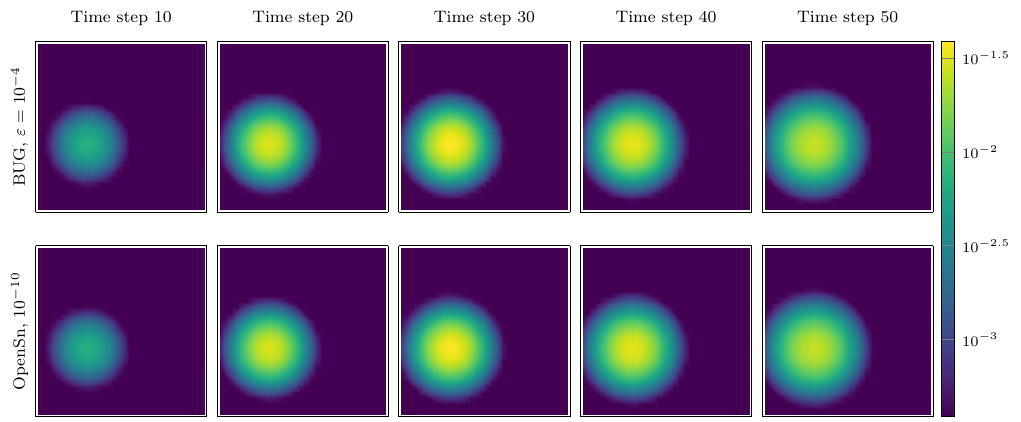}
\caption{Passive burst, $N = 32$, $M = 128$, $G = 16$. Neutron density
summed over the groups in the plane $z = 1/2$ at time steps $10$, $20$,
$30$, $40$ and $50$ of $160$. Top: BUG at $\varepsilon = 10^{-4}$.
Bottom: OpenSn at the tolerance $10^{-10}$. One logarithmic color scale
for all panels, two decades below the largest density, shared with
Figure~\ref{fig:aniso-frames}.}
\label{fig:passive-frames}
\end{figure}

Figure~\ref{fig:passive-frames} follows the passive burst in time in
the plane $z = 1/2$, for BUG at $\varepsilon = 10^{-4}$ and for OpenSn,
at the time steps and on the colour scale of
Figure~\ref{fig:aniso-frames}, so that the two bursts compare side by
side. The pulse builds up while the source is on and spreads after
the source stops at step $32$. The reported maxima from BUG and OpenSn
differ by about two percent at the plotted steps. The passive source
is wider than the active source. The images alone cannot isolate the
effects of source width, anisotropy, group speed and multiplication.

\FloatBarrier

\begin{figure}[!htbp]
\centering
\includegraphics[width=\textwidth]{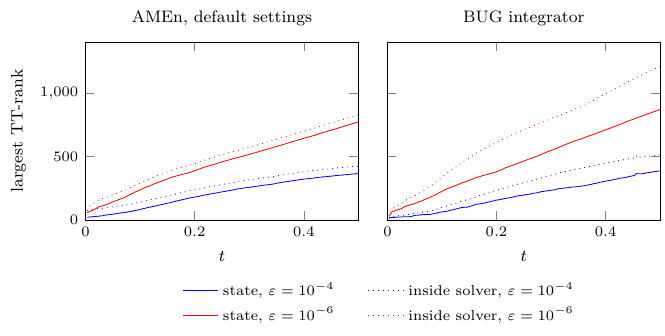}
\caption{Passive burst, largest TT-rank over the run at
$\varepsilon = 10^{-4}$ (blue) and $10^{-6}$ (red), the state after
rounding (solid) and the rank inside the solver (dotted), as in the
tables. Left: AMEn with the default settings, where the rank inside the
solver is that of the iterate it delivers before rounding. Right: the
BUG integrator, where it is that of the augmented bases.}
\label{fig:passive-ranks}
\end{figure}

Figure~\ref{fig:passive-ranks} shows the ranks over the time domain,
the rank of the state after the rounding and the rank inside the
solver, for AMEn and for BUG at $\varepsilon = 10^{-4}$ and $10^{-6}$.
The state ranks grow approximately linearly over the plotted interval
for both integrators in the tested layout. Comparable growth is
therefore not specific to AMEn. These measurements do not establish
minimal approximation ranks or exclude effects of rounding and
accumulated integration error. At $\varepsilon=10^{-6}$, the largest
internal ranks are $828$ for AMEn and $1207$ for BUG; the tighter BUG run
at $\varepsilon=5\times10^{-7}$ reaches $1342$. Thus BUG's lower runtime
does not come from
uniformly smaller intermediate ranks.

\FloatBarrier

\begin{figure}[!htbp]
\centering
\includegraphics[width=\textwidth]{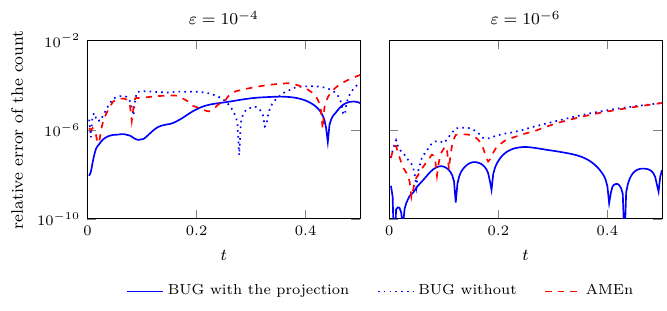}
\caption{Passive burst, relative error of the neutron count against the
reference at every step, for BUG with the group balance projection
(solid) and without it (dotted) and for AMEn, which runs without it
(dashed). Left: $\varepsilon = 10^{-4}$. Right: $\varepsilon = 10^{-6}$.}
\label{fig:passive-balance}
\end{figure}

Figure~\ref{fig:passive-balance} shows the effect of the group balance
projection of Section~\ref{sec:error-propagation} on the error of the
neutron count, for BUG with and without the projection and everything
else unchanged, next to AMEn, which runs without it. The projection
makes the balance of every group exact for the compressed step. What
remains is the count error $\mathcal{E}_{\mathcal{N}}$ against the
reference, which comes from the compression of the operators and from
the error of the flux, which changes the leakage of the later steps. At
$T$ the projection reduces the count error by a factor of eight at
$\varepsilon = 10^{-4}$ and by a factor of a thousand at $10^{-6}$, and
without it both methods show larger count errors. This supports
the usefulness of the algebraic balance correction without
identifying it with exact physical counts or local conservation.

\FloatBarrier
\subsubsection{Active burst}
\label{sec:nonseparable}

Table~\ref{tab:active} compares the active bursts at the common endpoint
$t=0.15625$, corresponding to $50$ steps. The isotropic block uses equal
group speeds; the anisotropic block uses the order-three kernel and
group speeds of Section~\ref{sec:bursts}. In the isotropic case, BUG at
$\varepsilon=5\times10^{-7}$ gives scalar-flux error $1.4\times10^{-6}$
in $102$\,s, compared with AMEn's $1.6\times10^{-6}$ in $1626$\,s,
approximately $15.9$ times faster with slightly smaller error.
In the anisotropic case, the corresponding errors are
$1.8\times10^{-6}$ and $2.3\times10^{-6}$, with times $123$ and
$1840$\,s, approximately a factor of $15$. The full-array reference
costs validate the accuracy targets; they do not measure the cost of
the full-array method at those same errors. The memory advantage is
substantial: BUG uses $6.2$--$7.2$\,GB, against $122$--$169$\,GB for
reference generation.

\begin{table}[!htbp]
\centering
\caption{Active burst, the isotropic block with isotropic scattering and
equal speeds, the anisotropic block with the anisotropic scattering and
group speeds of the third case of Section~\ref{sec:bursts}. Columns as
in Table~\ref{tab:passive}. The reference of the isotropic block ran on
one thread, that of the anisotropic block with multithreaded
MATLAB. Execution modes and reference limitations are as in
Table~\ref{tab:passive}.}
\label{tab:active}
\vspace{6pt}
\footnotesize
\setlength{\tabcolsep}{2.8pt}
\begin{tabular}{l ll cc c c rr}
\toprule
& \textbf{Solver} & \textbf{Tol} & $\bm{\mathcal{E}_\phi}$ & $\bm{\mathcal{E}_{\mathcal{N}}}$ & \textbf{Rank} & \textbf{Exec} & \textbf{Wall [s]} & \textbf{Memory [GB]} \\
\midrule
\multirow{11}{*}{\rotatebox{90}{\textbf{isotropic}}}
& full tensor & $10^{-10}$ (ref.) & -- & -- & -- & 1T & 13467 & 122 \\
& OpenSn & $10^{-10}$, reference & -- & -- & -- & 24R & 618 & 20.0 \\
& OpenSn & $10^{-2}$ & $1.7\times 10^{-5}$ & $3.3\times 10^{-5}$ & -- & 24R & 327 & 18.5 \\
& OpenSn & $10^{-3}$ & $1.6\times 10^{-5}$ & $2.9\times 10^{-5}$ & -- & 24R & 345 & 18.6 \\
& OpenSn & $10^{-6}$ & $3.0\times 10^{-8}$ & $3.8\times 10^{-8}$ & -- & 24R & 483 & 19.9 \\
& AMEn & $\varepsilon = 10^{-4}$ & $2.2\times 10^{-4}$ & $3.2\times 10^{-5}$ & 229 (284) & 1T & 317 & 7.5 \\
& AMEn & $\varepsilon = 10^{-6}$ & $1.6\times 10^{-6}$ & $1.1\times 10^{-6}$ & 500 (552) & 1T & 1626 & 11.1 \\
& BUG & $\varepsilon = 10^{-4}$ & $3.8\times 10^{-4}$ & $1.8\times 10^{-6}$ & 209 (337) & 1T & 16.2 & 6.3 \\
& BUG & $\varepsilon = 10^{-5}$ & $3.0\times 10^{-5}$ & $3.6\times 10^{-7}$ & 352 (542) & 1T & 34.6 & 6.6 \\
& BUG & $\varepsilon = 10^{-6}$ & $2.8\times 10^{-6}$ & $9.0\times 10^{-9}$ & 513 (748) & 1T & 85.2 & 6.8 \\
& BUG & $\varepsilon = 5\times 10^{-7}$ & $1.4\times 10^{-6}$ & $1.5\times 10^{-8}$ & 563 (810) & 1T & 102 & 6.9 \\
\midrule
\multirow{11}{*}{\rotatebox{90}{\textbf{anisotropic}}}
& full tensor & $10^{-10}$ (ref.) & -- & -- & -- & MT & 14138 & 169 \\
& OpenSn & $10^{-10}$, reference & -- & -- & -- & 24R & 1312 & 32.5 \\
& OpenSn & $10^{-2}$ & $1.5\times 10^{-5}$ & $2.7\times 10^{-5}$ & -- & 24R & 716 & 26.0 \\
& OpenSn & $10^{-3}$ & $1.4\times 10^{-5}$ & $2.3\times 10^{-5}$ & -- & 24R & 768 & 26.0 \\
& OpenSn & $10^{-6}$ & $8.4\times 10^{-8}$ & $5.3\times 10^{-8}$ & -- & 24R & 998 & 32.4 \\
& AMEn & $\varepsilon = 10^{-4}$ & $3.4\times 10^{-4}$ & $6.2\times 10^{-5}$ & 227 (282) & 1T & 403 & 7.4 \\
& AMEn & $\varepsilon = 10^{-6}$ & $2.3\times 10^{-6}$ & $5.7\times 10^{-7}$ & 496 (547) & 1T & 1840 & 11.1 \\
& BUG & $\varepsilon = 10^{-4}$ & $5.0\times 10^{-4}$ & $2.2\times 10^{-6}$ & 214 (338) & 1T & 14.3 & 6.2 \\
& BUG & $\varepsilon = 10^{-5}$ & $3.4\times 10^{-5}$ & $3.6\times 10^{-8}$ & 348 (549) & 1T & 39.9 & 6.6 \\
& BUG & $\varepsilon = 10^{-6}$ & $3.3\times 10^{-6}$ & $2.9\times 10^{-8}$ & 514 (767) & 1T & 98.4 & 7.1 \\
& BUG & $\varepsilon = 5\times 10^{-7}$ & $1.8\times 10^{-6}$ & $1.6\times 10^{-8}$ & 566 (836) & 1T & 123 & 7.2 \\
\bottomrule
\end{tabular}
\end{table}

The full-array reference implementation uses substantial working
storage in addition to the angular-flux state itself. We restrict
the reference comparisons to $N=32$; a full state at $N=64$ would
occupy approximately $34$\,GB, but reference-solver memory depends
also on the implementation and its working vectors. No accuracy or
performance conclusion for $N=64$ is drawn here.

\begin{figure}[!htbp]
\centering
\includegraphics[width=\textwidth]{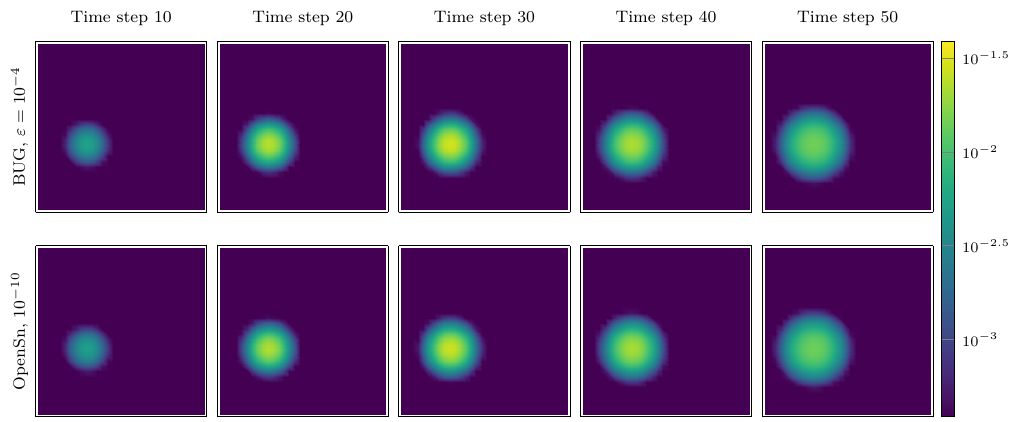}
\caption{Active burst with anisotropic scattering and group speeds,
$N = 32$, $M = 128$, $G = 16$. Neutron density summed over the groups in the plane
$z = 1/2$. Top: BUG at
$\varepsilon = 10^{-4}$. Bottom: OpenSn at the tolerance $10^{-10}$. One
logarithmic color scale for all panels, two decades below the largest
density, shared with Figure~\ref{fig:passive-frames}.}
\label{fig:aniso-frames}
\end{figure}

Figure~\ref{fig:aniso-frames} shows neutron density, the scalar flux
divided groupwise by speed and summed over groups, in the plane
$z=1/2$. BUG and OpenSn show similar localized profiles on the common
colour scale. No group-resolved diagnostic is supplied to attribute
this localization to slow-group dominance. The anisotropic experiment
is active; these results do not establish performance for anisotropic
passive transport or realistic material cross sections.

\FloatBarrier

\section{Conclusions}
\label{sec:conclusions}

We have constructed a hybrid TT/QTT representation of the upwind DG
discretization for three-dimensional transient multigroup neutron
transport. The central contribution is the fixed-matrix upwind
construction: the same two one-dimensional DG matrices serve all
spatial axes, directions and energy groups, with the signed and absolute
direction cosines supplying the angular dependence. This hybrid TT/QTT
representation preserves the Kronecker structure and admits explicit
operator-rank bounds for uniform Cartesian grids and homogeneous
coefficients. These bounds are
independent of mesh size, group count and direction count at fixed
scattering order; they do not bound solution ranks or runtime. The
AMEn-based Crank--Nicolson verification exhibits second-order
convergence. Its finest separable manufactured state takes $13$\,kB
in the hybrid format against $35$\,TB as a full array, a final-state
compression factor of approximately $2.7\times10^9$. Peak process memory
is $2.64$\,GB.

The passive and active burst tests validate compressed states against
full-array discrete references. At comparable scalar-flux errors,
selected TT-BUG runs are approximately $15$--$20$ times faster than the
tested AMEn baseline and use less peak memory. The complete reference
runs use $113$--$169$\,GB, compared with $6.2$--$9.9$\,GB for BUG.
Reference-generation runtimes do not establish speedups over a
full-array method at matched accuracy. OpenSn provides an independent
discretization comparison. For the reported anisotropic active-burst
configurations, TT-BUG computes the angular flux approximately $50$
times faster on one processor than OpenSn does on $24$ processors.
Measured as allocated processors times wall time, the processor-time
cost is approximately $1200$ times lower. This is a comparison of the
reported configurations: the errors against each code's own discrete
reference do not establish equal physical or angular-flux accuracy.

The group correction restores group-integrated algebraic balance of
the compressed step when its correction system is compatible. It does
not enforce cell-local conservation or positivity. The global
step-defect estimate separates errors before final rounding, rounding
changes and balance corrections. It is conditional on an energy bound
and includes the full projected-integration defect for BUG. The
source-onset and saturation modifications of TT-BUG are empirically
tested; their order, convergence and general tolerance-to-error
relation remain unproved. Comparable rank growth in the two integrators
does not establish minimal approximation ranks.

The present tests cover isotropic passive transport and isotropic and
anisotropic synthetic active transport. They do not establish
performance for anisotropic passive problems, realistic heterogeneous
materials or increasing group and direction counts. Those regimes,
matched-accuracy comparisons with optimized transport solvers, and
local conservation through a macro--micro formulation are directions
for further work.

\section*{Data and code availability}
The MATLAB code, scripts and outputs supporting Section~\ref{sec:examples}
are available from the corresponding author on reasonable request.

\appendix
\section{Computing environment}
\label{app:machine}

The manufactured-solution runs of Section~\ref{sec:manufactured} use
MATLAB R2026b and the TT-Toolbox on one computational thread of an
Apple M4 Max processor (MacBook Pro, $36$\,GB RAM). The burst runs use
MATLAB R2025b and the TT-Toolbox on an Apple M3 Ultra processor
(Mac Studio, $32$ CPU cores, $256$\,GB RAM). The TT runs and the isotropic
active full-array reference use one MATLAB computational thread. The
passive and anisotropic active full-array references use multithreaded
MATLAB, denoted MT in the tables; no exact thread-count normalization
is inferred from the machine's core count.

OpenSn is built from its public repository at commit 3cc5157 and runs
on the same Mac Studio with $24$ MPI ranks. Every time in
Tables~\ref{tab:passive} and~\ref{tab:active} is a wall time of a single
run, without multiplication by threads or ranks. These timings describe
the stated execution modes, not measured CPU time or a controlled
parallel-efficiency study. Peak memory is resident process memory,
summed over MPI ranks for OpenSn. Full-array reference-generation
times appear explicitly in the tables and are not matched-error
performance baselines.

\section{Solver diagnostics of the manufactured solution}
\label{app:tables}

Table~\ref{tab:mms-diag} lists the solver diagnostics behind
Table~\ref{tab:mms}, the manufactured solution of
Section~\ref{sec:manufactured}.

\begin{table}[htbp]
\centering
\caption{Solver diagnostics of the runs of Table~\ref{tab:mms}: mean and
largest number of AMEn sweeps per step and number of steps that stopped
at the limit of 30; largest true relative residual of the system of a
step after the solve and after the state is rounded; largest rank inside
AMEn and after rounding; number of steps in which the cap $r_{\max} = 16$
removed more than the tolerance; memory of the three matrices of a step
together (the operator and the two Crank--Nicolson matrices), largest
memory of the state returned by AMEn, and peak resident memory of the
MATLAB process ($1.01$\,GB without any computation). The diagnostics come
from instrumented repetitions of the configurations of
Table~\ref{tab:mms}, not from the timed runs of that table.}
\label{tab:mms-diag}
\vspace{6pt}
\footnotesize
\setlength{\tabcolsep}{2.6pt}
\begin{tabular}{rl rrr cc rr r rrr}
\toprule
& & \multicolumn{3}{c}{sweeps} & \multicolumn{2}{c}{residual} & \multicolumn{2}{c}{rank} & cap & \multicolumn{3}{c}{memory} \\
\cmidrule(lr){3-5}\cmidrule(lr){6-7}\cmidrule(lr){8-9}\cmidrule(lr){11-13}
$\bm{N}$ & \textbf{Form} & \textbf{Mean} & \textbf{Max} & \textbf{At 30} & \textbf{Solve} & \textbf{Rounded} & \textbf{AMEn} & \textbf{State} & \textbf{Acted} & \textbf{Matrices} & \textbf{AMEn state} & \textbf{Process} \\
\midrule
16 & TT & 27.5 & 30 & 98 & $4.1\times 10^{-8}$ & $5.9\times 10^{-6}$ & 76 & 16 & 153 & 14\,MB & 2.1\,MB & 1.52\,GB \\
16 & hybrid & 29.6 & 30 & 149 & $1.5\times 10^{-7}$ & $1.4\times 10^{-5}$ & 76 & 16 & 160 & 13\,MB & 1.2\,MB & 1.36\,GB \\
32 & TT & 16.7 & 20 & 0 & $2.4\times 10^{-8}$ & $1.0\times 10^{-6}$ & 50 & 14 & 0 & 18\,MB & 1.1\,MB & 1.07\,GB \\
32 & hybrid & 19.4 & 22 & 0 & $4.8\times 10^{-8}$ & $1.5\times 10^{-6}$ & 58 & 16 & 311 & 13\,MB & 0.58\,MB & 1.02\,GB \\
64 & TT & 8.2 & 11 & 0 & $1.5\times 10^{-8}$ & $6.6\times 10^{-7}$ & 23 & 4 & 0 & 34\,MB & 0.39\,MB & 1.09\,GB \\
64 & hybrid & 8.0 & 11 & 0 & $8.0\times 10^{-8}$ & $5.2\times 10^{-7}$ & 31 & 12 & 0 & 13\,MB & 0.21\,MB & 1.08\,GB \\
128 & TT & 7.0 & 9 & 0 & $2.0\times 10^{-8}$ & $8.9\times 10^{-8}$ & 17 & 4 & 0 & 96\,MB & 0.39\,MB & 1.32\,GB \\
128 & hybrid & 7.6 & 10 & 0 & $4.2\times 10^{-8}$ & $4.1\times 10^{-7}$ & 23 & 7 & 0 & 13\,MB & 0.11\,MB & 1.30\,GB \\
256 & TT & 4.0 & 6 & 0 & $1.1\times 10^{-8}$ & $1.3\times 10^{-7}$ & 10 & 4 & 0 & 346\,MB & 0.36\,MB & 2.64\,GB \\
256 & hybrid & 6.5 & 9 & 0 & $3.2\times 10^{-8}$ & $4.3\times 10^{-7}$ & 20 & 7 & 0 & 13\,MB & 0.088\,MB & 2.64\,GB \\
\bottomrule
\end{tabular}
\end{table}

\section{The modified rank-adaptive BUG step}
\label{app:bug}

Algorithm~\ref{alg:bug} specifies the implemented step. Write
$p=3d+2$ for the total number of cores; $d$ remains the number of digits
of one spatial axis. The scaled source is $\widetilde{\mathbf q}_s^n$,
consistent with the scaled matrices in Table~\ref{tab:notation}.
Define the auxiliary basis-enrichment predictor
\begin{equation*}
\mathbf Z^n=\operatorname{round}\!\left(
\widetilde{\mathbf B}_s(\widetilde{\bm\Psi}^n+
\widetilde{\mathbf q}_s^n),\varepsilon\right).
\end{equation*}
This predictor is not the physical right-hand side
$\widetilde{\mathbf B}_s\widetilde{\bm\Psi}^n+
\widetilde{\mathbf q}_s^n$; it supplies additional left bases only.
Each local frame $\mathbf F_k$ is an orthonormal embedding, so its
Galerkin matrices are $\mathbf F_k^\top\widetilde{\mathbf K}_s\mathbf F_k$
and $\mathbf F_k^\top\widetilde{\mathbf B}_s\mathbf F_k$.
The source and auxiliary tensors enter the left orthogonalization
with zero coefficient in the final core, preserving their left bases
without changing the input state. All norms in the algorithm are
Frobenius norms of coefficient tensors. Rounding uses tolerance
$\varepsilon$, subject to the rank-cap qualification in Section~\ref{sec:tt}.

\begin{algorithm}[htbp]
\DontPrintSemicolon
\caption{One step $\widetilde{\bm{\Psi}}^n \mapsto \widetilde{\bm{\Psi}}^{n+1}$ of the modified rank-adaptive BUG integrator, in the notation of Table~\ref{tab:notation}.}
\label{alg:bug}
\SetKwFunction{Sweep}{Sweep}\SetKwFunction{Sat}{Saturated}
\textbf{Function} \Sweep{$\bm{\Psi}$, $\mathbf{X}$}: left-orthogonalize $\bm{\Psi} + 0\cdot\widetilde{\mathbf q}_s^n + 0\cdot\mathbf{X}$ (fixed bases $\mathbf{Q}_{<k}$)\;
\Indp for $k = p, p-1, \dots, 2$: express $\bm{\Psi}$ in $\mathbf{Q}_{<k}$, leaf $k$ and the right basis $\mathbf{V}_{>k}$ as $\mathbf{x}_0$; solve the Galerkin
Crank--Nicolson system $(\mathbf F_k^\top\widetilde{\mathbf K}_s\mathbf F_k)\mathbf x_1 = (\mathbf F_k^\top\widetilde{\mathbf B}_s\mathbf F_k)\mathbf x_0 + \mathbf F_k^\top\widetilde{\mathbf q}_s^n$ by GMRES to $10^{-8}$;
$\mathbf{V}_{\ge k} \leftarrow$ orthonormal basis of the rows of $[\mathbf{x}_1; \mathbf{x}_0]$, dropping new directions with singular values below $0.1\,\varepsilon\|\mathbf{x}_1\|$\;
at the root ($k = 1$): Galerkin Crank--Nicolson step on $\mathbb{R}^{n_1}\otimes\operatorname{span}\mathbf{V}_{\ge 2}$; \Return the result with the augmented ranks\;
\Indm
\textbf{Function} \Sat{$\mathbf{Y}$}: true if at some bond below its full rank the augmentation added new directions and $\operatorname{round}(\mathbf{Y},\varepsilon)$ kept all augmented directions\;
\BlankLine
$\mathit{onset} \leftarrow$ a source component switches from zero to nonzero on this step or did so on the preceding step (including reactivation)\;
\eIf{$\mathit{onset}$}{
  $\mathbf{Y} \leftarrow$ \Sweep{$\widetilde{\bm{\Psi}}^n$, $\mathbf Z^n$}; $\mathit{iterate} \leftarrow$ true\;
}{
  $\mathbf{Y} \leftarrow$ \Sweep{$\widetilde{\bm{\Psi}}^n$, $\mathbf{0}$}; $\mathit{iterate} \leftarrow$ \Sat{$\mathbf{Y}$}\;
  \If{$\mathit{iterate}$}{$\mathbf{Y} \leftarrow$ \Sweep{$\widetilde{\bm{\Psi}}^n$, $\mathbf Z^n + \operatorname{round}(\mathbf Y,\varepsilon)$}\;}
}
\If{$\mathit{iterate}$}{
  \Repeat{$\|\mathbf Y\|=0$ \textbf{Or} $\|\mathbf{Y} - \mathbf{Y}_{\mathrm{prev}}\| < \varepsilon\,\|\mathbf{Y}\|$ \textbf{Or} $P_{\max}$ passes have been made}{
    $\mathbf{Y}_{\mathrm{prev}} \leftarrow \operatorname{round}(\mathbf{Y},\varepsilon)$; $\mathbf{Y} \leftarrow$ \Sweep{$\widetilde{\bm{\Psi}}^n$, $\mathbf{Y}_{\mathrm{prev}}$}\;
  }
}
$\widetilde{\bm{\Psi}}^{n+1} \leftarrow$ group balance projection of $\operatorname{round}(\mathbf{Y},\varepsilon)$ (Section~\ref{sec:error-propagation})\;
\end{algorithm}

\paragraph{Implementation notes}
The first implementation, which applied the start-up rule on the
first two steps only, failed our tests with a delayed source onset, and an implementation
without the rule has an error floor of order $10^{-3}$ that does not
decrease with the tolerance. A trigger based on the part of
$\operatorname{round}(\widetilde{\mathbf B}_s(\widetilde{\bm\Psi}^n+\widetilde{\mathbf q}_s^n), \varepsilon)$
outside the frames of $\widetilde{\bm{\Psi}}^n$ was also tried: on ordinary steps
that part is $1.5\times 10^{-2}$ to $10^{-1}$ of the norm, far above any
tolerance used, so the trigger fires on every step, which costs two to
three times as much for slightly smaller errors; it is not used. With at
most 16 sweeps the start-up iteration did not converge on the active
burst with $M = 128$; the limit is 64. On the passive burst the
rule acts on steps 1 and 2 and, for $\varepsilon = 10^{-4}$, on two
further steps triggered by saturation. Source onset is detected from
the prescribed component amplitudes in the endpoint-averaged step
source: a zero-to-nonzero transition triggers the rule, including after
an earlier pulse. The 64-pass limit is a safeguard, not a residual
acceptance condition. If it is reached first, the displayed algorithm
still rounds and projects the last iterate; the change criterion may
then remain unsatisfied. No convergence guarantee is inferred from
that exit. The per-run frequency of such exits is not reported here.

\section{Details of the rank bounds}
\label{app:rank-details}

The first two subsections give details of Proposition~\ref{prop:rank-operator}.
The third explains the hypothesis and scope of Proposition~\ref{prop:rank-solution}.
The final subsection records the Kronecker term count.

\subsection{What the operator bound counts}
\label{app:rank-counts}
\label{rem:rank-operator}
The count in \ref{prop:rank-operator-ii} treats the cores of $\mathbf{C}$ and $\mathbf{D}$ as
unrelated, while in fact they are built from the same three shift
states of~\cite{kazeev2012laplace}, so the bound need not be attained.
In the example of Section~\ref{sec:manufactured}, where the modes are
ordered group, direction, $x$, $y$, $z$, the operator assembled as a
tensor train with one core per index has largest rank 7 (6 for the
two coarsest tolerances); splitting its space cores into digits raises this
to between 7 and 10, and rounding to the tolerance brings it back to 7
on every mesh from $N = 16$ to $N = 256$. These are ranks after a
tolerance-dependent rounding, measured in the implemented order, to which
the counts of the proposition do not apply; the bound for that order is
the value $12$ of \ref{app:rank-order}.
The rank bound concerns the operator contribution; total work
still depends on state ranks, mode sizes and the integration algorithm.

\subsection{Other orders of the indices}
\label{app:rank-order}
\label{rem:rank-order}
The counts \ref{prop:rank-operator-ii} and Figure~\ref{fig:rank-links}
use the order $x$, $y$, $z$, direction, group, in which the rows of an
unfolding inside the $j$-th axis hold identities for every term that does
not act on a preceding axis. The implementation of
Section~\ref{sec:examples} uses the order group, direction, $x$, $y$,
$z$, and then the group and direction factors of every term lie in the
rows, so that this merging is not available. A weaker bound holds in
every order. Across a link between two different indices each of the $q$
products has rank one, so $r_k \le q$ by Lemma~\ref{lem:rank-facts}(a).
Inside the digits of an axis the two streaming products of that axis
contribute $r_\ell(\mathbf{C})$ and $r_\ell(\mathbf{D})$ and every other
product carries the identity on that axis and contributes one, by
Lemma~\ref{lem:rank-facts}(c) and (d), so that
$r_k \le r_\ell(\mathbf{C}) + r_\ell(\mathbf{D}) + q - 2 \le q + 4$.
The same holds for the one-sided matrices $\mathbf{C} + \mathbf{D}$ and
$\mathbf{D} - \mathbf{C}$ that the implementation stores, whose blocks
have the structure \eqref{eq:graphite-CD-blocks} and whose quantized
ranks are therefore at most three by the argument of
\ref{prop:rank-operator-iii}. For the $q = 8$ products of
Section~\ref{sec:manufactured} the bound is $12$, again independent of
$N$, $M$, $G$ and the step.

\subsection{Interpretation of the single-rounding bound}
\label{app:rank-solution-interpretation}
\label{rem:rank-solution}
Proposition~\ref{prop:rank-solution} is a statement about one rounding
of one tensor, not about a computed time iteration. To read it, take for
$\widehat{\bm{\Psi}}$ the state the solver delivers at some step, before
rounding, and for $\bm{\Psi}_\ast$
the coefficients of the projection \eqref{eq:graphite-proj} of the exact
flux at that time, or of any approximation of the exact flux with known
ranks. Then $\eta$ is the distance between the computed state and
$\bm{\Psi}_\ast$. It contains the discretization error together with
the errors of the inexact solves, of the compressed operators and
sources, and of all earlier roundings. The hypothesis
\eqref{eq:rank-solution-hyp} requires the relative tolerance $\varepsilon$ to
be at least $\sqrt{p-1}$ times $\eta/\|\widehat{\bm{\Psi}}\|$. In the hybrid
form $p$ grows like $3\log_2 N$, so this factor grows slowly with the
mesh. When the hypothesis holds, the rounded state has ranks at most
$r_k^\ast$ and lies within $\eta + \varepsilon\|\widehat{\bm{\Psi}}\|$ of
$\bm{\Psi}_\ast$, on every mesh. The ranks $r_k^\ast$ themselves must
be known or bounded separately. For the exact flux \eqref{eq:mms-exact}
of Section~\ref{sec:manufactured} they are one at every link between
two different indices and at most two at the links inside the digits of
an axis. The reason is the addition theorem: on a block of $2^{\ell-1}$
cells, the coefficients of $\sin(\pi x)$ are those of a function in the
span of $\sin(\pi s)$ and $\cos(\pi s)$, with $s$ the local coordinate
of the block, so the $\ell$-th unfolding has at most two independent
columns.

The runs of Table~\ref{tab:mms} use $\varepsilon = 10^{-6}$, which is below
their relative discretization error. The hypothesis therefore fails for
these runs, and the proposition does not explain the ranks they carry,
up to seven on the fine meshes. Nor does it bound the ranks inside AMEn
or the ranks of the state along a time history.
Section~\ref{sec:error-propagation} bounds the propagation of the errors
under explicit assumptions. Turning that estimate into a bound on the
ranks of a computed time history would require, in addition, suitable
low-rank comparison tensors at every step and a verification of the
hypothesis of Proposition~\ref{prop:rank-solution} at every step. We do
not establish these conditions. Finally, the mass matrix of the scheme
is a Kronecker product of diagonal matrices. Its square root therefore
acts core by core and leaves every rank unchanged, so the proposition
holds verbatim in the mean-square norm of the scheme when the rounding
is done in that norm.

\subsection{Kronecker term count}
\label{rem:term-count}
The physics matrix is the sum of $q$ Kronecker products of five factors,
with one factor for each of the five indices. The operators of
\eqref{eq:semidiscrete} give the products. The time term $\mathbf{V}^{-1}$
and the removal by collisions $\bm{\Sigma}_t$ share one product, the
streaming operator $\mathbf{L}$ gives six, and the scattering operator
$\mathbf{S}$ gives $\ell_{\max}+1$, where $\ell_{\max}$ is the Legendre
order of its expansion \eqref{eq:rank-K-scatter}. We explore scattering
only, so the term count is $q = 8 + \ell_{\max}$, and $q = 8$ for the
isotropic cases considered here. In principle one could include
fission, whose operator \eqref{eq:rank-K-fission} adds one product, so
that $q = 9 + \ell_{\max}$. The matrix $\mathbf{B}$ has similar structure. Proposition~\ref{prop:rank-operator} bounds the ranks of
$\mathbf{K}$ and $\mathbf{B}$ across every link by $q$ and by the
quantized ranks of $\mathbf{C}$ and $\mathbf{D}$, independently of $N$,
$M$, $G$ and the step.

\section{Energy estimate for anisotropic discrete scattering}
\label{app:anisotropic-energy}

The energy estimate can be extended without assuming that the discrete
angular matrices are mutually orthogonal projectors. In the
direction--group ordering, define
\begin{align*}
\mathbf C_\Omega
 &=\mathbf I_M\otimes\operatorname{diag}(\sigma^1,\ldots,\sigma^G)
   -\sum_{\ell=0}^{\ell_{\max}}\mathbf P_\ell\otimes\bm\Sigma_{s,\ell}
   -\mathbf P_0\otimes\bm\chi(\nu\bm\sigma_f)^\top,\\
\mathbf W_\Omega&=\operatorname{diag}(\mathbf w)\otimes\mathbf I_G,
\qquad
\mathbf H_\Omega=\operatorname{diag}(\mathbf w)\otimes\mathbf V_G^{-1}.
\end{align*}
Positive quadrature weights and positive speeds make $\mathbf H_\Omega$
positive definite. Set
\begin{equation}
\gamma_\Omega=\max\!\left\{0,-\lambda_{\min}\!\left(
\mathbf H_\Omega^{-1/2}
\operatorname{sym}(\mathbf W_\Omega\mathbf C_\Omega)
\mathbf H_\Omega^{-1/2}\right)\right\}.
\label{eq:anisotropic-gamma}
\end{equation}
By the Rayleigh quotient,
$\mathbf y^\top\mathbf W_\Omega\mathbf C_\Omega\mathbf y
\ge-\gamma_\Omega\mathbf y^\top\mathbf H_\Omega\mathbf y$ for every
$\mathbf y$. Tensoring this inequality with the positive definite
spatial mass matrix and adding the nonnegative streaming contribution
of Lemma~\ref{lem:energy}(i) yields
\begin{equation*}
\mathbf x^\top\mathbf W\mathcal A\mathbf x
\ge-\gamma_\Omega\|\mathbf x\|_E^2.
\end{equation*}
Consequently, Lemma~\ref{lem:onestep} and Proposition~\ref{prop:global}
hold with $\gamma_\Omega$ in place of $\gamma$, provided
$\theta\gamma_\Omega\Delta t<1$. This is a bound for the actual
discrete angular--group collision matrix, not a consequence of
non-multiplication. We do not evaluate its numerical size or all the
step-defect terms for the anisotropic benchmark, and it supplies no
independent convergence result for the modified BUG integrator.

\section*{Acknowledgments}
The authors appreciate the support of Eagle Nuclear Energy Corporation and the University of New Mexico Center for Advanced Research Computing, and 
Sandia National Laboratories, which is a multimission laboratory managed and
operated by National Technology \& Engineering Solutions of Sandia, LLC,
a wholly owned subsidiary of Honeywell International Inc., for the
U.S.\ Department of Energy's National Nuclear Security Administration
under contract DE-NA0003525. This paper describes objective technical
results and analysis. Any subjective views or opinions that might be
expressed in the paper do not necessarily represent the views of the
U.S.\ Department of Energy or the United States Government.

\bibliographystyle{plain}
\bibliography{ttrain}

\end{document}